\documentclass[a4paper,10pt]{article}
\usepackage[utf8]{inputenc}
\usepackage{amssymb}
\usepackage{amsthm}
\usepackage{tikz}
\usepackage{tikz-cd}
\usepackage{bm} 
\usepackage{color}
\usepackage{mathtools}
\usepackage[titletoc,page]{appendix}
\usepackage[export]{adjustbox}
\usepackage{enumitem}

\usepackage[mathscr]{euscript} 
\DeclareFontFamily{OT1}{pzc}{}
\DeclareFontShape{OT1}{pzc}{m}{it}{<-> s * [1.10] pzcmi7t}{}
\DeclareMathAlphabet{\mathpzc}{OT1}{pzc}{m}{it} 
     
\usetikzlibrary{calc}
\usetikzlibrary{matrix,arrows,decorations}
\newtheorem{theorem}{Theorem}[subsection]
\newtheorem{lemma}[theorem]{Lemma}
\newtheorem{remark}[theorem]{Remark}
\newtheorem{definition}[theorem]{Definition}
\newtheorem{construction}[theorem]{Construction}
\newtheorem{proposition}[theorem]{Proposition}

\newtheorem{example/theorem}[theorem]{Example/Theorem}
\newtheorem{construction/theorem}[theorem]{Construction/Theorem}

\newtheorem{convention/notation}[theorem]{Convention/Notation}
\newtheorem{corollary}[theorem]{Corollary}

\newtheorem{construction/definition}[theorem]{Construction/Definition} 
\newtheorem{proposition/definition}[theorem]{Proposition/Definition} 
\newtheorem{example/construction}[theorem]{Example/Construction}
\newtheorem{theorem/facts}[theorem]{Theorem/Facts}

\newcommand{\ConSum}[1]{\underset{#1}{\#}}

\usepackage{hyperref}

\title{Topological $K$-theory with coefficients and the $e$-invariant}
\author{Yi-Sheng Wang\footnote{Email:yswangbl@gmail.com}}
 
\begin{document}
\maketitle  
\section*{Abstract}
We compare the invariants of flat vector bundles defined by Atiyah et al. and Jones et al. and prove that, up to weak homotopy, they induce the same map, denoted by $e$, from the $0$-connective algebraic $K$-theory space of the complex numbers to the homotopy fiber of the Chern character. We examine homotopy properties of this map and its relation with other known invariants. In addition, using the formula for $\tilde{\xi}$-invariants of lens spaces derived from Donnelly's fixed point theorem and the $4$-dimensional cobordisms constructed via relative Kirby diagrams, we recover the formula for the real part of $e$-invariants of Seifert homology spheres given by Jones and Westbury, up to sign.
 
We conjecture that this geometrically defined map $e$ can be represented by an infinite loop map. The results in its companion paper \cite{Wang2} give strong evidence for this conjecture.    

\tableofcontents

\section{Introduction} 

Complex flat vector bundle theory is intimately linked to algebraic $K$-theory of the complex numbers via the (topological) group completion theorem and also plays an important role in differential topology through Chern-Simons theory. Invariants, such as Kamber-Tondeur classes, constructed via flat vector bundle theory are often useful in the study of both algebraic $K$-theory and differential topology. 

The aim of this paper and its companion \cite{Wang2} is to study maps from the $0$-connective algebraic $K$-theory space of the complex numbers, denoted by $K_{a}\mathbb{C}$, to the homotopy fiber of the Chern character, denoted by $F_{t,\mathbb{C/Z}}$. In the present paper, we examine carefully the properties of the map $e:K_{a}\mathbb{C}\rightarrow F_{t,\mathbb{C/Z}}$ induced from the construction of the topological index in \cite[p.87]{APS3}, and in the subsequent paper \cite{Wang2}, we provide evidence in favor of the conjecture that the map $e$ can be lifted to a map in the stable homotopy category. One important implication of this conjecture is a refined $\operatorname{BL}$ index theorem.

\subsection*{Main results}
In their study of the index theorem for flat vector bundles \cite[p.87]{APS3}, Atiyah and Patodi and Singer construct a geometric model for $\mathbb{C/Z}$-$K$-theory, denoted by $\mathcal{G}_{\operatorname{APS}}(M)$, and define a homomorphism (of semigroups) from the semigroup of stable flat vector bundles to the abelian group $\mathcal{G}_{\operatorname{APS}}(M)$:
\begin{equation}
\bar{e}_{\operatorname{APS}}:[M,\operatorname{BGL}(\mathbb{C}^{\delta})]\rightarrow \mathcal{G}_{\operatorname{APS}}(M),
\end{equation} 
for every pointed compact smooth manifold $M$. The abelian group $\mathcal{G}_{\operatorname{APS}}(M)$ is generated by the geometric cocycles $\{([\omega],(V,\phi))\}$ subject to certain relations, where $\omega$ is an odd closed differential form and $mV\xrightarrow{\phi}\epsilon^{m\operatorname{dim}V}$ is an isomorphism between the direct sum of $m$ copies of the vector bundle $V$ and the trivial bundle over $M$.

In \cite{JW}, another geometric model $\mathcal{G}_{\operatorname{JW}}(M)$ is considered and used to define a similar homomorphism,
\begin{equation}\label{JW}
\bar{e}_{\operatorname{JW}}:[M,\operatorname{BGL}(\mathbb{C}^{\delta})]\rightarrow \mathcal{G}_{\operatorname{JW}}(M),
\end{equation} 
where the abelian group $\mathcal{G}_{\operatorname{JW}}(M)$ is generated by the $5$-tuples $\{(V,\nabla_{v},W,\nabla_{w},\omega)\}$ subject to certain relations, where $\nabla_{v}$ (resp. $\nabla_{w}$) is a connection on the vector bundle $V$ (resp. $W$), and $\omega$ is an odd differential form whose differential measures the difference between the Chern characters of $\nabla_{v}$ and $\nabla_{w}$\footnote{This model and the construction of the homomorphism are credited to Karoubi \cite[Chapter $7$]{Ka1} by Jones and Westbury.}. 

Though both models satisfy the following exact sequence
\[\tilde{K}^{-1}(M)\xrightarrow{\operatorname{ch}} \tilde{H}^{odd}(M,\mathbb{C})\rightarrow \mathcal{G}_{\operatorname{APS/JW}}(M)\rightarrow \tilde{K}(M)\xrightarrow{\operatorname{ch}} \tilde{H}^{ev}(M,\mathbb{C}),\] 
it was not clear if they are isomorphic to each other and equivalent to the (topological) definition of complex topological $K$-theory with coefficients in $\mathbb{C/Z}$, where $\tilde{K}(M)$ is the reduced topological $K$-theory of $M$, $\operatorname{ch}$ is the Chern character, and 
\[\tilde{H}^{ev/odd}(M,\mathbb{C}):=[M,\prod_{\mathclap{i \text{ even/odd }}}K(\mathbb{C},i)].\] By $K(A,i)$, we understand the Eilenberg-Maclane space of the abelian group $A$ in degree $i$.
 
\begin{theorem}\label{Intro:Thm1} 
Let $[M,F_{t,\mathbb{C/Z}}]$ be the abelian group of homotopy classes of pointed maps from $M$ to $F_{t,\mathbb{C/Z}}$ and $[-,-]_{\operatorname{Ho}(\mathcal{P})}$ denote the abelian group of homotopy classes of maps in the category of prespectra $\mathcal{P}$. Then there are canonical isomorphisms
\[[\Sigma^{\infty}M,\Sigma^{-1}\mathbf{K}_{t}\wedge\mathbf{M}\mathbb{C/Z}
]_{\operatorname{Ho}(\mathcal{P})}\simeq [M,F_{t,\mathbb{C/Z}}]\xleftarrow{\sim} \mathcal{G}_{APS}(M)\overset{\sim}{\rightleftarrows}\mathcal{G}_{JW}(M),\]
for every pointed compact smooth manifold $M$, where $\mathbf{K}_{t}$ is the complex topological $K$-theory prespectrum and $\mathbf{M}\mathbb{C/Z}$ the Moore prespectrum.  
\end{theorem}
Utilizing the last isomorphism in Theorem \ref{Intro:Thm1}, we show that homomorphisms $\bar{e}_{\operatorname{JW}}$ and $\bar{e}_{\operatorname{APS}}$ induce the same map, up to weak homotopy, from $K_{a}\mathbb{C}$ to $F_{t,\mathbb{C/Z}}$\footnote{The name $e$ is due to Jones and Westbury \cite{JW}. We justify this name by proving it generalizes the Adams $e$-invariant in the subsequent paper \cite{Wang2}.}:  
\[e:K_{a}\mathbb{C}\rightarrow F_{t,\mathbb{C/Z}}.\]

By the construction of the map $e$ and Atiyah, Patodi and Singer's index theorem for flat vector bundles \cite[Theorem $5.3$]{APS3}, it is not difficult to see how $e$ is related to the $\tilde{\xi}$-invariant. However, we need yet another result in order to make effective use of the $\tilde{\xi}$-invariant to compute the $e$-invariant---the induced homomorphism of $e$:
\begin{theorem}\label{Intro:Thm2}
Let  
\begin{align*}
\rho_{1}:\pi_{1}(\Sigma_{1}^{n})&\rightarrow \operatorname{GL}_{N}(\mathbb{C}),\\
\rho_{2}:\pi_{1}(\Sigma_{2}^{n})&\rightarrow \operatorname{GL}_{N}(\mathbb{C}),
\end{align*}
be representations of homology spheres $\Sigma_{1}^{n}$ and $\Sigma_{2}^{n}$. Suppose that, via the plus map $\Sigma^{n}_{i}\rightarrow S^{n}$, $i=1$ or $2$ (see homomorphism \eqref{Eq:flatonhomosphere}), $(\Sigma_{1}^{n},\rho_{1})$ and $(\Sigma_{2}^{n},\rho_{2})$ yield the same element in $[S^{n},K_{a}\mathbb{C}]$.  
Then $\tilde{\xi}(\Sigma_{1}^{n},\rho_{1})=\tilde{\xi}(\Sigma_{2}^{n},\rho_{2})$.  
\end{theorem}
In fact, we conjecture that given any two plus maps between compact smooth manifolds $M_{1}\rightarrow M$ and $M_{2}\rightarrow M$ and two representations $\rho_{1}:\pi_{1}(M_{1})\rightarrow \operatorname{GL}_{N}(\mathbb{C})$ and $\rho_{2}:\pi_{1}(M_{2})\rightarrow \operatorname{GL}_{N}(\mathbb{C})$, if $(M_{1},\rho_{1})$ and $(M_{2},\rho_{2})$ induce the same element in $[M,K_{a}\mathbb{C}]$ via these plus maps, then $\tilde{\xi}(M_{1},\rho_{1})=\tilde{\xi}(M_{2},\rho_{2})$.

Next, we employ the notion of locally unipotent bundles over a manifold to show that the $e$-invariant is highly non-trivial. 
\begin{theorem}\label{Intro:Thm3}\footnote{This has been observed by Jones and Westbury, yet the argument presented there is not entirely clear to us. For more details, we refer to the footnote on p.\pageref{CommDiaforeTor}.} 
When restricted to the torsion part, the map $e$ induces an isomorphism:
\[e_{\ast}\vert_{\operatorname{Tor}}:\operatorname{Tor}[L,K_{a}\mathbb{C}]\xrightarrow{\sim} \operatorname{Tor}[L,F_{t,\mathbb{C/Z}}],\]
for every finite $CW$-complex $L$. 
\end{theorem}

The next theorem is about spin (resp. $\operatorname{spin}^{c}$) structures on certain $4$-dimensional cobordisms. The result is important for computing $e$-invariants of Seifert homology spheres, but it is of interest in its own right. 
\begin{theorem}\label{Intro:Thm4}
Given pairwise coprime integers $a_{1},a_{2},...,a_{n}$ and integers $b,b_{1},b_{2},...,b_{n}$ such that 
\[a_{1}...a_{n}(-b+\sum_{i=1}^{n}\frac{b_{i}}{a_{i}})=1,\]
then the $4$-dimensional cobordism represented by the following relative Kirby diagram
\begin{figure}[h]
\caption{}
\label{LD0}
\centering
\includegraphics[width=0.65\textwidth]{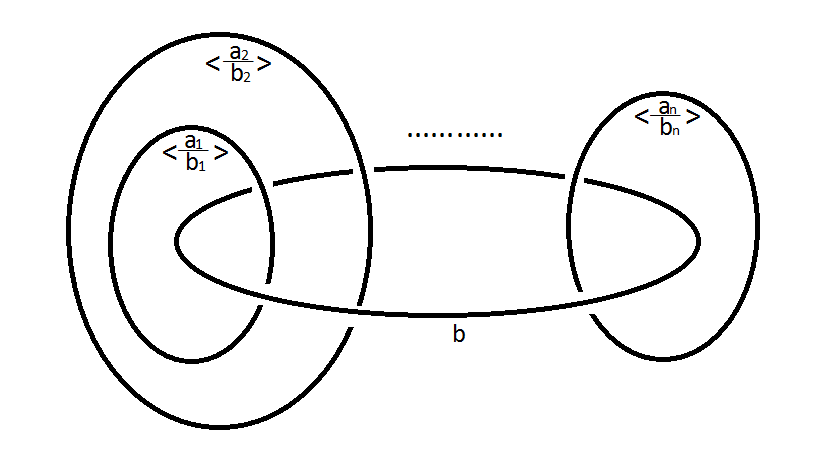}
\end{figure}

\noindent
is spin when the product $a_{1}...a_{n}$ is even. When it is odd, there is a $\operatorname{spin}^{c}$ structure that restricts to the canonical $\operatorname{spin}^{c}$ structure on the boundary. 
\end{theorem} 
 
\subsection*{Outline of the paper}
Section $2$ discusses the topological approach to cohomology theory with coefficients and reviews some basic notions and results in rational homotopy theory. In Section $3$, we recall some geometric models for topological $K$-theory with coefficients in the literature and prove that they are equivalent to the topological definitions given in Section $2$. Based on the results in Sections $2$ and $3$, we explain how $\bar{e}_{\operatorname{APS}}$ gives rise to the map $e$ in Section $4$. In the same section, the relation of the map $e$ with other invariants and its homotopy properties are also investigated. Large part of Section $5$ is devoted to the study of spin and $\operatorname{spin}^{c}$ structures on the $4$-cobordism represented by relative Kirby diagram \eqref{LD0}, which we use to compute $e$-invariants of Seifert homology spheres and recover, up to sign, the formula given in \cite{JW}\footnote{We refer to Remark \ref{CobordisminJW} as well as Remark \ref{XiofLensinJW} for the comparison between their approaches and ours.}. 

The primary model for the stable homotopy category in this paper is the category of prespectra $\mathcal{P}$. We sometimes, however, need some constructions in the Adams category $\mathcal{A}$. The appendix is meant to review some preliminaries on prespectra and give a detail comparison between $\mathcal{P}$ and $\mathcal{A}$. 

\subsection*{Relation with other papers}
Maps from algebraic $K$-theory to topological $K$-theory with coefficients in $\mathbb{C/Z}$ have also been studied in \cite{BNV}, \cite{Bu1}, \cite{Bu2} and \cite{ASS}. In \cite[Example $6.9$]{BNV} \cite{Bu1} and \cite{Bu2}, two such maps are constructed in terms of differential cohomology theories via sheaves of spectra. However, the relation between these three maps, the two maps and the map $e$, is not entirely clear \cite[p.5]{Wang2}. In fact, one might need a space level comparison in order to unveil the connection between them \cite[Remark 12.15]{Bu2}. In \cite{ASS}, on the other hand, they use von Neumann algebras to define a generalized topological $K$-theory with coefficients and associate each unitary flat connection with a class in their generalized $K$-theory with coefficients in $\mathbb{R/Z}$. Furthermore, they show that, in the case of commutative $C^{\ast}$-algebras, if considering only unitary flat vector bundles, one can identify $\bar{e}_{\operatorname{APS}}$ with this assignment \cite[Corollary $5.5$]{ASS}. 

This paper is greatly inspired by the work and the idea of Jones and Westbury in \cite{JW}. In particular, part of Sections $4$ and $5$ bears close resemblance to Sections $2$, $3$ and $4$ in \cite{JW}. Our approaches are different from theirs, however. The origin of the geometric model $\mathcal{G}_{\operatorname{JW}}(M)$ and homomorphism \eqref{JW} can be dated back to \cite[Chapter $7$]{Ka1} (see also \cite{Ka2} and \cite{Ka3}), where topological $K$-theory with coefficients in $\mathbb{C/Z}$ is called multiplicative $K$-theory.

\subsection*{Notation and conventions}
 
The convenient category of topological spaces used here is the symmetric monoidal category of pointed $k$-spaces $\mathpzc{Top}_{\ast}$ endowed with the Quillen model structure (see \cite[p.58-60]{Ho}). Given two topological spaces $X,Y$, $[X,Y]$ stands for the set of homotopy classes of maps from $X$ to $Y$ in $\mathpzc{Top}_{\ast}$. For the sake of convenience, all topological spaces considered here are assumed to have the homotopy type of $CW$-complexes. To avoid confusion, the bold letters are reserved for prespectra and maps of prespectra.  

By the $\tilde{\xi}$-invariant of a spin (resp. $\operatorname{spin}^{c}$) manifold, we understand the $\tilde{\xi}$-invariant associated to the Dirac operator induced from the spin (resp. $\operatorname{spin}^{c}$) structure.

\section{Cohomology theory with coefficients} 
In this section, we shall review some basic notions and constructions in stable homotopy theory. Especially, the definition of generalized cohomology theories with coefficients in an abelian group $A$ is given in the first part of the section. Our purpose here is to reinterpret Suslin's result \cite[Corollary $4.6$]{Su} in terms of stable homotopy theory. The second part of the section discusses some basic properties of rationalization of topological spaces and $CW$-prespectra. Throughout this paper, we work primarily in the category of prespectra $\mathcal{P}$. We have chosen $\mathcal{P}$ to be the model for the stable homotopy category because it is large enough to have both homotopy fiber and cofiber constructions and also because it admits a model structure and a Quillen adjunction 
\[\Omega^{\infty}:\mathcal{P}\leftrightarrows \mathpzc{Top}_{\ast}:\Sigma^{\infty}.\] 
However, since some constructions work only for or better with the $CW$-prespectra, we occasionally need to switch from $\mathcal{P}$ to the Adams category $\mathcal{A}$. Nevertheless, in most cases, there is not much difference between them. In particular, given $\mathbf{E}$, $\mathbf{F}$ two $\Omega$-$CW$-prespectra, there is a canonical isomorphism
\[[\mathbf{E},\mathbf{F}]_{\operatorname{Ho}(\mathcal{P})}\xrightarrow{\sim} [\mathbf{E},\mathbf{F}]_{\operatorname{Ho}(\mathcal{A})}.\]
For more details about prespectra and the model structure on $\mathcal{P}$, we refer to \ref{CPre}, whereas a detailed comparison between $\mathcal{A}$ and $\mathcal{P}$ can be found in \ref{AandPcomparison}.   

\subsection{Moore spaces and prespectra}\label{MSandMP}
Recall that a Moore space $M(G,n)$ of an abelian group $G$ in degree $n\in\mathbb{N}$ is a topological space in $\mathpzc{Top}_{\ast}$ equipped with an isomorphism $\phi: \tilde{H}_{n}(M(G,n),\mathbb{Z})\simeq G$ and having $\tilde{H}_{\ast}(M(G,n),\mathbb{Z})=0$ when $\ast\neq n$. If $n>1$, we further demand $M(G,n)$ to be simply connected. Notice that the isomorphism $\phi$ is part of the data, though it is usually dropped from the notation. Given an abelian group and $n\in\mathbb{N}$, there is a standard model for $M(G,n)$: 
Consider a free resolution of $G$, 
\[0\rightarrow R\xrightarrow{i} F\xrightarrow{q} G,\]
where $F$ and $R$ are the free abelian groups with generators $\{a_{\alpha}\}_{\alpha\in A}$ and $\{b_{\beta}\}_{\beta\in B}$, respectively. Assign a $n$-sphere $S^{n}$ to every generator in $\{a_{\alpha}\}_{\alpha\in A}$ and $\{b_{\beta}\}_{\beta\in B}$ and construct a map $f$
\[\underset{\beta\in B}{\vee}  S_{\beta}^{n}\xrightarrow{f} \underset{\alpha\in A}{\vee} S_{\alpha}^{n}\]
such that it realizes the homomorphism $i$. Then the mapping cone of $f$ gives the model for $M(G,n)$. The following lemma can be easily deduced from this construction.

\begin{lemma}\phantomsection\label{MapfromMoorespace}
For $n>1$, given any homomorphism $\phi:G\rightarrow \pi_{n}(X)$, one can find a continuous map 
\[f:M(G,n)\rightarrow X\] 
such that $f_{\ast}=\phi$.
\end{lemma}
This lemma implies the uniqueness of the homotopy type of $M(G,n)$, for $n\geq 2$. On the contrary, when $n=1$, the homotopy type of $M(G,1)$ is not unique. A fixed model is therefore needed. From now on, we let the cofiber of the degree $m$ map 
\[\varphi_{m}:S^{k}\rightarrow S^{k}\]
be our model for $M(\mathbb{Z}/m,k)$ for $k\geq 1$.

Moore spaces can be used to define homotopy groups with finite coefficients (see \cite[Chap.4]{We1} and \cite[Introduction]{Su}). 

\begin{definition}\phantomsection\label{Hogroupwithfcoefspaces}
Given a topological space $X$, its $n$-th homotopy group with coefficients in $\mathbb{Z}/m$ is given by  
\[\pi_{n}(X)_{\mathbb{Z}/m}:=[M(\mathbb{Z}/m,n-1),X], \text{ for } n\geq 3. \]
If $X$ is a group-like $H$-space, the set of homotopy classes \[\pi_{2}(X)_{\mathbb{Z}/m}:=[M(\mathbb{Z}/m,1),X] \] is endowed with a natural abelian group structure, and in this case, we define $\pi_{1}(X)_{\mathbb{Z}/m}$ to be the quotient of $\pi_{1}(X)\xrightarrow{m} \pi_{1}(X)$. When $X$ is a connected $H$-space, these homotopy groups with coefficients in $\mathbb{Z}/m$ fit into the following long exact sequence of abelian groups
\begin{multline*}
...\rightarrow \pi_{n}(X)\xrightarrow{m} \pi_{n}(X)\rightarrow \pi_{n}(X)_{\mathbb{Z}/m}\rightarrow \pi_{n-1}(X)\\
\rightarrow...\rightarrow \pi_{1}(X)\xrightarrow{m} \pi_{1}(X)\rightarrow \pi_{1}(X)_{\mathbb{Z}/m}\rightarrow 0.
\end{multline*}
\end{definition}

\textbf{Moore prespectra:} 
A Moore $CW$-prespectrum $\mathbf{M}G$ of an abelian group $G$ is a $(-1)$-connective $CW$-prespectrum equipped with an isomorphism $\phi:\pi_{0}(\mathbf{H}\mathbb{Z}\wedge\mathbf{M}G)\simeq G$ and having $\pi_{\ast}(\mathbf{H}\mathbb{Z}\wedge\mathbf{M}G)=0$, for $\ast\neq 0$. Via the Hurewicz isomorphism and an analogue of Lemma \ref{MapfromMoorespace}, we see the homotopy type of a Moore $CW$-prespectrum of an abelian group is uniquely determined. Similar to Moore spaces, there is a standard model for a Moore $CW$-prespectrum of a given abelian group. In fact, the construction is basically the same as that of the standard model of a Moore space. One only needs to replace $n$-spheres by sphere prespectra. For example, $\Sigma^{\infty}M(G,k)[-k]$ is a model of the Moore $CW$-prespectrum $\mathbf{M}G$, for every $k\geq 1$. Recall that, given a prespectrum $\mathbf{E}$, $\mathbf{E}[k]$ is the prespectrum whose $n$-th component is $E_{k+n}$, for $n\geq 0$ (p.\pageref{StabilityofP}).

In view of this construction and Lemma \ref{MapfromMoorespace}, given a homomorphism of abelian groups $\phi:H\rightarrow G$, there always exists a map of $CW$-prespectra
\[\bm{\phi}:\mathbf{M}H\rightarrow \mathbf{M}G\]
such that 
\[\phi=\pi_{0}(id\wedge \bm{\phi}):\pi_{0}(\mathbf{H}\mathbb{Z}\wedge \mathbf{M}H)\rightarrow \pi_{0}(\mathbf{H}\mathbb{Z}\wedge \mathbf{M}G).\]
The map $\bm{\phi}$ is however not unique in $\operatorname{Ho}(\mathcal{P})$ as the standard construction of Moore $CW$-prespectra does not constitute a functor from the category of abelian groups to the stable homotopy category. Like Moore spaces, Moore $CW$-prepectra are useful for defining homotopy groups with finite coefficients. 
\begin{definition}\label{Cohomowithcoeff}
Let $\mathbf{F}$ be a $CW$-prespectrum and call the associated cohomology theory $F$-theory. Then $F$-theory with coefficients in $A$ is defined to be the cohomology theory represented $\mathbf{F}\wedge \mathbf{M}A$.
\end{definition}
In this paper, we are only interested in the case where $A=\mathbb{Z}/m$, $\mathbb{Q/Z}$, $\mathbb{R/Z}$, $\mathbb{C/Z}$, $\mathbb{Q}$, $\mathbb{R}$ and $\mathbb{C}$.
 
Now, in view of this definition, given a $CW$-prespectrum $\mathbf{E}$, the $n$-the homotopy group of $\mathbf{E}$ with coefficients in $\mathbb{Z}/m$---in this case, $\mathbf{F}$ is the sphere prespectrum $\mathbf{S}$---is given by 
\[ [\mathbf{S}[n],\mathbf{E}\wedge\mathbf{M}\mathbb{Z}/m]_{\operatorname{Ho}(\mathcal{P})}.\] 
To see its relation with the definition of homotopy groups with finite coefficients for pointed topological spaces (Definition \ref{Hogroupwithfcoefspaces}), we recall the following lemma which is a corollary of the lemmas in \cite[14.31-33]{Sw}.

\begin{lemma}
The $CW$-prespectrum $\Sigma^{\infty}M(\mathbb{Z}/m,n)[-n-1]$ is the $S$-dual of the Moore $CW$-prespectrum $\mathbf{M}\mathbb{Z}/m=\Sigma^{\infty}M(\mathbb{Z}/m,n)[-n]$, which we denoted by $\mathbf{M}^{\ast}\mathbb{Z}/m$. 
\end{lemma} 
\begin{proof}
Use the cofiber sequence 
\[\mathbf{S}\xrightarrow{m} \mathbf{S}\rightarrow \Sigma^{\infty}M(\mathbb{Z}/m,n)[-n] \rightarrow \mathbf{S}[1]\rightarrow \mathbf{S}[1]\]
and Corollary $14.33$ in \cite{Sw}.
\end{proof}
As an immediate consequence of this lemma and the $S$-duality, we have the following isomorphisms:
\[[\mathbf{S}[n],\mathbf{E}\wedge \mathbf{M}\mathbb{Z}/m]_{\operatorname{Ho}(\mathcal{A})}\simeq [\mathbf{M}^{\ast}\mathbb{Z}/m[n],\mathbf{E}]_{\operatorname{Ho}(\mathcal{A})}=[\Sigma^{\infty}M(\mathbb{Z}/m,k)[n-k-1], \mathbf{E}]_{\operatorname{Ho}(\mathcal{A})}.\]

When $\mathbf{E}$ is a $0$-connective $\Omega$-$CW$-prespectrum and $n\geq 2$, there are natural isomorphisms 
\begin{multline*}
[\Sigma^{\infty}M(\mathbb{Z}/m,k)[n-k-1], \mathbf{E}]_{\operatorname{Ho}(\mathcal{A})}\\ \xleftarrow{\sim} [\Sigma^{\infty}M(\mathbb{Z}/m,k)[n-k-1], \mathbf{E}]_{\operatorname{Ho}(\mathcal{P})}\simeq [M(\mathbb{Z}/m,n-1),E_{0}], 
\end{multline*}
where the second isomorphism is induced by the adjunction
\[\Sigma^{\infty}:\operatorname{Ho}(\mathpzc{Top}_{\ast})\leftrightarrows \operatorname{Ho}(\mathcal{P}):\Omega^{\infty}\]
(see \ref{AandPcomparison} and Proposition \ref{QadjunctionPresandspace}).
For $n=1$, we have the following isomorphisms 
\begin{multline*}
[\Sigma^{\infty}M(\mathbb{Z}/m,k)[n-k-1], \mathbf{E}]_{\operatorname{Ho}(\mathcal{P})}\simeq [M(\mathbb{Z}/m,1),E_{1}]\\  
\simeq \operatorname{Coker}([S^{2},E_{1}]\xrightarrow{m}[S^{2},E_{1}])\\
\simeq \operatorname{Coker}([S^{1},E_{0}]\xrightarrow{m}[S^{1},E_{0}])\simeq\pi_{1}(E_{0})_{\mathbb{Z}/m}. 
\end{multline*}
As a result, we have proved the following corollary:
\begin{corollary}\label{Cor:Towdefforhgpwithcoeff}
Let $n\geq 1$ and $\mathbf{E}$ be a $0$-connective $\Omega$-$CW$-prespectrum. Then the notion of homotopy groups with finite coefficients given by Moore $CW$-prespectra coincides with Definition \ref{Hogroupwithfcoefspaces} in the following sense:
\[ [\mathbf{S}[n],\mathbf{E}\wedge\mathbf{M}\mathbb{Z}/m]_{\operatorname{Ho}(\mathcal{P})}=\pi_{n}(E_{0})_{\mathbb{Z}/m}. \]
\end{corollary}

With this corollary, Suslin's theorem \cite[Corollary.4.6]{Su} can be reformulated as follows:
\begin{corollary}\phantomsection\label{Suslinthm1}
The canonical map of $CW$-prespectra 
\begin{align*}
\mathbf{K}_{a}\wedge \mathbf{M}\mathbb{Q/Z}&\rightarrow \mathbf{K}_{t}\wedge \mathbf{M}\mathbb{Q/Z};\\
\mathbf{K}_{a}\wedge M(\mathbb{Q/Z},k)[-k]&\rightarrow \mathbf{K}_{t}\wedge M(\mathbb{Q/Z},k)[-k]\\
\end{align*}
are $\pi_{\ast}$-isomorphisms, where $\mathbf{K}_{a}$ and $\mathbf{K}_{t}$ are $\Omega$-$CW$-prespectra representing $0$-connective algebraic $K$-theory of the complex numbers (resp. real numbers) and $0$-connective complex (resp, real ) topological $K$-theory, respectively.  
\end{corollary}
\begin{proof}
It is proved in \cite[Cor.4.6]{Su} that the canonical map
\[K_{a}:=\operatorname{BGL}^{+}(\mathbb{F}^{\prime\delta})\rightarrow K_{t}:=\operatorname{BGL}(\mathbb{F}^{\prime})\]
induces an isomorphism 
\[\pi_{n}(K_{a})_{\mathbb{Z}/m} \xrightarrow{\sim} \pi_{n}(K_{t})_{\mathbb{Z}/m}, \text{ for } n\geq 1,\]
where $\mathbb{F}^{\prime}=\mathbb{R}$ or $\mathbb{C}$. This, in view of Corollary \ref{Cor:Towdefforhgpwithcoeff}, implies the following isomorphisms  
\begin{align*}
[\mathbf{S}[n],\mathbf{K}_{a}\wedge\mathbf{M}\mathbb{Z}/m]_{\operatorname{Ho}(\mathcal{P})}&\rightarrow [\mathbf{S}[n],\mathbf{K}_{t}\wedge\mathbf{M}\mathbb{Z}/m]_{\operatorname{Ho}(\mathcal{P})};\\
[\mathbf{S}[n],\mathbf{K}_{a}\wedge M(\mathbb{Z}/m,k)[-k]]_{\operatorname{Ho}(\mathcal{P})}&\rightarrow [\mathbf{S}[n],\mathbf{K}_{t}\wedge M(\mathbb{Z}/m,k)[-k]]_{\operatorname{Ho}(\mathcal{P})},\\
\end{align*}
for $n,m\in \mathbb{N}$. Having 
\[\operatorname*{hocolim}_{m}M(\mathbb{Z}/m,k)=M(\mathbb{Q/Z},k),\]
we can further deduce 
\[[\mathbf{S}[n],\mathbf{K}_{a}\wedge M(\mathbb{Q/Z},k)[-k]]_{\operatorname{Ho}(\mathcal{P})}\rightarrow [\mathbf{S}[n],\mathbf{K}_{t}\wedge M(\mathbb{Q/Z},k)[-k]]_{\operatorname{Ho}(\mathcal{P})}\]
is an isomorphism, for all $n\in\mathbb{N}$, and hence the canonical map of prespectra 
\[\mathbf{K}_{a}\wedge  M(\mathbb{Q/Z},k)[-k]\rightarrow \mathbf{K}_{t}\wedge  M(\mathbb{Q/Z},k)[-k]\]
is a $\pi_{\ast}$-isomorphism. As we have the naive smash product
$\mathbf{E}\wedge  M(\mathbb{Q/Z},k)[-k]$ is equivalent to $\mathbf{E}\wedge \mathbf{M}\mathbb{Q/Z}$ in $\operatorname{Ho}(\mathcal{A})$ \cite[p.258-259]{Sw} and therefore also in $\operatorname{Ho}(\mathcal{P})$ (see \ref{AandPcomparison}), the first statement follows. 
 
\end{proof}
\begin{remark}\label{rmk:OtherformulationofSu}
In view of Lemma \ref{NaturalisoFibandCofib}, Suslin's theorem can be rephrased as follows: Let $F_{a,m}$ and $F_{t,m}$ be the infinite loop space of the homotopy fiber of 
\begin{align*} 
\mathbf{K}_{a}\wedge S^{1}[-1]&\xrightarrow{id\wedge\varphi_{m}} \mathbf{K}_{a}\wedge S^{1}[-1]; \\ 
\mathbf{K}_{t}\wedge S^{1}[-1]&\xrightarrow{id\wedge\varphi_{m}} \mathbf{K}_{t}\wedge S^{1}[-1],
\end{align*}
respectively. Then the natural map 
\[\operatorname{Su}:F_{a,m}\rightarrow F_{t,m}\]
is a homotopy equivalence. 
\end{remark}

Now given a monomorphism $\phi:H\rightarrow G$, there is a map $\tilde{\phi}:M(H,k)\rightarrow M(G,k)$ that realizes $\phi$, and the homotopy cofiber of $\tilde{\phi}$ is a model for$M(\operatorname{coker}(\phi),k)$. In a similar way, with the construction of Moore prespectra given earlier, the cofiber of the induced map of $CW$-prespectra  
\[\bm{\phi}:\mathbf{M}H\rightarrow \mathbf{M}G,\]
which realizes $\phi$, is a model for the homotopy type of the Moore $CW$-prespectrum $\mathbf{M}(\operatorname{coker}(\phi))$. Particularly, there is a homotopy cofiber sequence in $\operatorname{Ho}(\mathcal{P})$ 
\[\mathbf{E}\wedge M(H,k)[-k]\rightarrow \mathbf{E}\wedge M(G,k)[-k]\rightarrow \mathbf{E}\wedge M(\operatorname{coker}(\phi),k)[-k].\]

Now the Serre theorem gives us the following $\pi_{\ast}$-isomorphism:  
\[\mathbf{M}\mathbb{Q}\simeq \mathbf{S}\wedge \mathbf{M}\mathbb{Q}\simeq \mathbf{H}\mathbb{Q},\]
so $\mathbf{M}\mathbb{Q}$ is rational, where $\mathbf{H}A$ is the Eilenberg-Maclane prespectrum of the abelian group $A$. By Lemma \ref{HomodetermapofSp} and Lemma \ref{RationalizationofSp}, we know there is a unique map of $CW$-prespectra in $\operatorname{Ho}(\mathcal{A})$
\[\mathbf{M}\mathbb{Z}=\mathbf{S}\rightarrow \mathbf{M}\mathbb{Q}=\mathbf{H}\mathbb{Q}\]
that realizes the canonical inclusion
\[\pi_{0}(\mathbf{M}\mathbb{Z})\simeq \mathbb{Z}\rightarrow \mathbb{Q}\simeq\pi_{0}(\mathbf{M}\mathbb{Q}),\]
and thus the following maps of $CW$-prespectra can be defined without ambiguity
\begin{align*}
\mathbf{E}\wedge \mathbf{M}\mathbb{Z}&\rightarrow \mathbf{E}\wedge \mathbf{M}\mathbb{Q},\\
\mathbf{E}\wedge M(\mathbb{Z},k)[-k]&\rightarrow \mathbf{E}\wedge M(\mathbb{Q},k)[-k].
\end{align*} 
Letting $\mathbf{E}=\mathbf{K}_{a}$ or $\mathbf{K}_{t}$, we obtain the following homotopy cofiber sequences:
\begin{align*}
\mathbf{K}_{a}\wedge M(\mathbb{Z},k)[-k] &\xrightarrow{\mathbf{f}} \mathbf{K}_{a}\wedge M(\mathbb{Q},k)[-k]\rightarrow \mathbf{K}_{a}\wedge M(\mathbb{Q/Z},k)[-k];\\
\mathbf{K}_{t}\wedge M(\mathbb{Z},k)[-k] &\xrightarrow{\mathbf{g}} \mathbf{K}_{t}\wedge M(\mathbb{Q},k)[-k]\rightarrow \mathbf{K}_{t}\wedge M(\mathbb{Q/Z},k)[-k],
\end{align*}
and Corollary \ref{Suslinthm1} says that the map
\[\mathbf{K}_{a}\wedge M(\mathbb{Q/Z},k)[-k]\rightarrow \mathbf{K}_{t}\wedge M(\mathbb{Q/Z},k)[-k],\]
or equivalently the map, 
\[\mathbf{Cofib}(\mathbf{f})\rightarrow \mathbf{Cofib}(\mathbf{g})\] 
is a $\pi_{\ast}$-isomorphism. In view of Lemma \ref{NaturalisoFibandCofib}, we obtain another version of Suslin's theorem \cite[Corollary.4.6]{Su}:

\begin{theorem}\phantomsection\label{Suslinthm2}
The induced map of prespectra
\[\mathbf{Fib}(\mathbf{f})\rightarrow \mathbf{Fib}(\mathbf{g})\]
is a $\pi_{\ast}$-isomorphism.
\end{theorem} 
\begin{proof}
The assertion follows quickly from Lemma \ref{NaturalisoFibandCofib} as we have the following commutative diagram:

\begin{center}
\begin{tikzpicture}
\node(Lu) at (0,2) {$\mathbf{Fib}(\mathbf{f})$};
\node(Ll) at (0,0) {$\Omega\mathbf{Cofib}(\mathbf{f})$};
\node(Ru) at (4,2) {$\mathbf{Fib}(\mathbf{g})$};
\node(Rl) at (4,0) {$\Omega\mathbf{Cofib}(\mathbf{g})$};
 
\path[->, font=\scriptsize,>=angle 90]

(Ll) edge (Rl)
(Lu) edge (Ru)
   
(Lu) edge node [right]{$\wr$}(Ll)
(Ru) edge node [right]{$\wr$}(Rl);

\end{tikzpicture}
\end{center}

\end{proof}






\subsection{Rationalization of simple spaces and $CW$-prespectra} 
\noindent 
\textbf{Rational spaces:}
A rational equivalence is a map $f:X\rightarrow Y\in \mathpzc{Top}_{\ast}$ such that the induced homomorphism 
\[f_{\ast}: \tilde{H}_{\ast}(X,\mathbb{Q})\rightarrow \tilde{H}_{\ast}(Y,\mathbb{Q})\] 
is an isomorphism, for all $\ast\in \mathbb{N}\cup\{0\}$.

Since in this paper we are mainly concerned with the rationalization of connected infinite loop spaces, we assume all spaces in this subsection are simple---a simple space is a connected space whose fundamental group acts trivially on its homotopy groups, even though most of the theorems stated here can be applied to more general spaces as well (see \cite[Chapter 5-6]{MP}).
 
\begin{lemma}
Given a simple space $Z$, the following statements are equivalent:
\begin{enumerate}
\item For any rational equivalence $f:X\rightarrow Y$, the induced map of sets
\[f^{\ast}:[Y,Z]\rightarrow [X,Z]\] is a bijection.
\item $\pi_{\ast}(Z)$ is rational, for $\ast\geq 1$.
\item $\tilde{H}_{\ast}(Z,\mathbb{Z})$ is rational, for $\ast\geq 1$.
\end{enumerate}
\end{lemma}
\begin{proof}
see \cite[Theorem 6.1.1]{MP}.
\end{proof} 
A space that satisfies any of these three conditions above is called a rational space.

\begin{definition}\label{Def:RationalizationTop}
Given a map of simple spaces $f:X\rightarrow Y$, $(Y,f)$ is called a rationalization of $X$ if and only if $Y$ is rational and $f$ is a rational equivalence.  
\end{definition}
The existence of the rationalization of a simple space is proved in \cite[Theorem 5.3.2]{MP}, and by Definition \ref{Def:RationalizationTop}, it also not difficult to see $Y$ is uniquely determined up to homotopy equivalence. We usually denote the rationalization of $X$ simply by $X_{\mathbb{Q}}$, dropping the map $X\rightarrow X_{\mathbb{Q}}$ from the notation.

The following statements characterize the rationalization of a simple space.
\begin{lemma}\phantomsection\label{Rationalizationofspaces}
The following statements are equivalent:
\begin{enumerate}
\item $f:X\rightarrow Y$ is a rationalization of $X$.
\item For any rational space $Z$, we have 
   \[f^{\ast}: [Y,Z]\rightarrow [X,Z]\] is a bijection.
\item The induced homomorphism
\[f_{\ast}:\pi_{\ast}(X)\rightarrow \pi_{\ast}(Y)\] is the rationalization of $\pi_{\ast}(X)$, for $\ast\geq 1$.
\item The induced homomorphism
\[f_{\ast}:\tilde{H}_{\ast}(X,\mathbb{Z})\rightarrow \tilde{H}_{\ast}(Y,\mathbb{Z})\] is the rationalization of $\tilde{H}_{\ast}(X,\mathbb{Z})$, for $\ast\geq 1$.
\end{enumerate}
\end{lemma}
\begin{proof}
See \cite[Theorem 6.1.2]{MP}.
\end{proof}

One simple observation is that, given a simple space $X$, if $\pi_{n}(X)$ is a torsion group, for every $n\geq 1$, then $X_{\mathbb{Q}}$ is contractible.

\vspace*{1.5em}
\noindent
\textbf{Rationalization of $CW$-prespectra}
\begin{definition}
Let $\mathbf{M}\mathbb{Q}$ be the Moore $CW$-prespectrum for rational numbers, $\mathbf{F}$ an $\Omega$-$CW$-prespectrum, and $f$ a map of $CW$-prespectra $\mathbf{E}\rightarrow \mathbf{F}$ in $\mathcal{P}$, then $(\mathbf{F},f)$ is a rationalization of $\mathbf{E}$ if and only if there is a $\pi_{\ast}$-isomorphism 
\[h:\mathbf{E}\wedge \mathbf{M}\mathbb{Q}\rightarrow \mathbf{F}\]
such that the following commutes

\begin{center}
\begin{tikzpicture}
\node(Lm) at (0,0)  {$\mathbf{E}$}; 
\node(Ru) at (2,1)  {$\mathbf{E}\wedge \mathbf{M}\mathbb{Q}$};
\node(Rl) at (2,-1) {$\mathbf{F}$};

\path[->, font=\scriptsize,>=angle 90]

(Lm) edge (Ru)  
(Lm) edge node [above]{$f$}(Rl) 
(Ru) edge node [right]{$h$}(Rl);
\end{tikzpicture}
\end{center}

A $CW$-prespectrum $\mathbf{E}$ is called rational if and only if the canonical map $\mathbf{E}\rightarrow \mathbf{E}\wedge\mathbf{M}\mathbb{Q}$ is a $\pi_{\ast}$-isomorphism. In this case, $\Omega^{\infty}\mathbf{E}$ is a rational space ( \ref{QadjunctionPresandspace}). In particular, $E_{0}$, the zero component of $\mathbf{E}$, is a rational space, if $\mathbf{E}$ is already an $\Omega$-prespectrum (fibrant object in $\mathcal{P}$).
\end{definition}

The following statements are equivalent and they all characterize a rationalization of a $CW$-prespectrum.
\begin{lemma}\phantomsection\label{RationalizationofSp}
\begin{enumerate}
\item $f:\mathbf{E}\rightarrow \mathbf{F}$ is a rationalization of $\mathbf{E}$.
\item For any $CW$-prespectrum $\mathbf{G}$, the induced homomorphism \[f_{\ast}:[\mathbf{S}[\ast],\mathbf{G}\wedge \mathbf{E}]_{\operatorname{Ho}(\mathcal{P})}\rightarrow [\mathbf{S}[\ast],\mathbf{G}\wedge\mathbf{F}]_{\operatorname{Ho}(\mathcal{P})}\] is a rationalization of the abelian group $[\mathbf{S}[\ast],\mathbf{G}\wedge\mathbf{E}]_{\operatorname{Ho}(\mathcal{P})}$.
\item Given any rational $CW$-prespectrum $\mathbf{G}$, we have 
\[f^{\ast}:[\mathbf{F},\mathbf{G}]_{\operatorname{Ho}(\mathcal{P})}\rightarrow [\mathbf{E},\mathbf{G}]_{\operatorname{Ho}(\mathcal{P})}\] is an isomorphism. 
\end{enumerate}
\end{lemma}
\begin{proof}
See \cite[5.4-5;5.8-9]{Ru} for the statements in terms of the Adams category $\mathcal{A}$, and use Theorem \ref{RelationbetweenAandP} to translate it into one in terms of $\mathcal{P}$.
\end{proof}

The next lemma says, any rational $CW$-prespectrum is equivalent to a graded Eilenberg-Maclane prespectrum.  
\begin{lemma}
Let $\mathbf{E}$ be a rational $CW$-prespectrum, then there is a $\pi_{\ast}$-isomorphism 
\[\mathbf{E}\xrightarrow{\sim} \prod_{n\in\mathbb{Z}}\mathbf{H}\pi_{n}(\mathbf{E})[n].\] 
\end{lemma}
\begin{proof}
See \cite[Theorem 7.11 (ii)]{Ru}.
\end{proof}
For instance, the rationalization of real topological $K$-theory prespectrum is given by
\[(\prod_{\mathclap{i\in\mathbb{N}\cup\{0\}}} \mathbf{H}\mathbb{Q}[4i],c\circ\operatorname{ch}),\] 
while   
\[(\prod_{\mathclap{i\in\mathbb{N}\cup\{0\}}} \mathbf{H}\mathbb{Q}[2i],\operatorname{ch})\]
is the rationalization of complex topological $K$-theory prespectrum, where $c$ is induced by the complexification of real vector bundles.

Another important feature of a rational $\Omega$-$CW$-prespectrum is that the homotopy classes of the maps into a rational $\Omega$-$CW$-prespectrum are determined completely by the induced homomorphisms between homotopy groups.  

\begin{lemma}\phantomsection\label{HomodetermapofSp}
Given two $CW$-prespectra $\mathbf{E}$ and $\mathbf{F}$ with $\mathbf{F}$ a rational $\Omega$-prespectrum, then
the canonical homomorphism 
\[[\mathbf{E},\mathbf{F}]_{\operatorname{Ho}(\mathcal{P})}\rightarrow \operatorname{Hom}^{0}(\pi_{\ast}(\mathbf{E})\otimes \mathbb{Q},\pi_{\ast}(\mathbf{F}))\]
is an isomorphism, where
\[\operatorname{Hom}^{0}(A_{\ast}.B_{\ast})\]
is the abelian group of graded homomorphisms of degree $0$.
\end{lemma}
\begin{proof}
See \cite[Theorem 7.11 (iii)]{Ru} where it is proved in terms of $\mathcal{A}$. The result follows from Theorem \ref{RelationbetweenAandP} and Remark \ref{Morenotionsonhomotopy} as we have the isomorphisms  
\[ [\mathbf{E},\mathbf{F}]_{\operatorname{Ho}(\mathcal{P})}\xrightarrow{\sim} [\mathbf{E},\mathbf{F}]_{l}\xrightarrow{\sim} [\mathbf{E},\mathbf{F}]_{\operatorname{Ho}(\mathcal{A})}.\]
\end{proof}

The next lemma shows the notions of homotopy and weakly homotopy are the same for maps from a $CW$-complex to a rational infinite loop space.
\begin{lemma}\phantomsection\label{Weakhtyandhty}
Given $\mathbf{E}$ a rational $\Omega$-prespectrum and $L$ a $CW$-complex, we have
\[[L,E_{0}]_{w}\simeq [L,E_{0}],\]
where $[L,Y]_{w}$ is the set of weakly homotopy classes of maps from $L$ to $Y$. Two maps $L$ to $Y$ are said to be weakly homotopic if and only if they are homotopic when restricted to every finite subcomplex of $L$.
\end{lemma} 
\begin{proof}
Recall first given two maps from a $CW$-prespectrum $\mathbf{F}$ to an $\Omega$-prespectrum $\mathbf{E}$ in $\mathcal{P}$, we say they are homotopic up to phantom maps if and only if they are homotopic when restricted to every finite $CW$-subprespectrum of $\mathbf{F}$. It is an equivalence relation on the set of maps from $\mathbf{E}$ to $\mathbf{F}$ and the corresponding set of equivalent classes is denoted by $[\mathbf{E},\mathbf{F}]_{\operatorname{Ho}(\mathcal{P}),w}$. One can check it is endowed with a natural abelian group structure.  

Secondly, we note there are isomorphisms:
\[[X,E_{0}]_{w}\simeq [\Sigma^{\infty}X,\mathbf{E}]_{\operatorname{Ho}(\mathcal{P}),w}\xleftarrow[\ref{HomodetermapofSp}]{\sim}[\Sigma^{\infty}X,\mathbf{E}]_{\operatorname{Ho}(\mathcal{P})}\simeq [X,E_{0}],\] 
The first isomorphism follows from the fact that every finite $CW$-subprespectrum of $\Sigma^{\infty}X$ is contained in $\Sigma^{\infty}C$, for some finite $CW$-complex $C$, the second isomorphism is the consequence of Lemma \ref{HomodetermapofSp}, and the last isomorphism results from the adjunction (Proposition \ref{QadjunctionPresandspace} 
\[\Sigma^{\infty}:\mathpzc{Top}_{\ast}\leftrightarrows \mathcal{P}:\Omega^{\infty}.\]
\end{proof}   
 
\section{Topological $K$-theory with coefficients} 
The purpose of this section is to show the geometric model defined in \cite[Section $5$]{APS2} and \cite[p.88]{APS3} and the geometric model considered in \cite{JW}\footnote{This geometric model is credited to Karoubi by Jones and Westbury.} are isomorphic. We shall also see both of them realize topological $K$-theory with coefficients in $\mathbb{C/Z}$ (Definition \ref{Cohomowithcoeff}) when they can be defined. The idea is to first identify the geometric model given by Atiyah, Patodi and Singer with the topological definition (Definition \ref{Cohomowithcoeff}) through Segal's construction (see \cite[p.88]{APS3} and the second model in Theorem \ref{CZKtheory}) and then construct an isomorphism between these two geometric models.

\subsection{Topological $K$-theory with coefficients in $\mathbb{Z}/m$:}
Here we review Atiyah, Patodi and Singer's geometric model of topological $K$-theory with coefficients in $\mathbb{Z}/m$. From now on, $\mathbf{K}_{t}$ stands only for the $\Omega$-$CW$-prespectrum that represents $0$-connective complex topological $K$-theory.  
    
\begin{lemma}\phantomsection\label{ZM}
\begin{enumerate}
\item 
Let $\mathbf{F}_{t,m}$ be the homotopy fiber of the map of prespectra
\[\mathbf{K}_{t}\wedge M(\mathbb{Z},k)[-k]=:\mathbf{K}_{t}\wedge_{na} \mathbf{M}\mathbb{Z}\xrightarrow{id\wedge \bm{\varphi}_{m}} \mathbf{K}_{t}\wedge_{na} \mathbf{M}\mathbb{Z}:=\mathbf{K}_{t}\wedge M(\mathbb{Z},k)[-k],\] where $\bm{\varphi}_{m}:\mathbf{M}\mathbb{Z}\rightarrow \mathbf{M}\mathbb{Z}$ is induce by the degree $m$ map of $S^{1}$. Then there is a natural $\pi_{\ast}$-isomorphism:
\[\mathbf{F}_{t,m}\rightarrow \Omega(\mathbf{K}_{t}\wedge_{na} \mathbf{M}\mathbb{Z}/m).\] 
  
\item 
There is a $1$-$1$ correspondence 
\[ [\Sigma^{\infty}X\wedge \mathbf{M}^{\ast}\mathbb{Z}/m,\mathbf{K}_{t}]_{\operatorname{Ho}(\mathcal{P})} \xleftrightarrow{1-1} [\Sigma^{\infty}X,\mathbf{K}_{t}\wedge\mathbf{M}\mathbb{Z}/m]_{\operatorname{Ho}(\mathcal{P})},\]
for any topological space $X$, where $\mathbf{M}^{\ast}\mathbb{Z}/m$ is the $S$-dual of $\mathbf{M}\mathbb{Z}/m$. 
\end{enumerate}
\end{lemma}

\begin{proof}
For the first assertion, we refer to Lemma \ref{NaturalisoFibandCofib} (see also \cite[Section $7$]{MMSS} or \cite[p.128]{LMS}). The second assertion follows from the definition of the $S$-duality \cite[Definition 14.20]{Sw}.  
\end{proof}
We use $\tilde{K}^{\ast}(X)_{\mathbb{Z}/m}$ to denote the abelian group $[\Sigma^{\infty}X,\mathbf{K}_{t}\wedge \mathbf{M}\mathbb{Z}/m]_{\operatorname{Ho}(\mathcal{P})}$.

In the following, we explain how to see geometric model defined by Atiyah, Patodi and Singer \cite[Section $5$]{APS2} indeed models topological $K$-theory with coefficients in $\mathbb{Z}/m$. 

\begin{theorem}\phantomsection\label{ZMG}
Let $L$ be a finite $CW$ complex. Then we have 
\begin{multline*}
\tilde{K}^{-1}(X)_{\mathbb{Z}/m}\simeq\{ (V,\phi) \mid  V \text{ is a vector bundle over } L\\ 
\text{ and } \phi:mV\oplus\epsilon^{\star} \xrightarrow{\sim} \epsilon^{\star} \text{ is an isomorphism}\}/\sim, 
\end{multline*}
where 
$(V,\phi) \sim (V^{\prime},\phi^{\prime})$ if and only if there is an isomorphism of vector bundles
\[\psi :V\oplus\epsilon^{\star}\xrightarrow{\sim} V^{\prime}\oplus\epsilon^{\star},\]
such that the following diagram commutes up to homotopy (of isomorphisms of vector bundles)
\begin{center}  
\begin{tikzpicture} 
\node (lupper)  at (0,3) {$mV\oplus \epsilon^{\star}$};
\node (rupper)  at (4,3) {$mV^{\prime}\oplus \epsilon^{\star}$};
\node (mlow)    at (2,0) {$\epsilon^{\star}$};
  
\path [->, font=\scriptsize, >=angle 90]
   
(lupper)  edge node[above]{m$\psi\oplus id$}(rupper);

\draw (lupper) to [out=-90,in=145]node[right]{\scriptsize $\phi\oplus id$} (mlow);
\draw (rupper) to [out=-90,in=35]node[left]{\scriptsize $\phi^{\prime}\oplus id$}(mlow);
\end{tikzpicture}
\end{center}  
where by $\epsilon^{\star}$ we understand a trivial vector bundle of suitable dimension.
\end{theorem}

The additive operation is the direct sum of vector bundles, the zero element is the trivial bundle with identity, and the inverse element of $(V,\phi)$ is $(V^{\perp},\phi^{\perp})$, where $V^{\perp}$ is the complement of $V$, namely, $V\oplus V^{\perp}=\epsilon^{\star}$, with $dim V^{\perp}> dim L$. $\phi^{\perp}$ is constructed as follows: Without loss of generality, we can assume $\phi:V\xrightarrow{\sim} \epsilon^{\star}$
\begin{center}  
\begin{tikzpicture} 
\node (lupper)  at (0,3) {$m (V\oplus V^{\perp})$};
\node (rupper)  at (4,3) {$m \epsilon^{dim V+dim V^{\perp}}$};
\node (mlow)    at (2,0) {$\epsilon^{dim V}$};

\draw [->] (lupper) to node [above]{\scriptsize $=$} (rupper);
\draw [left hook->] (mlow) to [out=150,in=-90] node [right]{\scriptsize $\phi^{-1}$}(lupper);
\draw [right hook->] (mlow) to [out=30,in=-90] (rupper); 
\end{tikzpicture}
\end{center}  
Because the dimension of $V^{\perp}$ is large enough, there is always an isomorphism $\phi^{\perp}:mV^{\perp}\rightarrow \epsilon^{\star}$ such that $\phi^{\perp}\oplus \phi$ homotopic to $id$ (see \cite[IX.1.1]{Ko}).

\begin{proof}
Note first $\tilde{K}^{-1}(X)_{\mathbb{Z}/m}=K^{-1}(X)_{\mathbb{Z}/m}$ as $\tilde{K}^{-1}(S^{0})=0$. Observe also that the second statement of Lemma \ref{ZM} and the Bott periodicity imply the following isomorphisms 
\begin{multline*}
[\Sigma^{\infty}X_{+}\wedge \Sigma^{3}\mathbf{M}^{\ast}\mathbb{Z}/m,\mathbf{K}_{t}]\xleftrightarrow{\text{Bott}\hspace*{.4em}\text{periodicity}}[\Sigma^{\infty}X_{+}\wedge \mathbf{M}^{\ast}\mathbb{Z}/m,\Sigma^{-1}\mathbf{K}_{t}]\\
\xleftrightarrow{S-\text{duality}} [\Sigma^{\infty}X_{+},\Sigma^{-1}\mathbf{K}_{t}\wedge\mathbf{M}\mathbb{Z}/m]=:K^{-1}(X)_{\mathbb{Z}/m}. 
\end{multline*} 
The first abelian group is exactly the one Atiyah, Patodi and Singer use to derive their geometric model and hence the assertion follows from the argument in \cite[p.428-429]{APS2}. One thing to be mindful of is in \cite{APS2} they consider a pair of vector bundles $(E,F)$ over $X$ with an isomorphism $\phi:mE\rightarrow mF$, whereas in the model presented here, $F$ is always trivial. This kind of normalization can always be achieved as every vector bundle over a finite $CW$-complex has a complement (Remark \ref{ExtendedbareAPS}).   
\end{proof}

\subsection{Topological $K$-theory with coefficients in $\mathbb{Q/Z}$}
We recall first the homotopy colimit 
\[\operatorname*{hocolim}_{m}\mathbf{M}\mathbb{Z}/m=\mathbf{M}\mathbb{Q/Z}\]
gives a model for the Moore $CW$-prespectrum of $\mathbb{Q/Z}$ and 
\[\mathbf{K}_{t}\wedge_{na}\mathbf{M}
\mathbb{Q/Z}=\operatorname*{hocolim}_{m}{\mathbf{K}_{t}\wedge_{na}\mathbf{M}\mathbb{Z}/m},\]
where the direct system is given by 
\[\mathbf{M}\mathbb{Z}/m  \xrightarrow{n} \mathbf{M}\mathbb{Z}/nm\]
with $n$ is the map induced by the degree $n$ map $\varphi_{n}:S^{1}\rightarrow S^{1}$.  
Let $\tilde{K}^{-1}(X)_{\mathbb{Q/Z}}$ be $[\Sigma^{\infty}X,\Sigma^{-1}\mathbf{K}_{t}\wedge \mathbf{M}\mathbb{Q/Z}]$, we have the following geometric model for topological $K$-theory with coefficients in $\mathbb{Q/Z}$.

\begin{corollary}\phantomsection\label{QZ} 
 
Let $L$ be a finite $CW$-complex. Then $\tilde{K}^{-1}(L)_{\mathbb{Q/Z}}$ is isomorphic to the following abelian group.
\[\{ (V,\phi)\mid \text{ There exists } m\in\mathbb{N}  \text{ s.t. } \phi: mV\oplus\epsilon^{\ast}\xrightarrow{\sim} \epsilon^{\ast}\}/\sim.\]
$(V,\phi)\sim (V^{\prime},\phi^{\prime})$ if and only if there exist $n$ and $n^{\prime}$ such that $nm=n^{\prime}m^{\prime}$ and $(V,\prescript{}{n}\phi)$, $(V^{\prime},\prescript{}{n^{\prime}}\phi^{\prime})$ represent the same element in $\tilde{K}^{-1}(L)_{\mathbb{Z}/nm}$, where $(V,\prescript{}{n}\phi)$ (resp. $(V,\prescript{}{n}\phi)$) is the image of $(V,\phi)$ (resp. $(V,\phi^{\prime})$) under the homomorphism
\[\tilde{K}^{-1}(L)_{\mathbb{Z}/m}\rightarrow \tilde{K}^{-1}(L)_{\mathbb{Z}/nm}\]
which is induced by 
\begin{equation}\label{MapinQZ}
\mathbf{M}\mathbb{Z}/m \xrightarrow{n} \mathbf{M}\mathbb{Z}/nm.
\end{equation}
 
\end{corollary}
\begin{proof}
The isomorphism $\prescript{}{n}\phi$ is constructed as follows: Firstly, observe $\phi$ induces an isomorphism
\[(1_{X}\times S^{1}\varphi_{n})^{\ast}(mV\otimes(H\oplus (n-1)\epsilon^{1}))\xrightarrow{\sim}(1_{X}\times S^{1}\varphi_{n})^{\ast}(\epsilon^{\star}\otimes(H\oplus (n-1)\epsilon^{1})),\]
where $H$ is the Hopf bundle over $S^{2}$, and $\varphi_{n}$ the degree $n$ map from $S^{1}$ to itself. So $mV\otimes (H\oplus (n-1)\epsilon^{1})$ is a vector bundle over $L\times S^{2}$. Secondly, recall that there is a unique isomorphisms (up to homotopy) 
\[h:H^{n}\oplus (n-1)\epsilon^{1}\simeq nH.\] Then $\prescript{}{n}\phi$ is the isomorphism $nmV\rightarrow \epsilon^{\star}$ such that the following commutes
\begin{center}
\begin{equation}\label{Diag:Realizingmtonm}
\begin{tikzpicture}[baseline=(current bounding box.center)] 
\node(Lu) at (0,2) {$(1_{L}\times S^{1}\varphi_{n})^{\ast}(mV\otimes(H\oplus (n-1)\epsilon^{1}))$};
\node(Ru) at (7,2) {$(1_{L}\times S^{1}\varphi_{n})^{\ast}(\epsilon^{\star}\otimes(H\oplus (n-1)\epsilon^{1}))$};
\node(Ll) at (0,0) {$V\otimes nmH$};
\node(Ml1) at (2,0) {$nm V\otimes H$};
\node(Ml2) at (5,0) {$\epsilon^{\star}\otimes H$};
\node(Rl) at (7,0) {$\epsilon^{\star}\otimes nmH$};

\path[->, font=\scriptsize,>=angle 90]

 (Ru) edge node [right]{$id\otimes h$}(Rl)
(Lu) edge node [left]{$id\otimes h$}(Ll);

\draw [->](Ml1) to [out=-30,in=-150] node [below]{$\prescript{}{n}\phi\otimes id$}(Ml2);
\draw [->](Lu) to [out=30,in=150] node [above]{\scriptsize $(1_{L}\times S^{1}\varphi_{n})^{\ast}(\phi\otimes id)$} (Ru);
\draw[double equal sign distance] (Rl) to (Ml2);

\draw[double equal sign distance](Ll) to (Ml1);
\end{tikzpicture}
\end{equation}
\end{center}
Comparing with the construction in \cite[p.428-9]{APS2}, we see this assignment
\[(V,\phi)\mapsto (V,\prescript{}{n}\phi)\]
realizes the homomorphism 
\[\tilde{K}^{-1}(L)_{\mathbb{Z}/m}\simeq \tilde{K}(L_{+}\wedge M(\mathbb{Z}/m,2))\rightarrow \tilde{K}(L_{+}\wedge M(\mathbb{Z}/nm,2))\simeq \tilde{K}^{-1}(X)_{\mathbb{Z}/nm}\]
induced by the dual of map \eqref{MapinQZ}. The claim follows.
\end{proof} 
 
\begin{remark}\phantomsection\label{Moreoncompatibility}
\begin{enumerate}
\item One might hope $\prescript{}{n}\phi$ is homotopic (through isomorphisms of vector bundles) to $n\phi$, yet it is not clear if $n\phi$ makes diagram \eqref{Diag:Realizingmtonm} commute. Nevertheless, $n\phi$ and $\prescript{}{n}\phi$ are not far away. For instance, if we pick up a connection $\nabla_{v}$ of $V$, then we have $\frac{1}{n}\operatorname{Tch}(nm\nabla_{v},\prescript{}{n}\phi^{\ast}d)=\operatorname{Tch}(m\nabla_{v},\phi^{\ast}d)$---This can be seen from the commutative diagram in the proof.    
\item We denote the element in $\tilde{K}^{-1}(L)_{\mathbb{Q/Z}}$ by $[(V,\phi)_{m}]$ to indicate the element $(V,\phi)$ is living in $\tilde{K}^{-1}(L)_{\mathbb{Z}/m}$. 

\end{enumerate}
\end{remark}

\subsection{Topological $K$-theory with coefficients in $\mathbb{C/Z}$}
We are now ready to compare the geometric models defined in \cite{APS3} and \cite{JW}. We prove that these two models are indeed equivalent to topological $K$-theory with coefficients in $\mathbb{C/Z}$ when they are defined. A concrete isomorphism between these two geometric models is given in the proof of Theorem \ref{CZG}. Let \[\tilde{K}^{\ast}(X)_{\mathbb{C/Z}}:=[\Sigma^{\infty}X,\Sigma^{\ast}\mathbf{K}_{t}\wedge \mathbf{M}\mathbb{C/Z}],\]
for $X\in \mathpzc{Top}_{\ast}$.

\begin{theorem}\phantomsection\label{CZKtheory}
\begin{enumerate}
\item  There is a natural $\pi_{\ast}$-isomorphism in $\mathcal{P}$:
\[\mathbf{K}_{t}\wedge \mathbf{M}{\mathbb{C/Z}}\simeq \Sigma^{1}\mathbf{Fib}(\mathbf{ch}),\] 
where $\mathbf{Fib}(\mathbf{ch})$ is the homotopy fiber of the Chern character \[\mathbf{ch}:\mathbf{K}_{t}\rightarrow \prod\limits_{i}\mathbf{H}\mathbb{C}[2i],\] 
and by $\mathbf{E}[j]$ we understand the $k$-fold shifting of $\mathbf{E}$. In particular, we have
\[\tilde{K}^{\ast-1}(X)_{\mathbb{C/Z}}\simeq[\Sigma^{\infty}X,\Sigma^{\ast}\mathbf{Fib}(\mathbf{ch})]\simeq [X,F_{t,\mathbb{C/Z}}].\]

\item There is an isomorphism
 \[\operatorname{Coker}[\tilde{K}^{\ast}(X)_{\mathbb{Q}}\xrightarrow{(-j_{\ast},p_{\ast})} \tilde{K}^{\ast}(X)_{\mathbb{C}}\oplus \tilde{K}^{*}(X)_{\mathbb{Q/Z}}]\simeq\tilde{K}^{\ast}(X)_{\mathbb{C/Z}},\]
for any topological space $X$,
where $j_{\ast}$ is the map induced by $\mathbf{M}\mathbb{Q}\rightarrow\mathbf{M}\mathbb{C}$, and $p_{\ast}$ the natural map from $\mathbf{M}\mathbb{Q}$ to the cofiber of $\mathbf{M}\mathbb{Z}\rightarrow\mathbf{M}\mathbb{Q}$. By $\tilde{K}^{\ast}(X)_{\mathbb{F}}$ we understand $\tilde{K}^{\ast}(X)\otimes\mathbb{F}$ when $\mathbb{F}=\mathbb{Q},\mathbb{R}$ or $\mathbb{C}$.\footnote{This construction is credited to Segal by Atiyah et al. \cite{APS3}.}


\end{enumerate}

\end{theorem}


\begin{proof}
 
Observe first the natural inclusions 
\[\mathbb{Z}\hookrightarrow \mathbb{Q}\hookrightarrow \mathbb{C}\] 
give us the following maps of Moore $CW$-prespectra  
\[\mathbf{M}\mathbb{Z}\rightarrow \mathbf{M}\mathbb{Q}\rightarrow \mathbf{M}\mathbb{C}.\]
These maps are uniquely determined as both prespectra $\mathbf{M}\mathbb{Q}\cong \mathbf{H}\mathbb{Q}$ and $\mathbf{M}\mathbb{C}\cong \mathbf{H}\mathbb{C}$ are rational (see Lemma \ref{HomodetermapofSp}).

Recall also, given an injective homomorphism $A\rightarrow B$, the homotopy cofiber of the induced map of Moore $CW$-prespectra 
\[\mathbf{M}A\rightarrow\mathbf{M}B\]
is a model of the Moore $CW$-prespectrum $\mathbf{M}B/im(A)$. Combining these facts with the following commutative diagram:
\begin{center}
\begin{tikzpicture} 
\node(Ll) at (0,0) {$\mathbf{K}_{t}\wedge\mathbf{M}\mathbb{Z}$};
\node(Mu) at (3,2) {$\prod\limits_{i}\mathbf{H}\mathbb{Q}[2i]$};
\node(Ml) at (3,0) {$\mathbf{K}_{t}\wedge\mathbf{M}\mathbb{Q}$};
\node(Ru) at (6,2) {$\prod\limits_{i}\mathbf{H}\mathbb{C}[2i]$};
\node(Rl) at (6,0) {$\mathbf{K}_{t}\wedge\mathbf{M}\mathbb{C}$};

\path[->, font=\scriptsize,>=angle 90]

(Ll) edge (Ml)
(Ml) edge (Rl)
(Ll) edge node [above]{$\mathbf{ch}$}(Mu)
(Mu) edge (Ru)  
(Ml) edge node [right]{$\mathbf{ch}_{\otimes \mathbb{Q}}$} node[left]{$\wr$}(Mu) 
(Rl) edge node [right]{$\mathbf{ch}_{\otimes \mathbb{C}}$} node[left]{$\wr$}(Ru);

\end{tikzpicture}
\end{center}
we obtain 
\[\mathbf{Cofib}(\mathbf{ch}:\mathbf{K}_{t}\rightarrow \prod\limits_{n}\mathbf{H}\mathbb{C}[2i])\simeq \mathbf{K}_{t}\wedge \mathbf{M}\mathbb{C/Z},\]
where $\mathbf{Cofib}(\mathbf{ch})$ is the homotopy cofiber of $\mathbf{ch}$. The first assertion then follows from Lemma \ref{HofibercofiberseqinP} (see also \cite[p.126]{LMS}), which says there is a natural $\pi_{\ast}$-isomorphism from $\mathbf{Fib}(\mathbf{f})\rightarrow \Omega \mathbf{Cofib}(\mathbf{f})$ for any map of prespectra $\mathbf{f}$.
   
Secondly, as any map $\mathbf{M}\mathbb{Z}\rightarrow \Omega\mathbf{M}\mathbb{C/Z}$ is null-homotopic, there is a unique map of prespectra 
\[i :\mathbf{M}\mathbb{Q/Z}\rightarrow\mathbf{M}\mathbb{C/Z}\]
making the following commutes
\begin{center}
\begin{equation}\label{ZQCMoorespectra}
 \begin{tikzpicture}[baseline=(current bounding box.center)] 
\node(Lu) at (0,1) {$\mathbf{M}\mathbb{Z}$};
\node(Ll) at (0,0) {$\mathbf{M}\mathbb{Z}$};
\node(Mu) at (3,1) {$\mathbf{M}\mathbb{Q}$};
\node(Ml) at (3,0) {$\mathbf{M}\mathbb{C}$};
\node(Ru) at (6,1) {$\mathbf{M}\mathbb{Q/Z}$};
\node(Rl) at (6,0) {$\mathbf{M}\mathbb{C/Z}$};

\path[->, font=\scriptsize,>=angle 90]

(Ll) edge (Ml)
(Ml) edge node [above]{$\pi$}(Rl)
(Lu) edge (Ll)
(Lu) edge node [above]{$k$}(Mu)
(Mu) edge node [above]{$p$}(Ru)  
(Mu) edge node [right]{$j$}(Ml)  
(Ru) edge node [right]{$i$}(Rl);

\end{tikzpicture} 
\end{equation}
 \end{center}
Hence, to see the second claim, it suffices to show the following sequence is exact 
  
\[\tilde{K}^{\ast}(X)_{\mathbb{Q}}\xrightarrow{(-j_{\ast},p_{\ast})} \tilde{K}^{\ast}(X)_{\mathbb{C}}\oplus \tilde{K}^{\ast}(X)_{\mathbb{Q/Z}}\xrightarrow{\pi_{\ast}+i_{\ast}} \tilde{K}^{\ast}(X)_{\mathbb{C/Z}}\rightarrow 0.\]
For the surjectivity of $\pi_{\ast}+i_{\ast}$, we consider the following commutative diagram of exact sequences. We have dropped the topological space $X$ and $\tilde{(-)}$ from the notations to simplify the presentation. 
\begin{center}
\begin{tikzpicture} 
\node (A) at (-1,2) {$K^{*}_{\mathbb{C}}$};
\node (B) at (-1,0) {$K^{*}_{\mathbb{Q/Z}}$};
\node (C) at (1,1) {$K^{*}_{\mathbb{C/Z}}$};
\node (D) at (3,0) {$K^{*+1}_{\mathbb{Z}}$};
\node (E) at (5,0) {$K^{*+1}_{\mathbb{Q}}$};
\node (F) at (5,-1) {$K^{*+1}_{\mathbb{C}}$};
\node (G) at (5,-2) {$0$};

\path [->,font=\scriptsize,>=angle 90]
(A) edge node [above]{$\pi_{\ast}$}(C)
(B) edge node [above]{$i_{\ast}$}(C)
(B) edge (D)
(C) edge (D)
(D) edge (E)
(D) edge (F)
(E) edge node[right]{i}(F)
(F) edge (G);
  
\end{tikzpicture}
\end{center}  
Given an element $x\in K^{\ast}_{\mathbb{C/Z}}$, we have its image in $K^{\ast+1}_{\mathbb{C}}$ is zero. Since $K^{\ast+1}_{\mathbb{Q}}\rightarrow K^{\ast+1}_{\mathbb{C}}$ is injective, we know the image of $x$ in $K^{\ast+1}_{\mathbb{Q}}$ is also zero. Therefore there is an element $y\in K^{\ast}_{\mathbb{Q/Z}}$ such that its image in $K^{\ast+1}_{Z}$ coincides with the image of $x$. Now the difference $x-i_{\ast}(y)$ in $K^{\ast}_{\mathbb{C/Z}}$ can be corrected by an element $z\in K^{\ast}_{\mathbb{C}}$. That is $\pi_{\ast}(z)+i_{\ast}(y)=x$, and hence $\pi_{\ast}+j_{\ast}$ is onto.

Next, we note commutative diagram \eqref{ZQCMoorespectra} implies 
\[im(-j_{*}, p_{*}) \subset \operatorname{ker}(\pi_{*} + i_{*}).\]  
As to the other direction, it can be deduced from the following commutative diagram of exact sequences: 
\begin{center}  
\begin{tikzpicture}[scale=1]
\node (Z)  at (0,3) {$K^{*}_{\mathbb{Z}}$};
\node (Q)  at (0,1) {$K^{*}_{\mathbb{Q}}$};
\node (C)  at (2,2) {$K^{*}_{\mathbb{C}}$};
\node (Q/Z)at (2,0) {$K^{*}_{\mathbb{Q/Z}}$};
\node (C/Z)at (4,1) {$K^{*}_{\mathbb{C/Z}}$};
\node (Z1) at (6,0) {$K^{*+1}_{\mathbb{Z}}$};
   
\path [->, font=\scriptsize, >=angle 90]
   
(Z)  edge (C)
(Z)  edge node [left]{$k_{\ast}$}(Q)
(Q)  edge node [above]{$-j_{\ast}$}(C)
 
(Q/Z)edge node [below]{$i_{\ast}$}(C/Z)
(C)  edge node [above]{$\pi_{\ast}$}(C/Z)
(C/Z)edge (Z1)
(Q/Z)edge (Z1);
\draw[->](Q) to [out=-90,in=180]node [above]{\tiny $p_{\ast}$} (Q/Z);  
\end{tikzpicture}
\end{center}   
More precisely, suppose $(z,y)$ is in the kernel of $\pi_{\ast} + i_{\ast}$, namely, $\pi_{\ast}(z)= -i_{\ast}(y)$. Then by the exactness of 
\[K^{\ast}_{\mathbb{Q}}\rightarrow K^{\ast}_{\mathbb{Q/Z}}\rightarrow K^{\ast+1}_{\mathbb{Z}},\] 
there is an element $u$ in $K^{\ast}_{\mathbb{Q}}$ which has its image equal to $y$ in $K^{\ast}_{\mathbb{Q/Z}}$. Therefore the element $z+j_{\ast}(u)$ has its image in $K^{\ast}_{\mathbb{C/Z}}$ is zero. Now, by the exactness of 
\[K^{\ast}_{\mathbb{Z}}\rightarrow K^{\ast}_{\mathbb{C}}\rightarrow K^{\ast}_{\mathbb{C/Z}},\] 
one can find an element $v$ in $K^{*}_{\mathbb{Z}}$ such that its image in $K^{\ast}_{\mathbb{C}}$ equal to $z+j_{\ast}(u)$. Then one checks $u-k_{\ast}(v)$ has its image equal to $(z,y)$ in $K^{\ast}_{\mathbb{C}}\oplus K^{\ast}_{\mathbb{Q/Z}}$ under $(-j_{\ast},p_{\ast})$ and hence we have shown 
\[im(-j_{*}, p_{*}) \supset ker(\pi_{*} + i_{*}).\]
  
\end{proof}

The next theorem is the main result of this subsection. We identify the geometric models defined in \cite{APS3} and \cite{JW}
and show they indeed model topological $K$-theory with coefficients in $\mathbb{C/Z}$ when they can be defined.

\begin{theorem}\phantomsection\label{CZG}
\begin{enumerate}
\item
Given a finite $CW$-complex $L$, we have  
\begin{multline*}
\tilde{K}^{-1}(L)_{\mathbb{C/Z}}\simeq\mathcal{G}_{\operatorname{APS}}(L):=\{(\omega,[(V,\phi)_{m}])\mid \phi: mV\oplus \epsilon^{\star} \xrightarrow{\sim} \epsilon^{\star}; \hspace{.5em}\\
 \omega \text{ a closed odd cocycle} \} / \sim,
\end{multline*}
where two elements are equivalent $(\omega,[(V,\phi)_{m}]) \sim (\omega^{\prime},[(V^{\prime},\phi^{\prime})_{m^{\prime}}])$ if and only if 
there exists two natural numbers $l$ and $l^{\prime}$ and an isomorphism $\psi:\epsilon^{\star}\rightarrow \epsilon^{\star}$ such that 
\[\begin{cases}
ml=m^{\prime}l^{\prime};\\  
(V,\prescript{}{l}\phi)_{ml}\sim (V^{\prime},\psi\circ \prescript{}{l^{\prime}}\phi^{\prime})_{m^{\prime}l^{\prime}}\in \tilde{K}^{-1}(L)_{\mathbb{Z}/ml};\\
[\omega ]-[\omega^{\prime}]=-[\frac{1}{ml}\operatorname{Tch}(\psi^{\ast}d,d)]\in \tilde{H}^{odd}(X,\mathbb{C}). 
\end{cases}\]

Given two elements $(\omega,[(V,\phi)_{m}])$ and $(\omega^{\prime},[(V^{\prime},\phi^{\prime})_{m^{\prime}}])$---without loss of generality, one can assume $m=m^{\prime}$, the addition is give by 
\[(\omega,[(V,\phi)_{m}])+(\omega^{\prime},[(V^{\prime},\phi^{\prime})_{m^{\prime}}]):=(\omega+\omega^{\prime},V\oplus V^{\prime},\phi\oplus\phi^{\prime})_{m}.\]
$(0,(\epsilon^{\star},id)_{1})$ represents the zero element, and the inverse element of $([\omega],[(V,\phi)_{m}])$ is $(-\omega,[(V^{\perp},\phi^{\perp})_{m}])$.

\item If $L=M$, a compact smooth manifold, then  
\begin{multline*}
\tilde{K}^{-1}(M)_{\mathbb{C/Z}}\simeq\mathcal{G}_{\operatorname{JW}}(M):=\{(V,\nabla,\omega)\mid \nabla \text{ is a connection on } V \text{ and }\\
\omega \text{ an odd differntial form such that }
d\omega =ch(\nabla)\} / \sim,
\end{multline*} 
where $(V,\nabla,\omega)\sim (V^{\prime},\nabla^{\prime},\omega^{\prime})$ 
if and only if there exists $(U,\nabla_{U})$, a vector bundle with connection over $X\times I$, such that 
$(U,\nabla_{U})$ restricts to $(V,\nabla)$ and $(V^{\prime},\nabla^{\prime})$ stably on $X\times\{1\}$ and $X\times\{0\}$ and 
\[\int_{t\in I} ch(\nabla_{U})=\omega-\omega^{\prime}+d\eta,\] where $\eta$ is an even differential form. 

Given two elements $(V,\nabla,\omega)$, $(V^{\prime},\nabla^{\prime},\omega^{\prime})$, the addition is given by
\[(V\oplus V^{\prime},\nabla\oplus\nabla^{\prime},\omega+\omega^{\prime}).\]       
The zero element is the trivial bundle with trivial connection and zero form, $(\epsilon^{\ast},d,0)$, and the inverse element of $(V,\nabla,\omega)$ is $ (V^{\perp},\nabla^{\perp},\omega^{\perp})$, where $\nabla^{\perp}$ is any connection on $V^{\perp}$, and $\omega^{\perp}:=-\omega+\int ch(\nabla_{U})$ with $U$ being any vector bundle over $X\times I$ that restricts to $V\oplus V^{\perp}$ and $\epsilon^{\star}$ on $X\times 1$ and $X\times 0$ and $\nabla_{U}$ a connection on $U$ that restricts to $\nabla\oplus \nabla^{\perp}$ and $d$ on $X\times 1$ and $X\times 0$. 

\end{enumerate} 
\end{theorem}

\begin{proof} 
Step $1.$: Observe that the geometric model $\mathcal{G}_{\operatorname{APS}}(L)$ can be deduced from the second topological definition\footnote{It is credited to Segal by Atiyah, Patodi and Singer \cite{APS3}.} of topological $K$-theory with coefficients in $\mathbb{C/Z}$ in Theorem \ref{CZKtheory} and the geometric model in Corollary \ref{QZ}. This shows the first assertion. 

Step $2.$: Recall that, given two connections $\nabla$ and $\nabla^{\prime}$, the associated Chern-Simon class, denoted by $\operatorname{Tch}(\nabla,\nabla^{\prime})$ is well-defined in $\Omega^{odd}(M)/d(\Omega^{even}(M))$ and has $\operatorname{ch}(\nabla)-\operatorname{ch}(\nabla^{\prime})=d\operatorname{Tch}(\nabla,\nabla^{\prime})$, where $\Omega^{even(odd)}(M)$ is the graded group of even (odd) differential forms of $M$ and $d$ is the differential. The explicit construction for $\operatorname{Tch}(\nabla,\nabla^{\prime})$ used here is given by the following integral
\[\int_{t}\operatorname{ch}(\nabla_{t}),\]
where $\nabla_{t}=t\nabla+(1-t)\nabla^{\prime}$.

Chern-Simon classes satisfy the following equality in $\Omega^{odd}(M)/d(\Omega^{even}(M))$: 
\begin{equation}\label{IdentityforTch}
\operatorname{Tch}(\nabla,\nabla^{\prime})+\operatorname{Tch}(\nabla^{\prime},\nabla^{\prime\prime})=\operatorname{Tch}(\nabla,\nabla^{\prime\prime}).
\end{equation}
 
Step $3.$: When $L=M$, a compact smooth manifold, we construct a homomorphism from the first geometric model to the second:  
\begin{align}\label{eq:IsomPhi}
\Phi : \mathcal{G}_{\operatorname{APS}}(M) &\rightarrow \mathcal{G}_{\operatorname{JW}}(M)\\
(\omega, [(V,\phi)_{m}])&\mapsto [(V,\nabla_{v},\frac{1}{m}\operatorname{Tch}(m\triangledown_{v},\phi^{*}d)-\omega)]\nonumber 
\end{align}
where a connection $\nabla_{v}$ on $V$ is chosen.
    
To see it is well-defined, one observes that choosing different connections on $V$ does not change the class. That is because we can define 
\[(U:=V\times I,\nabla_{u}:=t\nabla_{v}+(1-t)\nabla^{\prime}_{v})\]
which restricts to $(V,\nabla_{v})$ and $(V,\nabla_{v}^{\prime})$ on $M\times \{1\}$ and $M\times\{0\}$.
By \eqref{IdentityforTch}, we obtain the following equality in $\Omega^{odd}(M)/d(\Omega^{even}(M))$:
\[\operatorname{Tch}(\nabla_{u}):=\operatorname{Tch}(\nabla_{v},\nabla_{v}^{\prime})=\frac{1}{m}\operatorname{Tch}(m\nabla_{v},\phi^{\ast}d)-\frac{1}{m}\operatorname{Tch}(m\nabla_{v}^{\prime},\phi^{\ast}d).\]
Hence the equivalent class of $\Phi(\omega,[(V,\phi)_{m}])$ is independent of the choice of connections on $V$.

Now let $(\omega^{\prime},[(V^{\prime},\phi^{\prime})_{m^{\prime}}])$ be another representative---without loss of generality, one can assume $m^{\prime}=m$. By the definition of $\mathcal{G}_{APS}(M)$, we know they differ by an element $(\epsilon^{\star},\psi)$ in $\tilde{{K}}^{-1}(M)_{\mathbb{Q}}$. That is we have an equivalence:
\begin{equation}\label{Eq2inCZKtheory}
(V,\phi)\sim (V^{\prime},\psi\circ\phi^{\prime})\in \tilde{K}^{-1}(M)_{\mathbb{Z}/m}, \text{ for some } m,
\end{equation}
and 
\[\frac{1}{m}\operatorname{Tch}(\psi^{\ast}d,d)+\omega=\frac{1}{m}ch(\psi)+\omega=\omega^{\prime}\text{ in  } \Omega^{odd}(M)/d(\Omega^{even}(M)).\]
Equivalence \eqref{Eq2inCZKtheory} implies an isomorphism of vector bundles $\rho:V^{\prime}\rightarrow V$ with $m\rho$ homotopic (through isomorphisms of vector bundles) to  
\[\phi^{-1}\circ\psi\circ\phi^{\prime}:V^{\prime}\rightarrow V.\] 
Now observe that $\rho$ induces a vector bundle with connection $(U,\nabla_{U})$ which restricts to $(V,\nabla_{v})$ and $(V^{\prime},\nabla_{V^{\prime}})$ on $M\times\{1\}$ and $M\times\{0\}$, respectively, and note also, as $m\rho$ is homotopic to $\phi^{-1}\circ \psi\circ \phi^{\prime}$, we can deduce the following equality in $\Omega^{odd}(M)/d(\Omega^{even}(M))$:
\begin{multline*} 
\int_{t} ch(\nabla_{u})\\
=\frac{1}{m}ch(m\rho)=\frac{1}{m}\operatorname{Tch}((m\rho)^{\ast}m\nabla_{v},m\nabla_{v^{\prime}})=\frac{1}{m}\operatorname{Tch}(\phi^{\prime,\ast}\circ\psi^{\ast}\circ\phi^{\ast,-1}m\nabla_{v},m\nabla_{v^{\prime}})\\
=\frac{1}{m}(\operatorname{Tch}(\phi^{\prime,\ast}\circ\psi^{\ast}\circ\phi^{\ast,-1}m\nabla_{v},\phi^{\prime,\ast}\circ\psi^{\ast}d)+\operatorname{Tch}(\phi^{\prime,\ast}\circ\psi^{\ast}d, \phi^{\prime,\ast}d)+\operatorname{Tch}(\phi^{\prime,\ast}d,m\nabla_{v^{\prime}}))\\
=\frac{1}{m}(\operatorname{Tch}(m\nabla_{v},\phi^{\ast}d)-\operatorname{Tch}(m\nabla_{v^{\prime}},\phi^{\prime,\ast}d))+\omega^{\prime}-\omega+d\eta.
\end{multline*} 
Hence these two elements
\[(V,\nabla_{v},\frac{1}{m}\operatorname{Tch}(m\nabla_{v},\phi^{\ast}d)-\omega)\]
and
\[(V^{\prime},\nabla_{v^{\prime}},\frac{1}{m}\operatorname{Tch}(m\nabla_{v^{\prime}},\phi^{\prime,\ast}d)-\omega^{\prime})\]
are equivalent in $\mathcal{G}_{JW}(M)$.

Step $4.$ The inverse map of $\Phi$ is given by the assignment:
\begin{align}\label{eq:IsomPsi}
\Psi : \mathcal{G}_{\operatorname{JW}}(M) & \rightarrow \mathcal{G}_{\operatorname{APS}}(M)\\  
[(V,\nabla_{v},\omega)] &\mapsto (\frac{1}{m}\operatorname{Tch}(m\triangledown_{v},\phi^{*}d)-\omega,[(V,\phi)_{m}])\nonumber
\end{align}
where a (stable) trivialization 
$\phi:mV\oplus\epsilon^{\star}\xrightarrow{\sim} \epsilon^{\star}$ is chosen.  
 
Suppose $\phi^{\prime}:m^{\prime}V\oplus\epsilon^{\star}\xrightarrow{\sim} \epsilon^{\star}$ is another trivialization---without loss of generality, we can assume $m=m^{\prime}$. Then the following identities:
\begin{multline*}
\frac{1}{m}\operatorname{Tch}(m\nabla_{v},\phi^{\ast}d)-\omega-(\frac{1}{m}\operatorname{Tch}(m\nabla_{v^{\prime}},\phi^{\prime\ast}d)-\omega)
=\frac{1}{m}\operatorname{Tch}(\phi^{\prime\ast}d,\phi^{\ast}d)\\
=\frac{1}{m}\operatorname{Tch}(d,(\phi\circ\phi^{\prime,-1})^{\ast}d)=-\frac{1}{m}\operatorname{ch}(\phi\circ\phi^{\prime,-1}) \text{ in } \Omega^{odd}(M)/d(\Omega^{even}(M)) 
\end{multline*}
and 
\[(V,\phi^{\prime})+(\epsilon^{\star},\phi \circ\phi^{\prime -1})=(V,\phi)\]
show the equivalent class of $\Psi(V,\nabla_{v},\omega)$ is independent of the choice of trivializations.

Now suppose $(V^{\prime},\nabla^{\prime},\omega^{\prime})$ is another representative. Then, by the definition of $\mathcal{G}_{\operatorname{JW}}(M)$, there exists $(U,\nabla_{u})$ over $M\times I$ such that it restricts to $(V,\nabla_{v})$ and $(V^{\prime},\nabla_{v^{\prime}})$ on $M\times\{1\}$ and $M\times\{0\}$, respectively, with 
\[d\int_{t\in I}  ch(\nabla_{u})=\omega-\omega^{\prime}+d\eta,\]
where $\eta$ is an even differential form. 
Now let $\psi:V^{\prime}\rightarrow V$ be the isomorphism induced from the homotopy. Then we have the following equivalence:
\[(\epsilon^{\star},\phi\circ m\psi \circ \phi^{\prime,-1})+(V^{\prime},\phi^{\prime})\simeq (V^{\prime},\phi\circ m\psi)\simeq (V,\phi)\text{ in } \tilde{K}^{-1}(M)_{\mathbb{Z}/m}\]
On the other hand, we have the following identities in $\Omega^{odd}(M)/d(\Omega^{even}(M))$: 
\begin{multline*}
\frac{1}{m}\operatorname{Tch}(\phi\circ m\psi\circ \phi^{-1})=\frac{1}{m}(\operatorname{Tch}(\phi)-\operatorname{Tch}(\phi^{\prime}))+\operatorname{Tch}(\psi)\\
=\frac{1}{m}\operatorname{Tch}(m\nabla_{v^{\prime}},\phi^{\prime,\ast}d)-\frac{1}{m}\operatorname{Tch}(m\nabla_{v},\phi^{\ast}d)+\omega-\omega^{\prime}.  
\end{multline*}
Thus, the two elements
\begin{align*}
&(\operatorname{Tch}(m\nabla_{v},\phi^{\ast}d),[(V,\phi)_{m}]);\\
&(\operatorname{Tch}(m\nabla_{v^{\prime}},\phi^{\prime,\ast}d),[(V^{\prime},\phi^{\prime})_{m}])
\end{align*}
are equivalent in $\mathcal{G}_{\operatorname{APS}}(M)$.

\end{proof}

In the following, we recall some basic properties of topological $K$-theory with coefficients in $\mathbb{C/Z}$. Let $\mathbf{K}_{t,\mathbb{C/Z}}\cong \mathbf{K}_{t}\wedge \mathbf{M}\mathbb{C/Z}$ be the $\Omega$-$CW$-prespectrum that represents complex topological $K$-theory with coefficients in $\mathbb{C/Z}$.
   
\begin{proposition}\label{property} 
\begin{enumerate}
\item Given $X$ a topological space, there is a long exact sequence:
\[...\rightarrow \tilde{K}^{-1}(X) \rightarrow  \tilde{H}^{odd}(X) \rightarrow \tilde{K}^{-1}(X)_{\mathbb{C/Z}}\rightarrow \tilde{K}^{0}(X)\rightarrow \tilde{H}^{even}(X)\rightarrow....\]

\item 
$\mathbf{K}_{t,\mathbb{C/Z}}$ is a module object over the ring object $(\mathbf{K}_{t},\nu,\iota)$ in the stable homotopy category $\operatorname{Ho}(\mathcal{P})$. More precisely, there is a map of $CW$-prespectra in $\mathcal{P}$: 
\[m:  \mathbf{K}_{t} \wedge   \mathbf{K}_{t,\mathbb{C/Z}}\rightarrow \mathbf{K}_{t,\mathbb{C/Z}}\] 
such that 
\begin{center}
\begin{tikzpicture}
\node(Lu) at (0,2) {$(\mathbf{K}_{t}\wedge\mathbf{K}_{t})\wedge \mathbf{K}_{t,\mathbb{C/Z}}$};
\node(Ll) at (0,0) {$\mathbf{K}_{t}\wedge(\mathbf{K}_{t}\wedge\mathbf{K}_{t,\mathbb{C/Z}})$}; 
\node(Ml) at (3.5,0) {$\mathbf{K}_{t}\wedge \mathbf{K}_{t,\mathbb{C/Z}}$};
\node(Ru) at (6,2) {$\mathbf{K}_{t}\wedge \mathbf{K}_{t,\mathbb{C/Z}}$};
\node(Rl) at (6,0) {$\mathbf{K}_{t,\mathbb{C/Z}}$};

\path[->, font=\scriptsize,>=angle 90]

(Lu) edge node [above]{$\nu\wedge id$}(Ru)  
(Lu) edge (Ll)
(Ll) edge node [above]{$id \wedge m$}(Ml)
(Ml) edge node [above]{$m$}(Rl) 
(Ru) edge node [right]{$m$}(Rl);
\end{tikzpicture}
\end{center}
and 
\begin{center}
\begin{tikzpicture}
\node(Lu) at (0,2) {$\mathbf{S}\wedge\mathbf{K}_{t,\mathbb{C/Z}}$};
\node(Ml) at (2,0) {$\mathbf{K}_{t,\mathbb{C/Z}}$}; 
\node(Ru) at (4,2) {$\mathbf{K}_{t}\wedge\mathbf{K}_{t,\mathbb{C/Z}}$};
 
\path[->, font=\scriptsize,>=angle 90]

(Lu) edge node [above]{$\iota\wedge id$}(Ru)  
(Lu) edge node [below]{$l$}(Ml) 
(Ru) edge node [right]{$m$}(Ml);

\end{tikzpicture}
\end{center}
commute in $\operatorname{Ho}(\mathcal{P})$, where the multiplication $\nu:\mathbf{K}_{t}\wedge \mathbf{K}_{t}$ is induced by the tensor product of vector bundles.
 
\end{enumerate}    
\end{proposition}\label{FactsofKCZ}
\begin{proof}
The long exact sequence in the first assertion follows from Proposition \ref{QadjunctionPresandspace} because the sequence 
\[\mathbf{K}_{t}\rightarrow \mathbf{K}_{t}\wedge\mathbf{M}\mathbb{C}\rightarrow \mathbf{K}_{t}\wedge \mathbf{M}\mathbb{C/Z}\]
is isomorphic to a cofiber sequence in $\mathcal{P}$. The second statement results from the ring structure on $\mathbf{K}_{t}$ (see \cite[Chapter $13$]{Sw}) and the definition of topological $K$-theory with coefficients in $\mathbb{C/Z}$. Namely, $\mathbf{K}_{t,\mathbb{C/Z}}\cong \mathbf{K}_{t}\wedge \mathbf{M}\mathbb{C/Z}$. 
\end{proof}
\begin{corollary}\label{Cor:ThomisoforKCZ}
\begin{enumerate}
\item We have the abelian group
 \[ 
 \tilde{K}^{-1}(S^{n})_{\mathbb{C/Z}}=\begin{cases} 
   \mathbb{C/Z} \hspace{2em}\text{when } n=\text{ odd };\\
   0 \hspace{3.5em}\text{when } n=\text{ even }.
  \end{cases}
\]                                                      
\item  
Given $M$ an odd dimensional spin manifold, there is a Thom isomorphism:
\[\tilde{K}^{-1}(M)_{\mathbb{C/Z}}\xrightarrow{\sim} \tilde{K}^{-1}(\tilde{\operatorname{Th}}(\nu_{M}))_{\mathbb{C/Z}},\]
where $\tilde{\operatorname{Th}}(\nu_{M})$ is the reduced Thom space of $\nu_{M}$, the normal bundle of $M$.
\end{enumerate} 
\end{corollary}

\begin{proof}
The first assertion follows from the long exact sequence in Proposition \ref{FactsofKCZ}. As for the second statement, we consider an embedding $M\hookrightarrow S^{N}$ with $N$ even and let $D$ be the associated Dirac operator of the spin structure on the normal bundle. Then the assertion follows from the fact that $D$ induces a Thom class $[D]$ in $\tilde{K}^{0}(\operatorname{Th}(\nu_{M}))$ (\cite[Appendix C]{LM}) and $\mathbf{K}_{t,\mathbb{C/Z}}$ is a module $CW$-prespectrum over $\mathbf{K}_{t}$, where $\operatorname{Th}(\nu_{M})$ is the unreduced Thom space of $\nu_{M}$. For more on Thom isomorphisms of module $CW$-prespectra, we refer to \cite[Section 5.1]{Ru}. 
\end{proof}

 
\section{Invariants of algebraic and relative algebraic $K$-theory}
In this section, we recall Atiyah, Patodi and Singer's construction of the homomorphism (see \cite[p.88-9]{APS3})
\[\bar{e}_{\operatorname{APS}}:[M,\operatorname{BGL}(\mathbb{C}^{\delta})]\rightarrow [M,F_{t,\mathbb{C/Z}}]\simeq \mathcal{G}_{\operatorname{APS}}(M),\]
where $M$ is a compact smooth manifold, and explain how it gives rise to an invariant of algebraic K-theory of the complex numbers, which can be identified with the $e$-invariant in \cite[p.940]{JW} via isomorphisms \eqref{eq:IsomPhi} and \eqref{eq:IsomPsi} in the proof of Theorem \ref{CZG}. How the $\tilde{\xi}$-invariant is used to compute the $e$-invariant is also explained. In the second subsection, we review the relation between the $e$-invariant, the Borel regulator (see \cite[Section $3$]{JW}) and the Chern-Simons classes. Lastly, employing the notion of locally unipotent bundles over a compact smooth manifold, we show the $e$-invariant detects the torsion part of the abelian group $[L,K_{a}\mathbb{C}]$, for every finite $CW$-complex $L$.  
 
 
\subsection{The $e$-invariant}
The following lemma comes in handy when we want to lift a map $X\rightarrow Y\in\mathpzc{Top}_{\ast}$ to $X^{+}\rightarrow Y$. Throughout the paper, $X\rightarrow X^{+}$ always means the plus map with respect to the perfect radical of $\pi_{1}(X)$.
 
\begin{lemma}\phantomsection\label{Lm:homoisoandweakhty}
Let $X,X^{\prime}$ and $Y$ be three pointed and connected $CW$-complexes and $f:X\rightarrow X^{\prime}$ a map that induces an isomorphism on integral homology. Suppose further $Y$ is homotopy equivalent to a loop space. Then $f_{\ast}:[X^{\prime},Y]_{w}\rightarrow [X,Y]_{w}$ is an isomorphism. 
\end{lemma}
\begin{proof}
As already indicated in \cite{CCMT}, the approach used in \cite[Lemma $1$]{CCMT} can be applied to the present case. Namely, we let $MY$ be the James construction on $Y$, which is homotopy equivalent to $\Omega (S^{1}\wedge Y)$ as $Y$ is connected. As James construction is universal in the category of loop spaces, there is a retraction $MY\rightarrow Y$. Combining these facts with the following adjunction
\[[X,S^{1}\wedge Y]_{w}\leftrightarrow [\Omega X, Y]_{w},\]
we obtain the commutative diagram below
\begin{center}
\begin{tikzpicture}
\node(Lu) at (0,3) {$[X^{\prime},Y]_{w}$};
\node(Lm) at (0,1.5) {$[X^{\prime},MY]_{w}$};
\node(Ll) at (0,0) {$[X^{\prime},Y]_{w}$}; 
\node(Ru) at (3,3) {$[X,Y]_{w}$};
\node(Rm) at (3,1.5) {$[X,MY]_{w}$};
\node(Rl) at (3,0) {$[X,Y]_{w}$};

\path[->, font=\scriptsize,>=angle 90] 
(Lu) edge node [above]{$f_{\ast}$}(Ru)
(Lu) edge (Lm)
(Lm) edge (Ll) 
(Ll) edge node [above]{$f_{\ast}$}(Rl) 
(Ru) edge (Rm)
(Rm) edge (Rl); 
\draw [->] (Lm) to node [above]{\scriptsize $f_{\ast}$} node [below]{$\sim$}(Rm);
\draw [transform canvas={xshift=-0.5ex},->] (Lu) to [out=-150,in=150] node [left] {\scriptsize $id$}(Ll);
\draw [transform canvas={xshift=0.5ex},->] (Ru) to [out=-30,in=30]   node [right] {\scriptsize $id$}(Rl); 
\end{tikzpicture}
\end{center}
The middle horizontal homomorphism is an isomorphism because $S^{1}\wedge f$ is a homotopy equivalence. By diagram chasing, we see the upper horizontal homomorphism is an injection and the bottom one a surjection. Hence $f_{\ast}$ must be a bijection.  
\end{proof}

We now recall Atiyah, Patodi and Singer's construction of $\bar{e}_{\operatorname{APS}}$. 
\begin{construction}\phantomsection\label{emap}
For any compact smooth manifold $M$, one can define an assignment
\begin{align}\label{APSconstr}
[M, \operatorname{BGL}(\mathbb{C}^{\delta})]&\rightarrow  \mathcal{G}_{\operatorname{APS}}(M)\\
(V,\nabla_{v})&\longmapsto ([\frac{1}{m}\operatorname{Tch}(m\nabla_{v},\phi^{\ast}d)][(V,\phi)_{m}]),
\end{align}
where $\phi$ is a bundle isomorphism from $mV$ to $\epsilon^{\star}$. One can check that it is independent of the choice of trivializations and is a homomorphism of commutative semigroups, for every compact smooth manifold $M$. Furthermore, it is functorial with respect to smooth maps. Since every finite $CW$-complex is homotopy equivalent to a compact smooth manifold and every continuous map of smooth manifolds is homotopic to a smooth map, one obtains a natural transformation of semigroup-valued functors on the category of finite $CW$-complexes      
\begin{equation}\label{Emapashomofsemigroups}
[-, \operatorname{BGL}(\mathbb{C}^{\delta})] \rightarrow  [-,F_{t,\mathbb{C/Z}}].
\end{equation} 
The representability theorem in \cite[Addendum 1.5]{Ad1} then implies that there exists a map unique up to weak homotopy,  
\begin{equation}\label{bareasamap}
\bar{e}: \operatorname{BGL}(\mathbb{C}^{\delta})\rightarrow F_{t,\mathbb{C/Z}},
\end{equation} 
where $F_{t,\mathbb{C/Z}}$ is the zero component of the $\Omega$-prespectrum $\mathbf{Fib}(\mathbf{ch})$ (see Theorem \ref{CZKtheory}). Since $\operatorname{BGL}(\mathbb{C}^{\delta})\rightarrow K_{a}\mathbb{C}$ induces an isomorphism on integral homology, Lemma \ref{Lm:homoisoandweakhty} tells us there exists a unique map (up to weak homotopy) 
\[e:  K_{a}\mathbb{C} \rightarrow F_{t,\mathbb{C/Z}}.\]  
\end{construction}

\begin{remark}\label{ExtendedbareAPS}
\textbf{1.} Given a compact smooth manifold $M$, we can also consider the following $\mathbb{Z}/m$-type $K$-theory model (see \cite[Section $5$]{APS2})
\begin{multline*}
\check{\mathcal{G}}_{\operatorname{APS},\mathbb{Z}/m}(M):=\{(V,W,\phi)\mid V,W \text{ are vector bundles over }M \text{ and }\\\phi:mV\xrightarrow{\sim} mW \text{ is a vector bundle isomorphism}\}/\sim,
\end{multline*} 
where $(V,W,\phi)$ and $(V^{\prime},W^{\prime},\phi^{\prime})$ are equivalent if and only if there exists (stably) an isomorphism of vector bundles $\psi$ from $V\oplus W^{\prime}$ to $W\oplus V^{\prime}$ such that $m\psi$ is homotopic to $\phi\oplus \phi^{\prime,-1}$ through isomorphisms of vector bundles stably. This gives a variant of the geometric model $\mathcal{G}_{\operatorname{APS}}(M)$ defined as the abelian group
\begin{multline*}
\check{\mathcal{G}}_{\operatorname{APS}}(M):= \{([\omega],[(V,W,\phi)]_{m})\mid \phi:mV\oplus\epsilon^{\star}\xrightarrow{\sim} mW\oplus\epsilon^{\star}\\
\text{ is an isomorphism of vector bundles}\}/\sim, 
\end{multline*}
where $([\omega],[(V,W,\phi)]_{m})\sim ([\omega^{\prime}],[(V^{\prime},W^{\prime},\phi^{\prime})]_{m})$ if and only if  there exists $l,l^{\prime}$ and $\psi:\epsilon^{\star}\rightarrow \epsilon^{\star}$ such that $lm=l^{\prime}m^{\prime}$, 
\[(V\oplus W^{\prime}\oplus \epsilon^{\star}, W\oplus V^{\prime}\oplus \epsilon^{\star},\phi\oplus \phi^{\prime,-1}\oplus \psi^{-1})\sim 0\in \check{\mathcal{G}}_{\operatorname{APS},\mathbb{Z}/m}(M),\]
and 
\[[\omega]-[\omega^{\prime}]=-\frac{1}{ml}\operatorname{Tch}(\psi^{\ast}d,d).\]
It is not difficult to see the inclusion 
\[\mathcal{G}_{\operatorname{APS}}(M)\rightarrow \check{\mathcal{G}}_{\operatorname{APS}}(M)\]
is an isomorphism of abelian groups. This variant of $\mathcal{G}_{\operatorname{APS}}(M)$ helps us extend the domain of $\bar{e}_{\operatorname{APS}}$. Let $\tilde{K}(M,\mathbb{C})$ be the abelian group of virtual flat vector bundles of dimension zero. Namely, the kernel of the dimension map 
\[\operatorname{dim}:K(M,\mathbb{C})\rightarrow \mathbb{Z},\] 
where $K(M,\mathbb{C})$ is the Grothendieck group of flat vector bundles over $M$. Then one can consider the following homomorphism, denoted by $\bar{e}_{\operatorname{APS}}$ also,
\begin{align*}
\bar{e}_{\operatorname{APS}}:\tilde{K}(M,\mathbb{C})&\rightarrow [M,F_{t,\mathbb{C/Z}}]\\
 (V,W)&\mapsto ([\frac{1}{m}\operatorname{Tch}(m\nabla_{v},\psi^{\ast}m\nabla_{w})],[(V,W,\psi)_{m}])
\end{align*}
where $\psi$ is any chosen isomorphism of vector bundles from $mV$ to $mW$. It is not difficult to see it fits into the following commutative diagram:
\begin{center}
\begin{tikzpicture}
\node (Lu) at (0,3) {$[M,\operatorname{BGL}(\mathbb{C}^{\delta})]$} ;
\node (Lm) at (0,2) {$\tilde{K}(M,\mathbb{C})$};
\node (Ll) at (0,1) {$[M,K_{a}\mathbb{C}]$};
\node (Rm) at (5,2) {$[M,F_{t,\mathbb{C/Z}}]$};

\draw [->] (Lu) to (Lm);
\draw [->] (Lm) to (Ll);
\draw [->](Lu) to [out=0,in=150]node [above]{\scriptsize $\bar{e}_{\operatorname{APS}}$}(Rm);
\draw [->](Lm) to  node [above]{\scriptsize $\bar{e}_{\operatorname{APS}}$}(Rm);
\draw [->](Ll) to [out=0,in=-150]node [above]{\scriptsize $e_{\ast}$}(Rm);
\end{tikzpicture}
\end{center}
 
\noindent
\textbf{2.} In view of isomorphism \eqref{eq:IsomPhi}:
\begin{align*}
\Phi : \mathcal{G}_{\operatorname{APS}}(M) &\rightarrow \mathcal{G}_{\operatorname{JW}}(M)\\
(\omega, [(V,\phi)_{m}])&\mapsto [(V,\nabla_{v},\frac{1}{m}\operatorname{Tch}(m\triangledown_{v},\phi^{*}d)-\omega)], 
\end{align*}
where $\nabla_{v}$ is any connection on $V$, we see assignment \eqref{APSconstr} is equivalent to the homomorphism given in  \cite[p.940]{JW}:
\begin{align*}
[M, \operatorname{BGL}(\mathbb{C}^{\delta})]&\rightarrow  \mathcal{G}_{\operatorname{JW}}(M)\\
(V,\nabla_{v})&\longmapsto (V,\nabla_{v}, 0).
\end{align*} 
\end{remark} 
 
\begin{definition}\phantomsection\label{einvariantdef}
Given any pointed finite $CW$-complex $L$, we call the induced homomorphism 
\[e_{\ast}: [L,K_{a}\mathbb{C}] \rightarrow [L,F_{t,\mathbb{C/Z}}]\]
the $e$-invariant. 
\end{definition}

\subsection{The $\tilde{\xi}$-invariant of spin manifolds}
Recall that, given a spin manifold $M$, we have the Thom isomorphism (see Corollary \ref{Cor:ThomisoforKCZ})  
\[T_{D}:[M,F_{t,\mathbb{C/Z}}] \xrightarrow{\otimes [D]}[\tilde{\operatorname{Th}}(\nu_{M}),F_{t,\mathbb{C/Z}}],\] 
where $D$ is the Dirac operator associated to the spin structure on $\nu_{M}$ the normal bundle of $M$, a tubular neighborhood of the embedding $M\hookrightarrow S^{N}$ with $N$ odd. Notice also the collapsing map $S^{N}\rightarrow \tilde{\operatorname{Th}}(\nu_{M})$ induces a homomorphism
\[c^{\ast}:[M,F_{t,\mathbb{C/Z}}]\rightarrow [S^{N},F_{t,\mathbb{C/Z}}].\]
The following is a special case of Atiyah, Patodi and Singer's index theorem for flat vector bundles.  
\begin{corollary}\phantomsection\label{Xiandeidentified}
Given $M$ a spin manifold and a representation (or flat vector bundle)  
\[\rho:\pi_{1}(M)\rightarrow \operatorname{GL}_{N}(\mathbb{C}),\]
we have the equality
\[c^{\ast}\circ T_{D}\circ \bar{e}_{\operatorname{APS}}(M,\rho)=\tilde{\xi}(\rho,M)\text{ in } \mathbb{C/Z}.\] 
\end{corollary} 
\begin{proof}
This is essentially the index theorem for flat vector bundles in \cite[5.3]{APS3} (see the remarks in \cite[p.87,p. 89-90]{APS3} for the non-unitary case).
\end{proof} 
The next lemma is crucial when we want to use the previous corollary to compute the $e$-invariant. We first observe the degree one map of a homology $n$-sphere onto a $n$-sphere is a plus map, and hence a representation or flat vector bundle on a homology $n$-sphere $\Sigma^{n}$ 
\[f:\Sigma^{n}\rightarrow \operatorname{BGL}(\mathbb{C}^{\delta})\]
induces a map, unique up to homotopy,
\[f^{+}:S^{n}\rightarrow K_{a}\mathbb{C}.\]
Meaning there is a homomorphism (of semigroups)
\begin{equation}\label{Eq:flatonhomosphere}
\alpha:[\Sigma^{n},\operatorname{BGL}(\mathbb{C}^{\delta})]\rightarrow [S^{n},K_{a}\mathbb{C}].
\end{equation}
\begin{lemma}\label{Differentrepthesamevalue}
Let $\Sigma_{1}^{n}$ and $\Sigma_{2}^{n}$ be two homology $n$-spheres and
\begin{align*}
\rho_{1}:\pi_{1}(\Sigma_{1}^{n})&\rightarrow \operatorname{GL}_{N}(\mathbb{C}),\\
\rho_{2}:\pi_{1}(\Sigma_{2}^{n})&\rightarrow \operatorname{GL}_{N}(\mathbb{C}),
\end{align*}
be representations of their fundamental groups. Suppose also $(\Sigma_{1}^{n},\rho_{1})$ and $(\Sigma_{2}^{n},\rho_{2})$ have the same image $x$ in $[S^{n},K_{a}\mathbb{C}]$ under homomorphism \eqref{Eq:flatonhomosphere}.  
Then $\tilde{\xi}(\Sigma_{1}^{n},\rho_{1})=\tilde{\xi}(\Sigma_{2}^{n},\rho_{2})$ in $\mathbb{C/Z}$. 
\end{lemma}
\begin{proof}
Recall that the $\tilde{\xi}$-invariant is a spin cobordism invariant and any $H_{\ast}$-cobordism is a spin cobordism. For $n\geq 5$, we know that $(\Sigma_{1}^{n},\rho_{1})$ and $(\Sigma_{2}^{n},\rho_{2})$ bound a $H_{\ast}$-cobordism $(H,F)$ by Hausmann and Vogel's theory (see \cite[Corollary $4.2$]{HV}). In particular, $(\Sigma_{1}^{n},\rho_{1})$ and $(\Sigma_{2}^{n},\rho_{2})$ are spin cobordant via $(H,F)$. For $n=3$, $\Sigma^{n}_{1}$ and $\Sigma_{2}^{n}$ are spin cobordant as the third spin cobordism group $\Omega^{spin}_{3}(\ast)$ is trivial. Let $W$ be the spin cobordism and $s:W\rightarrow \operatorname{Th}(\operatorname{Spin})$ the associated map into the Thom space of the canonical vector bundle over $\operatorname{BSpin}$, given by the canonical map $\operatorname{BSpin}\rightarrow \operatorname{BO}$. Observe also that there is a plus map
\[\pi:W\rightarrow S^{n}\times I\]  
that restricts to the plus maps $\Sigma_{1}^{n}\rightarrow S^{n}$
and $\Sigma_{2}^{n}\rightarrow S^{n}$ on the boundary of $W$. Suppose $x\in[S^{n},K_{a}\mathbb{C}]$
is represented by the map $f:S^{n}\rightarrow \operatorname{BGL}^{+}(\mathbb{C})$. Then the following composition 
\begin{multline*}
S^{k+n}\wedge I_{+}\rightarrow \operatorname{Th}(\nu_{W})\wedge W_{+}\xrightarrow{id\wedge\pi} \operatorname{Th}(\nu_{W})\wedge (S^{n}\times I)_{+}\\
\xrightarrow{id\wedge p}\operatorname{Th}(\nu_{W})\wedge S^{n}_{+}\xrightarrow{s\wedge f} \operatorname{Th}(\operatorname{Spin})\wedge \operatorname{BGL}^{+}(\mathbb{C}^{\delta})
\end{multline*}
implies
$(\Sigma_{1}^{n},\rho_{1})$ and $(\Sigma_{2}^{n},\rho_{2})$
are spin cobordant in $\Omega^{spin}_{3}(\operatorname{BGL}^{+}(\mathbb{C}^{\delta}))$ and hence spin cobordant in $\Omega^{spin}_{3}(\operatorname{BGL}(\mathbb{C}^{\delta}))$, where $\operatorname{Th}(\nu_{W})$ is the Thom space of the normal bundle of $W$ and $p:S^{n}\times I\rightarrow S^{n}$ is the obvious projection. We have used the fact that there are equivalences  
\begin{multline*}
\Sigma^{\infty}\operatorname{Th}(\operatorname{Spin})\wedge \operatorname{BGL}(\mathbb{C}^{\delta})\simeq \operatorname{Th}(\operatorname{Spin})\wedge \Sigma^{\infty}\operatorname{BGL}(\mathbb{C}^{\delta})\\
\xrightarrow{\sim} \operatorname{Th}(\operatorname{Spin})\wedge \Sigma^{\infty}\operatorname{BGL}^{+}(\mathbb{C}^{\delta})\simeq \Sigma^{\infty}\operatorname{Th}(\operatorname{Spin})\wedge \operatorname{BGL}^{+}(\mathbb{C}^{\delta}) 
\end{multline*}
in $\operatorname{Ho}(\mathcal{P})$.
\end{proof}
Now, observe the following commutative diagram 
\begin{center}
\begin{tikzpicture}
\node (As) at (0,0) {$[S^{n},K_{a}\mathbb{C}]$};
\node (Ts) at (3,0) {$[S^{n},F_{t,\mathbb{C/Z}}]$};
\node (CZ) at (7,0) {$\mathbb{C/Z}$};
 
\node (Ah) at (0,-1.5) {$[\Sigma^{n},\operatorname{BGL}(\mathbb{C}^{\delta})]$};
\node (Th) at (3,-1.5) {$[\Sigma^{n},F_{t,\mathbb{C/Z}}]$};
\node (Tth)at (6,-1.5) {$[\tilde{\operatorname{Th}}(\nu_{\Sigma^{n}}),F_{t,\mathbb{C/Z}}]$};
\node (S2) at (9,-1.5) {$[S^{N},F_{t,\mathbb{C/Z}}]$};

\path[->,font=\scriptsize, >= angle 90]
(As)  edge node [above] {$e_{\ast}$}(Ts)
 
(Ah)  edge node [above] {$\bar{e}_{\ast}$}(Th)

(Ts)  edge node [above] {$\sim$}(CZ)
(Ah)  edge node [right] {$+$}(As)
(Th)  edge node [right] {$+$}(Ts);

\draw[->] (Th)  -- node [above] {\tiny $T_{D}$} node[below]{$\sim$} ++(Tth); 
\draw[->] (Tth)  -- node [above] {\tiny $c^{\ast}$} node[below]{$\sim$} ++(S2);

\end{tikzpicture}
\end{center}
where $+$ stands for the homomorphism induced by the universal property of the plus map $\Sigma^{n}\rightarrow S^{n}$. Notice also, in this case, $T_{D}$ and $c^{\ast}$ are isomorphisms. Hence, in view of Lemma \ref{Differentrepthesamevalue} and Corollary \ref{Xiandeidentified}, we have the following theorem:
\begin{theorem}\phantomsection\label{Xiandeidentified1}
Given an element $x\in [S^{n},K_{a}\mathbb{C}]$, then $e_{\ast}(x)$ can be computed by the $\tilde{\xi}$-invariant of any representation of a homology sphere $\rho:\pi_{1}(\Sigma^{n})\rightarrow \operatorname{GL}_{N}(\mathbb{C})$ whose image under homomorphism \eqref{Eq:flatonhomosphere} is $x$.
\end{theorem}

\subsection{Relative $K$-theory, the Chern-Simons classes and the Borel regulator}\label{RelKtheory}

Here we explain how the Chern-Simons classes give rise to a map from the relative $K$-theory space to the product of Eilenberg-Maclane spaces of odd degrees: 
\[\operatorname{ch}^{rel}:K^{rel}\mathbb{C}\rightarrow H^{odd}\mathbb{C}:=\prod_{\mathclap{i\text{ odd }}}K(\mathbb{C},i).\]
The relation of $\operatorname{ch}^{rel}$ with the Borel regulator and the $e$-invariant is also reviewed.  

\begin{definition}
Relative $K$-theory of the complex numbers is the generalized cohomology theory represented by the homotopy fiber of the canonical map from the $0$-connective algebraic $K$-theory prespectrum of the complex numbers to the $0$-connective complex topological $K$-theory prespectrum 
\[\mathbf{K}_{a}\mathbb{C}\rightarrow \mathbf{K}_{t}.\]
Its infinite loop space is denoted by $K^{rel}\mathbb{C}$.
\end{definition}

In order to find an explicit model of $K^{rel}\mathbb{C}$, we recall the following lemma \cite[Theorem $1.1(a)$]{Be2}.   

\begin{lemma}\label{Lm:Plusconstrandfib}
Denote the perfect radical of the group $G$ by $\mathfrak{P}G$. Then, given a homotopy fiber sequence 
\[F\rightarrow E\rightarrow B,\] 
with $\mathfrak{P}\pi_{1}B$ is trivial,   
the sequence
\[F^{+}\rightarrow E^{+}\rightarrow B\]
is also a homotopy fiber sequence.
\end{lemma} 
\begin{proof}
See \cite[p.150-151]{Be2}. Berrick considers $F^{\prime}$ the homotopy fiber of $E^{+}\rightarrow B$ and shows the induced map $F\rightarrow F^{\prime}$ is a plus map.
\end{proof}
Denote the canonical (topology-changing) map by   
\[\iota:\operatorname{BGL}(\mathbb{C}^{\delta})\rightarrow K_{t}\]
and consider the homotopy fiber sequence
\[\operatorname{BGL}(\mathbb{C}^{\delta})\times_{\iota}K_{t}^{I_{\ast}}\rightarrow \operatorname{BGL}(\mathbb{C}^{\delta}) \rightarrow K_{t},\]
where $K_{t}^{I_{\ast}}$ is the function space of pointed maps from $(I,0)$, the interval $[0,1]$ with base point $0$, to $(K_{t},\ast)$, and we define 
\[\operatorname{BGL}(\mathbb{C}^{\delta})\times_{\iota}K_{t}^{I_{\ast}}:=\{(x,\gamma(t))\in \operatorname{BGL}(\mathbb{C}^{\delta}))\times K_{t}^{I_{\ast}} \text{ with }\gamma(1)=\iota(x) \}.\]
Applying Berrick's theorem to this case, we obtain the following map of homotopy fiber sequences:  

\begin{center}
\begin{tikzpicture}
\node(Lu) at (0,2) {$\operatorname{BGL}(\mathbb{C}^{\delta})\times_{\iota}K_{t}^{I_{\ast}}$};
\node(Ll) at (0,0) {$\operatorname{BGL}^{+}(\mathbb{C}^{\delta})\times_{\iota}K_{t}^{I_{\ast}}$};
\node(Mu) at (4,2) {$\operatorname{BGL}(\mathbb{C}^{\delta})$};
\node(Ml) at (4,0) {$\operatorname{BGL}^{+}(\mathbb{C}^{\delta})$};
\node(Ru) at (8,2) {$K_{t}$};
\node(Rl) at (8,0) {$K_{t}$};

\path[->, font=\scriptsize,>=angle 90]

(Ll) edge (Ml)
(Ml) edge (Rl)
(Lu) edge (Mu)
(Mu) edge (Ru)  
(Lu) edge (Ll)
(Mu) edge (Ml)  
(Ru) edge node [right]{$\mid$}(Rl);

\end{tikzpicture}
\end{center}
where all vertical maps are the plus maps. In particular, we see $K^{rel}\mathbb{C}:=\operatorname{BGL}^{+}(\mathbb{C})\times K_{t}^{I_{\ast}}$ is the infinite loop space of relative $K$-theory.

Now, the space $\operatorname{BGL}(\mathbb{C}^{\delta})\times_{\iota} K_{t}^{I_{\ast}}$ classifies flat vector bundles with a trivialization, and given a flat vector bundle with a trivialization over a compact smooth manifold $M$, $(V,\nabla_{v},\phi:V\xrightarrow{\sim} \epsilon^{\star})$ say, one can defined the associated Chern-Simons class: 
\[\tilde{\operatorname{ch}}(\nabla_{v},\phi^{\ast}d):=\int_{t}ch(\nabla_{v,t}),\] 
where $\nabla_{v,t}$ is the connection on $V\times I$ given by  
\[\nabla_{v,t}:=t\nabla_{v}+(1-t)\phi^{\ast}d.\]
It determines a cohomological class as $\nabla_{v}$ is flat.
This assignment induces a homomorphism of semigroups
\[[M,\operatorname{BGL}(\mathbb{C}^{\delta})\times_{\iota}K_{t}^{I_{\ast}}]\rightarrow [M,H^{odd}\mathbb{C}], \]
and since it is functorial with respect to smooth maps, by Adams' representability theorem (see \cite{Ad}), we obtain a map 
\[\bar{\operatorname{ch}}^{rel}:\operatorname{BGL}(\mathbb{C}^{\delta})\times_{\iota}K_{t}^{I_{\ast}} \rightarrow H^{odd}\mathbb{C}.\] 
Furthermore, applying Lemma \ref{Lm:homoisoandweakhty}, we get a map from the relative $K$-theory 
\[\operatorname{ch}^{rel}: K^{rel}\mathbb{C}\rightarrow H^{odd}\mathbb{C},\] 
which is unique up to homotopy as $H^{odd}\mathbb{C}$ is a rational infinite loop space (see Lemma \ref{Weakhtyandhty}).
 
We now recall some properties of the map $\operatorname{ch}^{rel}$ (see \cite[Theorem $3.1$ and $3.2$]{JW} and compare with \cite[Remark $7.19$]{Ka1}): 
\begin{theorem}
The following diagram commutes up to weak homotopy

\begin{equation}\label{chrelecommutativediag}
\begin{tikzpicture}[baseline=(current  bounding  box.center)] 
 \node (OKt) at (0,6) {$\Omega K_{t}$};
 \node (Kr)  at (0,4) {$K^{rel}\mathbb{C}$};
 \node (Ka)   at (0,2) {$K_{a}\mathbb{C}$};
 \node (Kt)   at (0,0) {$K_{t}$};

 \node (OKt1) at (4,6) {$\Omega K_{t}$};
 \node (Pro2)  at (4,4) {$H^{odd}\mathbb{C}$};
 \node (Pro1)  at (4,2) {$F_{t,\mathbb{C/Z}}$};
 \node (Kt1)   at (4,0) {$K_{t}$};

 \path [->, font=\scriptsize, >=angle 90]
  (Kr)  edge node [above]{$\operatorname{ch}^{rel}$} (Pro2)
  (Ka)  edge node [above]{$e$} (Pro1) 
  
  (OKt)  edge node [right]{$i$}(Kr)
  (Kr)   edge node [right]{$\pi$}(Ka)
  (Ka)   edge node [right]{$\iota$}(Kt)
  (OKt1) edge node [right]{$\operatorname{ch}$}(Pro2) 
  (Pro2) edge (Pro1)
  (Pro1) edge (Kt1); 
  \draw [double equal sign distance](OKt) to (OKt1);
  \draw [double equal sign distance] (Kt) to (Kt1); 
\end{tikzpicture} 
\end{equation}
\end{theorem}
\begin{proof}
In view of Lemma \ref{Lm:homoisoandweakhty}, it suffices to show the following commutes up to weak homotopy
\begin{equation}\label{barchrelecommutativedia} 
\begin{tikzpicture}[baseline=(current  bounding  box.center)] 
 \node (OKt) at (0,6) {$\Omega K_{t}$};
 \node (Kr)  at (0,4) {$\operatorname{BGL}(\mathbb{C}^{\delta})\times_{\iota}K_{t}^{I_{\ast}}$};
 \node (Ka)   at (0,2) {$\operatorname{BGL}(\mathbb{C}^{\delta})$};
 \node (Kt)   at (0,0) {$K_{t}$};

 \node (OKt1) at (4,6) {$\Omega K_{t}$};
 \node (Pro2)  at (4,4) {$H^{odd}\mathbb{C}$};
 \node (Pro1)  at (4,2) {$F_{t,\mathbb{C/Z}}$};
 \node (Kt1)   at (4,0) {$K_{t}$};

 \path [->, font=\scriptsize, >=angle 90]
  (Kr)  edge node [above]{$\operatorname{ch}^{rel}$} (Pro2)
  (Ka)  edge node [above]{$e$} (Pro1) 
  
  (OKt)  edge node [right]{$i$}(Kr)
  (Kr)   edge node [right]{$\pi$}(Ka)
  (Ka)   edge node [right]{$\iota$}(Kt)
  (OKt1) edge node [right]{$\operatorname{ch}$}(Pro2) 
  (Pro2) edge (Pro1)
  (Pro1) edge (Kt1); 
  \draw [double equal sign distance](OKt) to (OKt1);
  \draw [double equal sign distance] (Kt) to (Kt1); 
\end{tikzpicture} 
\end{equation}

The commutativity of the upper square is clear as, given an isomorphism of the trivial bundle $\phi:\epsilon^{\star}\rightarrow \epsilon^{\star}$, we have $\tilde{\operatorname{ch}}(d,\phi^{\ast}d)$ computes the odd Chern character. Now let $M$ be a compact smooth manifold. The commutativity of the lower square follows from the fact that the homomorphism
\[
[M,F_{t,\mathbb{C/Z}}]\simeq \mathcal{G}_{\operatorname{APS}}(M)\rightarrow [M,K_{t}]
\]
can be realized by the assignment 
\[([\omega],[(V,\phi)_{m}]) \mapsto [V].\]
It is also straightforward to see the middle square commutes. One only needs to note the homomorphisms  
\[
[M,H^{odd}\mathbb{C}]\rightarrow [M,F_{t,\mathbb{C/Z}}]\simeq \mathcal{G}_{\operatorname{APS}}(M) 
\]
and 
\[
[M,\operatorname{BGL}(\mathbb{C}^{\delta})\times_{\iota}K_{t}]\rightarrow [M,\operatorname{BGL}(\mathbb{C}^{\delta})]
\] 
can be realized by 
\[[\omega]\mapsto [[\omega],[(\epsilon^{\star},id)_{1}]\] 
and  
\[(V,\nabla_{v};\phi:V\xrightarrow{\sim}\epsilon^{\star})\mapsto (V,\nabla_{v}),\]
respectively.
\end{proof}
\begin{remark}
In view of isomorphisms \eqref{eq:IsomPhi} and \eqref{eq:IsomPsi}, the proof of commutative diagram \eqref{barchrelecommutativedia} is equivalent to the one given in \cite[p.945]{JW}.  
\end{remark}

The following describes the relation between the Borel regulator $\operatorname{Bo}$ and the map $e$:

\begin{theorem}
The composition
\begin{equation}\label{eq:CompositionchJe}
K_{a}\mathbb{C}\rightarrow F_{t,\mathbb{C/Z}}\xrightarrow{\operatorname{Im}} \Omega K_{\mathbb{R}} \xrightarrow{\operatorname{ch}_{\otimes\mathbb{R}}} H^{odd}\mathbb{R}
\end{equation} 
is homotopic to the infinite loop space map
\begin{equation}\label{eq:Borelregulator}
\operatorname{Bo}:=\sum_{k}(-1)^{k-1}(\frac{1}{2\pi})^{k}\frac{(k-1)!}{(2k-1)!}b_{2k-1}:K_{a}\mathbb{C}\rightarrow H^{odd}\mathbb{R}, 
\end{equation}
where $b_{2k-1}$ is the Borel class of degree $2k-1$, $H^{odd}\mathbb{R}:=\prod\limits_{\mathclap{i \text{ odd }}}K(\mathbb{R},i)$, $\operatorname{Im}:F_{t,\mathbb{C/Z}}\rightarrow \Omega K_{t,\mathbb{R}}$ is induced by the homomorphism
\begin{align*}
\mathbb{C/Z}&\rightarrow \mathbb{R}\\
a+bi&\mapsto b
\end{align*}
and $\operatorname{ch}_{\otimes \mathbb{R}}:\Omega K_{t,\mathbb
R}\rightarrow H^{odd}\mathbb{R}$ is given by the Chern character. 
\end{theorem}
\begin{proof}
The construction of $b_{2k-1}$ is given in \cite[p.943-4]{JW}, and it has been proved in \cite[Theorem $3.1$]{JW} maps \eqref{eq:CompositionchJe} and \eqref{eq:Borelregulator} are weakly homotopic. Since the target is a rational infinite loop space, they are actually homotopic.   
\end{proof}

\subsection{The $e$-invariant and torsion subgroups} 
The $e$-invariant is highly non-trivial. In fact, we have the following theorem (compare with \cite[Theorem $7.20$]{Ka1} and \cite[Corollary $2.4$]{JW}):
\begin{theorem}\phantomsection\label{IsomonTor}  
The $e$-invariant restricts to an isomorphism on the torsion subgroup of $[L,K_{a}\mathbb{C}]$:  
\[e_{\ast}\vert_{\operatorname{Tor}}:\operatorname{Tor}([L,K_{a}\mathbb{C}])\xrightarrow{\sim} \operatorname{Tor}([L,F_{t,\mathbb{C/Z}}]),\]
for every finite $CW$-complex $L$.
\end{theorem}
We shall identify the homomorphism $e_{\ast}\vert_{\operatorname{Tor}}$ with Suslin's isomorphism $\operatorname{Su}_{\ast}$ (see Theorem \ref{Suslinthm1} and \cite[Corollary $4.6$]{Su}) by proving the lemma below\footnote{This idea is due to Karoubi \cite[p.115]{Ka1}.}. Then the theorem follows quickly from the argument in \cite[Proof of Corollary $2.4$]{JW} (see also \cite[p.115]{Ka1}). 

\begin{lemma}\phantomsection\label{CommDiaforeTor}\footnote{It is essentially Theorem $2.3$ in \cite{JW}. The argument presented there is however not clear to us. We are not sure how to see $\operatorname{Rep}(\pi_{1}(X),\operatorname{GL}(\mathbb{C}))_{\mathbb{Z}/m}$ generates $K^{-1}_{alg}(X,\mathbb{C})_{\mathbb{Z}/m}$ as claimed in \cite[p.941]{JW} since $\operatorname{Rep}(\pi_{1}(X),\operatorname{GL}(\mathbb{C}))_{\mathbb{Z}/m}$ should be zero all the time.}
Let $F_{a,m}$ and $F_{t,m}$ be the infinite loops spaces defined in Lemma \ref{rmk:OtherformulationofSu}. Then, given a finite $CW$-complex $L$, the following diagram of abelian groups commutes\footnote{The same commutative diagram but with $L=S^{n}$ has been claimed without proof in \cite[p.115]{Ka1}.}:   
\begin{center} 
\begin{equation}\label{Diag:erestrictedtotor} 
\begin{tikzpicture}[baseline=(current bounding box.center)] 
\node (Lu) at (0,2){$[L,F_{a,\mathbb{Z}/m}]$};
\node (Ll) at (0,0){$[L,K_{a}\mathbb{C}]$};
\node (Ru) at (4,2){$[L,F_{t,\mathbb{Z}/m}]$};
\node (Rl) at (4,0){$[L,F_{t,\mathbb{C/Z}}]$};

\path [->, font=\scriptsize, >=angle 90]
 
 (Lu) edge node[right]{$b$}(Ll)
 (Ll) edge node[above]{$e_{\ast}$}(Rl)
 (Lu) edge node[above]{$\operatorname{Su}_{\ast}$}(Ru)
 (Ru) edge node[right]{$j_{\ast}$}(Rl);

\end{tikzpicture}
\end{equation} 
\end{center}
where $j$ is induced by the canonical inclusions
\[\mathbb{Z}/m\hookrightarrow \mathbb{Q/Z}\rightarrow \mathbb{C/Z},\]
$\operatorname{Su}$ is Suslin's map (Remark \ref{rmk:OtherformulationofSu}), 
and $b$ is the Bockstein homomorphism.
\end{lemma} 
\begin{proof}
Denote the homotopy fiber of the composition
\[\operatorname{BGL}(\mathbb{C}^{\delta})\xrightarrow{m}\operatorname{BGL}(\mathbb{C}^{\delta})\rightarrow K_{a}\mathbb{C}\]
by $\operatorname{BGL}(\mathbb{C}^{\delta})\times_{m}K_{a}^{I_{\ast}}$, where the map $m$ is induced by taking the direct sum of $m$ copies of a flat vector bundle. Then, by Lemma \ref{Lm:Plusconstrandfib}, the canonical map
\[\operatorname{BGL}(\mathbb{C}^{\delta})\times_{m}K_{a}^{I_{\ast}}\rightarrow F_{a,\mathbb{Z}/m}\]
is a plus map, and hence, in view of Lemma \ref{Lm:homoisoandweakhty} , it suffices to show the diagram below commutes up to weak homotopy
\begin{center}
\begin{equation}\label{diag:Inpfisoontor}
\begin{tikzpicture}[baseline=(current bounding box.center)] 
\node (Lu) at (0,3){$\operatorname{BGL}(\mathbb{C}^{\delta})\times_{m}K_{a}^{I_{\ast}}$};
\node (Lm) at (0,1.5){$\operatorname{BGL}(\mathbb{C}^{\delta})$};
\node (Ll) at (0,0){$K_{a}\mathbb{C}$};
\node (Mu) at (3,3){$F_{a,\mathbb{Z}/m}$};
\node (Ru) at (6,3){$F_{t,\mathbb{Z}/m}$};
\node (Rl) at (6,0){$F_{t,\mathbb{C/Z}}$};

\path [->, font=\scriptsize, >=angle 90]
 
 (Lu) edge node [right]{$\pi$}(Lm)
 (Lm) edge node [right]{$+$}(Ll)
 (Ll) edge node [above]{$e$}(Rl)
 (Lu) edge node [above]{$+$}(Mu)
 (Mu) edge node [above]{$\operatorname{Su}$}(Ru)
 (Ru) edge node [right]{$j$}(Rl);

\end{tikzpicture} 
\end{equation}
\end{center}
where $+$'s stand for plus maps.
Now it is clear that, given a flat vector bundle $(V,\nabla_{v})$ over a compact smooth manifold $M$ and a homotopy $M\times I\rightarrow K_{a}\mathbb{C}$ that restricts to $(mV,\nabla_{v})$ and $(\epsilon^{\star},d)$ on $\partial (M\times I)$, we have 
\begin{align*}
(j\circ \operatorname{Su}\circ +)_{\ast}(V,\nabla_{v};\phi:mV\rightarrow \epsilon^{\star})&=(0,[V,\phi]_{m});\\
(e\circ +\circ \pi)_{\ast} (V,\nabla_{v};\phi:mV\rightarrow \epsilon^{\star})&=([\frac{1}{m}\tilde{\operatorname{ch}}(m\nabla_{v},\phi^{\ast}d)],[V,\phi]_{m}),
\end{align*}
where $\phi$ is induced by the given homotopy $M\times I\rightarrow K_{a}\mathbb{C}$.
Hence, if we can show 
\begin{equation}\label{eq:vanishingoftildech1}
[\tilde{\operatorname{ch}}(m\nabla_{v},\phi^{\ast}d)]=0,
\end{equation}
then we are done. 

To see this, we first observe that, given a connection $\nabla_{u}$ on a vector bundle $U$ over $M\times I$ with $(U,\nabla_{u})\vert_{M\times\{1\}}$ flat and $(U,\nabla_{u})\vert_{M\times\{0\}}=(\epsilon^{\star},d)$, we have 
\[\tilde{\operatorname{ch}}(\nabla_{u}):=\int_{t\in I} \operatorname{ch}(\nabla_{u})=\tilde{\operatorname{ch}}(\nabla_{u}\vert_{M\times\{1\}},\phi^{\ast}d)\in \Omega^{odd}(M)/\operatorname{Im}(d),\] 
where $\phi:U\vert_{M\times\{1\}}\xrightarrow{\sim} \epsilon^{\star}$ is an isomorphism induced by the vector bundle $U$ over $M\times I$, and $\operatorname{Im}(d)$ is the image of the differential $d:\Omega^{even}(M)\rightarrow \Omega^{odd}(M)$.

In view of this observation, it is sufficient to prove the following claim: Given a smooth manifold $W$ with boundary $\partial W$ and the commutative diagram
\begin{center}
\begin{tikzpicture}
\node(Lu) at (0,2) {$\partial W$};
\node(Ll) at (0,0) {$W$}; 
\node(Ru) at (2,2) {$\operatorname{BGL}(\mathbb{C}^{\delta})$};
\node(Rl) at (2,0) {$K_{a}\mathbb{C}$};

\path[->, font=\scriptsize,>=angle 90]

(Lu) edge node [above]{$f$}(Ru)  
(Ll) edge node [above]{$f^{\prime}$}(Rl) 
(Ru) edge node [right]{$+$}(Rl);
\draw[right hook ->] (Lu) to (Ll);
\end{tikzpicture}
\end{center}
then there exists a locally unipotent bundle $(V,\nabla_{v})$ on $W$ such that it restricts to the flat vector bundle on $\partial W$. We shall explain what we mean by locally unipotent bundles. A connection $\nabla_{E}$ on a vector bundle $E$ over $W$ is locally unipotent if and only if, for every point $p$, there exists an open neighborhood $\mathfrak{U}_{p}$ and a filtration 
\[E\vert_{\mathfrak{U}_{p}}=F_{p}^{n}\supsetneq F_{p}^{n-1}\supsetneq...\supsetneq F_{p}^{0}=\{0\}\]
such that 
\[\nabla^{2}(\Gamma(F_{p}^{i}))\subset \Omega^{2}(W;F_{p}^{i-1}),\] 
where, given a vector bundle $F$ over $W$, $\Gamma(F)$ denotes the space of sections, $\Omega^{2}(W;F_{p}^{i-1})$ is defined to be $\Gamma(\Lambda^{2}TW\otimes_{\mathbb{C}} F_{p}^{i-1})$, and $\nabla^{2}$ is the curvature $2$-form. A vector bundle equipped with a locally unipotent connection is called a locally unipotent vector bundle (see \cite[Sections $1$ and $2$]{Kr} and \cite[Remark 2.16]{Sch} for more details).

Return to our case, we let $W=M\times I$, $f^{\prime}$ be the homotopy, and $f$ represent the flat vector bundles $(mV,m\nabla_{v})$ and $(\epsilon^{\star},d)$ on $M\times\{1\}$ and $M\times\{0\}$, respectively. If the claim is true, then there exists a locally unipotent bundle $(U,\nabla_{u})$ on $M\times I$ that restricts to the flat vector bundles $(mV,m\nabla_{v})$ and $(\epsilon^{\star},d)$ on $\partial (M\times I)$. In particular, this implies $\int_{t\in I} \operatorname{ch}(\nabla_{u})=0$ and hence $[\tilde{\operatorname{ch}}(m\nabla_{v},\phi^{\ast}d)]=0$, where $\phi$ is an isomorphism from $mV$ to $\epsilon^{\star}$ induced from the homotopy $f^{\prime}$.
 
To show the claim, we observe the following diagram:
\begin{center}
\begin{tikzpicture}
\node(Lu) at (0,2) {$\partial W$};
\node(Ll) at (0,0) {$W$};
\node(Mu) at (2,2) {$L$};
\node(Ml) at (2,0) {$L^{\prime}$};
\node(Ru) at (4,2) {$\operatorname{BGL}(\mathbb{C}^{\delta})$};
\node(Rl) at (4,0) {$K_{a}\mathbb{C}$};

\path[->, font=\scriptsize,>=angle 90]

(Ll) edge node [above]{$f^{\prime}$}(Ml)
(Lu) edge node [above]{$f$}(Mu)   
(Ru) edge node [right]{$+$}(Rl);
\draw[right hook->](Lu) to (Ll);
\draw[right hook->](Mu) to (Ml);
\draw[right hook->](Mu) to (Ru);
\draw[right hook->](Ml) to (Rl);
\end{tikzpicture}
\end{center}
where $L$ and $L^{\prime}$ are finite $CW$-subcomplexes of $\operatorname{BGL}(\mathbb{C}^{\delta})$ and $K_{a}\mathbb{C}$ that contain the images of $f$ and $f^{\prime}$, respectively. Note also that $L^{\prime}$ can be obtained by attaching some $2$- or $3$-cells to $L$. By thickening $L$ and $L^{\prime}$, we get two smooth manifolds $N(L)$ and $N(L^{\prime})$ with $N(L^{\prime})$ given by attaching some $2$- or $3$-handles to $N(L)$ (see the diagram below). 
\begin{center}
\begin{tikzpicture}
\node(Lu) at (0,2) {$\partial W$};
\node(Ll) at (0,0) {$W$};
\node(Mu) at (2,2) {$N(L)$};
\node(Ml) at (2,0) {$N(L^{\prime})$};
\node(Ru) at (4,2) {$\operatorname{BGL}(\mathbb{C}^{\delta})$};
\node(Rl) at (4,0) {$K_{a}\mathbb{C}$};

\path[->, font=\scriptsize,>=angle 90]

(Ll) edge node [above]{$f^{\prime}$}(Ml)
(Lu) edge node [above]{$f$}(Mu)   
(Ru) edge node [right]{$+$}(Rl);
\draw[right hook->](Lu) to (Ll);
\draw[right hook->](Mu) to (Ml);
\draw[->](Mu) to (Ru);
\draw[->](Ml) to (Rl);
\end{tikzpicture}
\end{center}
Now there is a canonical flat vector bundle on $N(L)$ and, by the construction used in \cite[Section $3$]{Kr} (see also \cite[Section $2$]{Sch}), $N(L^{\prime})$ can be endowed with a locally unipotent bundle $(\bar{U},\nabla_{\bar{u}})$ that restricts to the canonical flat vector bundle on $N(L)$. Hence, $(f^{\prime,\ast}\bar{U},f^{\prime,\ast}\nabla_{\bar{u}})$ gives the required locally unipotent bundle on $W$. Thus, we have proved the claim.


 
\end{proof}



 
\section{$e$-invariants of Seifert homology spheres}

 
Following Jones and Westbury's idea, we derive a formula for the real part of $e$-invariants of Seifert homology spheres. We recall the formula for $\tilde{\xi}$-invariants of lens spaces obtained from Donnelly's fixed point theorem \cite[Section $4.5$]{Gi2} in the first subsection. In the second subsection, we review some properties of representations of the fundamental group of a Seifert homology sphere. A family of $4$-dimensional cobordisms between the connected sum of the lens spaces $\ConSum{\mathclap{1\leq i\leq n}} L(a_{i},b_{i})$ and the Seifert homology sphere $\Sigma(a_{1},...a_{n})$ is constructed via relative Kirby diagram in the third subsection. In the last subsection, using the results of the previous sections, we compute $e$-invariants of Seifert homology spheres, and particularly, we recover the formula for the real part of $e$-invariants of Seifert homology spheres in \cite[Theorem $C$]{JW} up to sign (see also \cite[Lemma $5.1$, Lemma $5.2$]{JW} and Remarks \ref{XiofLensinJW}, \ref{CobordisminJW} and \ref{ThmCDEinJW}).

\subsection{$\tilde{\xi}$-invariants of lens spaces}\label{Calculationforehomo}
Given a spin manifold $M$ and a unitary representation $\rho:\pi_{1}(M)\rightarrow U(k)$. The associated $\tilde{\xi}$-invariant is defined to be the difference 
\[\tilde{\xi}(\rho,M):=\frac{\eta+h}{2}(D_{\rho})-k\frac{\eta+h}{2}(D)\in\mathbb{R/Z},\]
where $D$ is the Dirac operator associated to the spin structure, $D_{\rho}$ is the associated twisted Dirac operator by the representation $\rho$, $\eta$ is the $\eta$-invariant \cite{APS1} and $h(A)$ is the dimension of the kernel of the operator $A$.
 
Donnelly's fixed point theorem for manifolds with boundary (\cite{Don} and \cite[Theorem $4.5.8$]{Gi2}) gives us the following formula for $\tilde{\xi}$-invariants of odd dimensional spherical space forms (see also \cite[Lemma $2.1$]{Gi}).
 
\begin{theorem}\label{XinvforQuotientspaces}
Let $G$ be a finite group, $\tau:G\rightarrow U(l)$ a fixed point free representation and $\tilde{\tau}$ a lifting of $\tau$ to $\operatorname{Spin}(2l)$. Denote the orbit space $S^{2l-1}/\tau(G)$ by $M$. Then we have the following formula for the $\tilde{\xi}$-invariant of the Dirac operator associated to the spin structure on $M$ and twisted by a representation $\rho:\pi_{1}(M)=G\rightarrow U(k)$:
\[\tilde{\xi}(\rho,M)=\frac{1}{\vert G\vert}\sum_{g\in G;g\neq id}(\operatorname{Tr(\rho(g))}-k)\operatorname{def}(\tau(g),\operatorname{spin}),\]
where 
\[\operatorname{def}(\tau(g),\operatorname{spin}):=\prod_{i=1}^{l}\frac{\sqrt{\lambda_{i}}}{\lambda_{i}-1}\] 
and $\lambda_{i}$, $i=1,...,l$, are eigenvalues of $\tau(g)$ in $U(l)$.
\end{theorem}
\begin{proof}
This is a special case of Theorem $4.5.9$ in \cite{Gi2}. 
\end{proof}
Let $G=\mathbb{Z}/p=<\lambda>$, where $\lambda\in\mathbb{C}$ and $\lambda^{p}=1$. Define the representation $\tau$ to be
$\tau(\lambda):=\operatorname{diag}(\lambda^{q_{1}},...,\lambda^{q_{l}})$ and denote $(q_{1},...,q_{l})$ by $\overrightarrow{q}$. Then the orbit space $M=S^{2l-1}/\tau(G)$ is a generalized lens space, denoted by $L(p;\overrightarrow{q})$. When $l$ is even and $p$ is odd, $L(p;\overrightarrow{q})$ has a unique spin structure, and when both $l$ and $p$ are even, it has two different spin structures, which are determined by $\overrightarrow{q}$. The following is an easy application of Theorem \ref{XinvforQuotientspaces}.
 
\begin{corollary}\label{XiforgenLensspaces}
Suppose the representation $\rho_{s}$ is given by $\rho_{s}(\lambda):=\lambda^{s}$. Then 
\[\tilde{\xi}(\rho,L(p,\overrightarrow{q}))=\frac{1}{p}\sum_{\mathclap{\substack{\lambda\in\mathbb{Z}/p\\ 
\lambda\neq id}}}(\lambda^{s}-1)\prod_{i=1}^{l}\frac{\sqrt{\lambda^{q_{i}}}}{\lambda^{q_{i}}-1}.\]  
\end{corollary}

In \cite[Lemma $2.3$-$4$]{Gi}, using the residue theorem, Gilkey has shown, if $q_{i}$ is coprime to the product $p\mu(\operatorname{Td}_{l}(\ast,\ast))$ for every $i$, where $\mu(\operatorname{Td}_{l}(\ast,\ast))$ is the denominator of $\operatorname{Td}_{l}(\ast,\ast)$, then
\[\frac{1}{p}\sum_{\mathclap{\substack{\lambda\in\mathbb{Z}/p\\
\lambda\neq id}}}\lambda^{s}\prod_{i=1}^{l}\frac{\sqrt{\lambda^{q_{i}}}}{\lambda^{q_{i}}-1}\equiv -\frac{d}{p}\operatorname{Td}(s-\frac{1}{2}\sum_{i=1}^{l} q_{i};p,\overrightarrow{q})\text{ mod }\mathbb{Z},\]
where $d$ is chosen so that the product $dq_{1}...q_{l}\equiv 1$ mod $p\mu(\operatorname{Td}_{l}(\ast,\ast))$;
\begin{equation}\label{eq:Tdl}
\operatorname{Td}_{l}(s;\overrightarrow{q}):=\sum_{a+b=l}s^{a}\frac{\operatorname{Td}_{b}(\overrightarrow{q})}{a!};
\end{equation} 
and $\operatorname{Td}_{b}(\overrightarrow{q})$ is the polynomial in $q_{1},...,q_{l}$ generated by the Todd polynomial $x/(1-e^{-x})$.

Gilkey's result can be easily generalized to the relative case. Namely, if $q_{i}$ is coprime to the product $p\mu(\operatorname{Td}_{l}(\ast,\ast)-\operatorname{Td}_{l}(0,\ast))$ for every $i$, then
\begin{multline}\label{eq:FromdeftoTd}
\tilde{\xi}(\rho_{s},L(p;\overrightarrow{q}))
=\frac{1}{p}\sum_{\mathclap{\substack{\lambda\in\mathbb{Z}/p\\ \lambda\neq id}}}(\lambda^{s}-1)\prod_{i=1}^{l}\frac{\sqrt{\lambda^{q_{i}}}}{\lambda^{q_{i}}-1}\\
\equiv -\frac{d}{p}[\operatorname{Td}(s-\frac{1}{2}\sum_{i=1}^{l} q_{i};p,\overrightarrow{q})-\operatorname{Td}(0-\frac{1}{2}\sum_{i=1}^{l} q_{i};p,\overrightarrow{q})]\text{ mod }\mathbb{Z},
\end{multline}
where $d$ is chosen so that the product $dq_{1}...q_{l}\equiv 1$ mod $p\mu(\operatorname{Td}_{l}(\ast,\ast)-\operatorname{Td}_{l}(0,\ast))$.
 
Formula \eqref{eq:FromdeftoTd} and Corollary \ref{XiforgenLensspaces} imply the following:

\begin{lemma}\phantomsection\label{FL}
Let $l=2$ and denote $L(p;1,q)$ by $L(p,q)$. Assume also $q$ is coprime to $2p$. Then we have the following formula:
\[\tilde{\xi}(\rho_{s},L(p,q))\equiv-\frac{ds^{2}}{2p}-\frac{ds}{2}\text{ mod }\mathbb{Z},\]
where $d$ is an integer satisfying $dq\equiv 1$ mod $2p$.
\end{lemma}
Note that, given $L(p,q)$, we can always assume $q$ is coprime to $2p$, and when $p$ is even, $q$ and $q+p$ specify the two different spin structures on $L(p,q)$.
\begin{proof}
From formula \eqref{eq:Tdl}, we know the denominator of $\operatorname{Td}_{2}(\ast,\ast)-\operatorname{Td}_{2}(0,\ast)$ is $2$, and by equality \eqref{eq:FromdeftoTd}, we can deduce
\begin{multline*}
\tilde{\xi}(\rho_{s},L(p,q))\\
\equiv -\frac{d}{2p}[(s-\frac{1}{2}(1+q))^{2}+(s-\frac{1}{2}(1+q))(1+q+p)-(-\frac{1}{2}(1+q))^{2}-(-\frac{1}{2}(1+q))(1+q+p)]\\
\equiv -\frac{d}{2p}(s^{2}+sp)
\text{ mod }\mathbb{Z},
\end{multline*} 
which proves the lemma.
\end{proof}

\begin{remark}\label{XiofLensinJW} 
It is stated in \cite[Lemma 5.2]{JW} that 
\[\tilde{\xi}(\rho_{s},L(p,q))\equiv\frac{-ds^{2}}{2p} \text{ mod } \mathbb{Z},\]
where $dq\equiv 1$ mod $p$. 
However, the right hand side of the equality varies when $d$ changes but $\tilde{\xi}(\rho_{s},L(p,q))$ should not depend on $d$. 
\end{remark}

\subsection{Seifert homology $3$-spheres}

In this subsection, we review some basic properties of Seifert homology $3$-spheres. A Seifert homology sphere $\Sigma(a_{1},...,a_{n})$ is a homology sphere and at the same time a Seifert manifold with $a_{i}$ the order of the $i$-th exceptional fiber---a Seifert manifold is a manifold admitting a Seifert fibering (see \cite[p.13]{Hat4}). It can be constructed topologically as follows (see \cite[p.2-4 and Figure $1.1$]{Sa}):
Let 
\[F=S^{2}\setminus int(D^{2}_{1}\cup...\cup D^{2}_{n})\] 
be a $n$-punctured $2$-sphere and $W\rightarrow F$ a $S^{1}$ fiber bundle with Euler number $b$. Assume the bundle over $\partial F$ has a fixed trivialization such that the boundary of $W$ can be identified with the $n$-tori $\partial D^{2}_{k}\times S^{1}$, $k=1...n$. Now we past another $n$ solid tori $D^{2}_{k}\times S^{1}$, $k=1,...,n$ in such a way that $a_{k}(S^{1}\times\{1\}) + b_{k}(\{1\}\times S^{1})$ in the $k$-th boundary component of $W$ is null homotopic in $D^{2}_{k}\times S^{1}$ after pasting, where $\{(a_{k},b_{k})\}_{k=1}^{n}$ $n$ pairs of relatively prime numbers. The following figure is an expression of $\Sigma(a_{1},...,a_{n})$ in terms of Dehn surgery: 
 
\begin{figure}[h]
\caption{The linking diagram of $\Sigma(a_{1},a_{2},...,a_{n})$}
\label{LD1}
\centering
\includegraphics[width=0.6\textwidth]{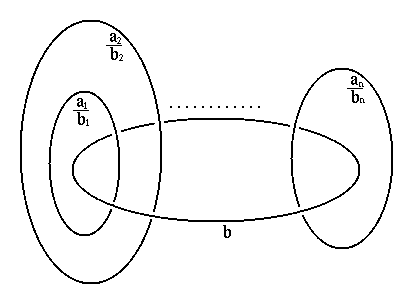}
\end{figure}

\noindent 
From the above figure we see the fundamental group of $\Sigma(a_{1},...,a_{n})$ can be expressed as follows:
\begin{multline}\label{Representationofpi1Sigma}
\pi_{1}(\Sigma(a_{1},a_{2},...,a_{n}))\\
=< h,x_{1},...x_{n}\vert 
 [x_{i},h]=1,x_{1}...x_{n}=h^{-b},x_{i}^{a_{i}}=h^{-b_{i}}, \forall i>.
\end{multline} 
Since it is a homology sphere, $a_{i},b_{i}$, $i=1,...n$ must satisfy the equation:
\[a_{1}\cdots a_{n}(-b+\sum_{i}\frac{b_{i}}{a_{i}})=1.\]
Now by the theorem in \cite{Ki1} (see also \cite[p.264,p.278 Remark I.7]{Rol}), we see that the diffeomorphic type of a Seifert homology sphere is completely determined by these pairwise coprime numbers $a_{1},...,a_{n}$.

\noindent
\textbf{Irreducible representations} \cite[p.932]{JW}:
Observe first that any complex representation of $\pi_{1}(\Sigma(a_{1},...,a_{n}))$ factors through $\operatorname{SL}_{N}(\mathbb{C})$ as $\pi_{1}(\Sigma(a_{1},...,a_{n}))$ is always perfect. Secondly, since we are only interested in the representations 
\[\rho:\pi_{1}(\Sigma(a_{1},...,a_{n}))\rightarrow \operatorname{SL}_{N}(\mathbb{C})\] that has $\rho(h)$ is a scalar multiple of the identity matrix, for example, irreducible representations, we may assume
\begin{align*}
\rho(h)&=\lambda_{h}I;\\
\lambda_{h}&=\zeta^{r_{h}}_{N};\\
\zeta_{N}&:= e^{2\pi i/N}. 
\end{align*} 
Let $\lambda_{1}(j),...,\lambda_{N}(j)$ are the eigenvalues of $\rho(x_{j})$, then they cane be expressed as follows: 
\begin{equation}\label{Sjksymbol}
\lambda_{k}(j)=\zeta^{Ns_{k}(j)-b_{j}r_{h}}_{Na_{j}},
\end{equation} 
where $s_{k}(j)\in\mathbb{Z}$ and $0\leqslant s_{k}(j)<a_{j}$, for every  
$1\leqslant j\leqslant n$ and $1\leqslant k\leqslant N.$ 
In this case, we say $\rho$ is the representation of type $(s_{k}(j))$.
 

\subsection{$4$-dimensional cobordisms via relative Kirby diagrams}  
Figure \ref{LD1} induces a relative Kirby diagram (Figure \ref{LD0}) whose corresponding $4$-dimensional cobordism, denoted by $W$, bounded by the disjoint union of the connected sum of lens spaces $\ConSum{\mathclap{1\leq i\leq n}}L(a_{i},b_{i})$ and the Seifert homology $3$-sphere $\Sigma(a_{1},...,a_{n})$. More precisely, $W$ can be constructed by attaching a $2$-handle along the circle with the framing coefficient $b$ to the product of $\ConSum{\mathclap{1\leq i\leq n}}L(a_{i},-b_{i})\times I$ (see \cite[Section $5.5$]{GS}). Namely, we have  
\begin{equation}\label{ThecobordismbetweenLandS}
W:=\ConSum{\mathclap{1\leq i\leq n}}L(a_{i},-b_{i})\times I\cup_{(\partial D^{2}) \times D^{2}} D^{2}\times D^{2}.
\end{equation}  
In order to use $W$ to compute the $\tilde{\xi}$-invariants of $\Sigma(a_{1},...,a_{n})$, we yet need to check if the $4$-dimensional cobordism $W$ is spin.

We first recall that every closed oriented $3$-manifold can be obtained by performing Dehn (rational) surgery on a link $L$ in $S^{3}$ (see \cite[p.273]{Rol}). We denote the resulting $3$-manifold by $M_{L}$. It is also proved by Kirby that $M_{L}$ and $M_{L^{\prime}}$ are diffeomorphic if and only if $L$ can be obtained from $L^{\prime}$ by performing a sequence of Rolfsen moves or introducing or deleting an unknot with coefficient $\infty$ (see \cite[p.278]{Rol}). In particular, utilizing Rolfsen moves and deleting an unknot with coefficient $\infty$, we have the following useful move, called slam-dunk by Cochren: 
\begin{figure}[h]
\caption{Slam-dunk}
\label{LD2}
\centering
\includegraphics[width=0.5\textwidth]{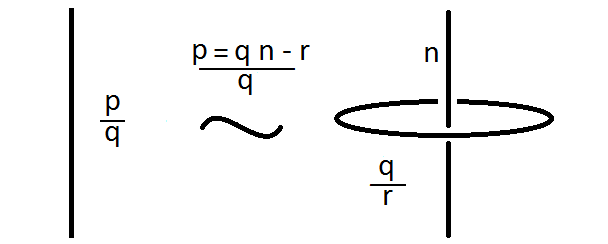}

\end{figure}

If $p,q$ both are larger than zero. We choose the unique integers $n$ and $r$ that satisfy $p=nq-r$ with $0\leq r <q$. In this way, we obtain an algorithm that turns a rational surgery on a knot into an integral surgery on a link when $p,q$ are larger than zero. By an integral surgery, we understand a Dehn surgery with surgery coefficients integers---it is also called framing coefficients. The following illustrates how via the slam-dunk algorithm a rational surgery diagram of the lens space $L(5,-2)$ is turned into an integral surgery diagram: 
 
\begin{figure}[h]
\caption{}
\label{LD3}
\centering
\includegraphics[width=0.6\textwidth]{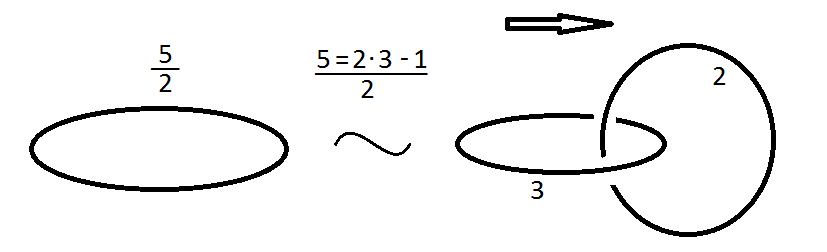}
\end{figure}

\noindent
Note the arrow indicates from where the algorithm starts. We also recall that, given an integral surgery diagram $L$, one can view it as a Kirby diagram, meaning, considering $M_{L}$ as the boundary of the $4$-manifold (or $2$-handlebody) $W_{L}$. Now Wu's formula says the second Stiefel-Whitney class $w_{2}(W_{L})$ is uniquely determined by the identity in $\mathbb{Z}/2$: 
\[<w_{2}(W_{L}),x>=x\cdot x\] 
for every $x\in H_{2}(W_{L};\mathbb{Z}/{2})$, where $x \cdot x$ means the self intersection number of $x$. This, when written in terms of linking numbers of the components in $L$, gives the following definition, which is the key to visualize spin structures on the $3$-manifold $M_{L}$ and obstructions to extending them over $W_{L}$.

\begin{definition}
Given a framed link (integral surgery diagram) $L$, a characteristic sublink $L^{\prime}\subset L$ is a sublink with $\operatorname{lk}(L^{\prime},K)\equiv \operatorname{lk}(K,K) (\operatorname{mod} 2)$ for every component $K$ in $L$, where $\operatorname{lk}(L_{1},L_{2})$ stands for the linking number of two links $L_{1}$ and $L_{2}$.
\end{definition}
We need the orientation to define linking numbers of links, but since we are working in $\mathbb{Z}/2$, choosing which orientations is not so important to us.

\begin{lemma}
Given a framed link $L$, there is a $1-1$ correspondence between the set of spin structures on $M_{L}$ and the set of characteristic sublinks of $L$.
\end{lemma}
\begin{proof}
The detailed proof can be found in \cite[Proposition 5.7.11]{GS}. We only sketch how the bijection is constructed. Given a spin structure $s$ on $M_{L}$, we define the relative obstruction class $w_{2}(W_{L},s)\in H^{2}(W_{L},M_{L};\mathbb{Z}/2)\simeq H_{2}(W_{L};\mathbb{Z}/2)$ to be the obstruction to expending the spin structure $s$ on $M_{L}$ over $W_{L}$. Since the classes in $H_{2}(W_{L};\mathbb{Z}/2)$ correspond bijectively to the sublinks of $L$ and the image of $w_{2}(W_{L},s)$ in $H^{2}(W_{L};\mathbb{Z}/2)$ must be the second Stiefel-Whitney class, we see the sublink corresponding to $w_{2}(W_{L},s)$ has to be a characteristic sublink, in view of Wu's formula. $s\mapsto w_{2}(W_{L},s)$ gives us the required bijection.
\end{proof}

\newpage
The following illustrates the unique spin structure on $L(5,2)$ as well as $L(7,3)$ and the two spin structures on $L(8,3)$---(red) bold circles denote the characteristic sublinks:
 
\begin{figure}[h]
\caption{}
\label{LD4}
\centering
\includegraphics[width=0.6\textwidth]{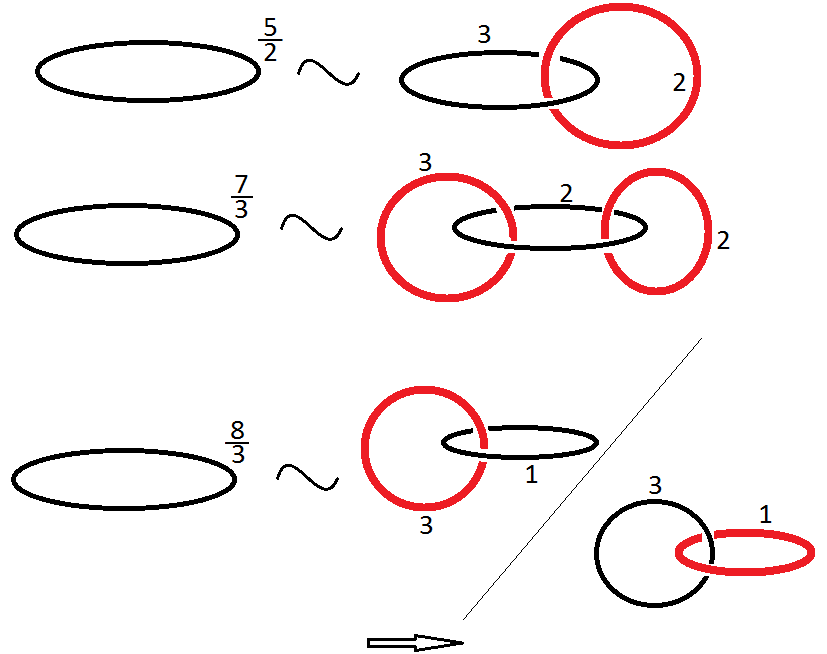}
\end{figure}

\begin{figure}[h]
\caption{}
\label{LD5}
\centering
\includegraphics[width=0.6\textwidth]{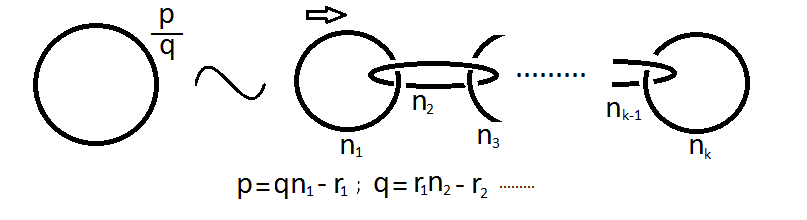}
\end{figure}

Given a lens space $L(p,-q)$, we call the Dehn surgery on an unknot with surgery coefficient $\frac{p}{q}$ the canonical rational surgery diagram of $L(p,-q)$ (the left-hand side of Figure \ref{LD5}) and call the integral surgery on $L_{\frac{p}{q}}$ the framed link obtained by applying the slam-dunk algorithm to the canonical surgery diagram of $L(p,-q)$ the canonical integral surgery diagram of $L(p,-q)$ (the right-hand side of Figure \ref{LD5}).
 
The following three lemmas describe some general properties of the canonical integral surgery diagram of $L(p,-q)$ when $p$ and $q$ both are larger than $0$.
\begin{lemma}\phantomsection\label{LemmaSpinonL1}
Given $L(p,-q)$ with $p$ even, then the two characteristic sublinks (spin structures) of its canonical integral surgery diagram can be distinguished by whether the first component is in the characteristic sublink.
\end{lemma}
\begin{proof}
Note that any closed oriented $3$-manifold $M$ is parallelizable and hence spin. Furthermore, spin structures on $M$ are parameterized by its first singular cohomology with coefficients in $\mathbb{Z}/2$. Since  
\[H^{1}(L(p,-q);\mathbb{Z}/2)=\mathbb{Z}/2,\] 
there are two spin structures on $L(p,-q)$. Suppose we are given a characteristic sublink $L^{\prime}$ and color it. Then we claim that, in order to obtain another characteristic sublink $L^{\prime\prime}$, the color of the first component $K_{1}$ must be changed. If it is not the case, the color of the second component should also stay unchanged because otherwise we have the contradiction:
\[\operatorname{lk}(L^{\prime\prime},K_{1})\equiv\operatorname{lk}(L^{\prime},K_{1})+1\not\equiv \operatorname{lk}(K_{1},K_{1})\hspace*{2em}(\operatorname{mod} 2).\] 
By induction, we assume, for every $i\leq k$, the color of the $i$-th component remains unaltered, and we want to show the color of the $(k+1)$-th component should stay the same. This can be seen easily since if the color of the $(k+1)$-th component is changed, the following gives a contradiction 
\[\operatorname{lk}(L^{\prime\prime},K_{k})\equiv\operatorname{lk}(L^{\prime},K_{k})+1\not\equiv \operatorname{lk}(K_{k},K_{k})\hspace*{2em}(\operatorname{mod}2).\]
Thus, the first component must have its color changed otherwise $L^{\prime\prime}=L^{\prime}$. 
\end{proof}
For the sake of convenience, we continue to assume components in a given characteristic sublink are colored. 
\begin{lemma}\phantomsection\label{LemmaSpinonL2}
Given $L(p,-q)$ with $p$ odd and $q$ even, then the first component of $L_{\frac{p}{q}}$ is not in the characteristic sublink.
\end{lemma}
\begin{proof}
Suppose $p=qn-r$ with $0<r<q$. Then the first component of the canonical surgery diagram has coefficient $n$. Since $q$ is even, we know from Lemma \ref{LemmaSpinonL1} the characteristic sublinks of $L_{\frac{q}{r}}$ are distinguished by the color of its first component. Now if $n$ is even, we choose the characteristic sublink of $L_{\frac{q}{r}}$ that has its first component uncolored and observe that the characteristic sublink of $L_{\frac{q}{r}}$ can be extended to the characteristic sublink of $L_{\frac{p}{q}}$, thinking of $L_{\frac{q}{r}}$ as a sublink of $L_{\frac{p}{q}}$. If $n$ is odd, we pick up the characteristic sublink of $L_{\frac{q}{r}}$ whose first component is colored and it is also easy to see this characteristic sublink of $L_{\frac{q}{r}}$ is extended to the characteristic sublink of $L_{\frac{p}{q}}$. In either case, we have the first component of the canonical integral surgery diagram of $L(p,-q)$ is not colored.
\end{proof}
Before showing the third lemma, we first recall the effect of handle sliding on characteristic sublinks (\cite[p.190]{GS}): Given $K$ and $K^{\prime}$ two components in a colored link, if we slide $K$ over $K^{\prime}$, then $K^{\prime}$ has its color changed if and only if $K$ is colored.
 
\begin{lemma}\phantomsection\label{LemmaSpinonL3}
Given $L(p,-q)$ with both $p$ and $q$ odd, then the first component of $L_{\frac{p}{q}}$ is in the characteristic sublink.
\end{lemma}
\begin{proof}
We first claim: Given the canonical integral surgery diagram of $L(p,-q)$ as follows:  
 
\begin{figure}[h]
\caption{}
\label{LD6}
\centering
\includegraphics[width=0.5\textwidth]{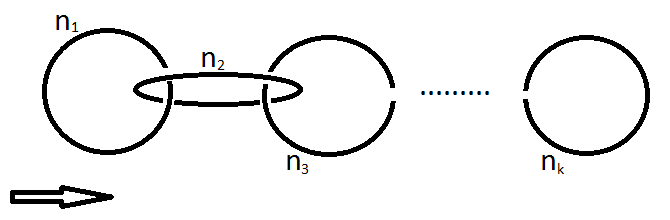}  
\end{figure}

\noindent
we have the expression of the canonical integral surgery diagram of $L(p,-q-p)$ looks like: 
 
\begin{figure}[h]
\caption{}
\label{LD7}
\centering
\includegraphics[width=0.6\textwidth]{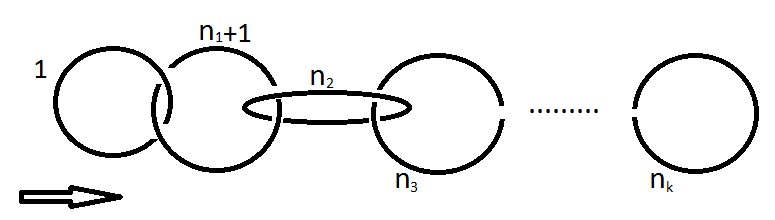}
\end{figure}

\noindent
This can be seen by the following computation: Suppose $p=n_{1}q-r$ with $0<r<q$. Meaning, after applying slam-dunk once to the canonical rational surgery diagram of $L(p,-q)$, the surgery coefficient of the second component is $\frac{q}{r}$. On the other hand, considering the canonical rational surgery diagram of $L(p,-q-p)$, we have the first slam-dunk gives us an unknot with the surgery coefficient $1$ because $p=(p+q)1-q$ and $0<q<p+q$. Applying slam-dunk again, we have $qm - s=p+q$ with $0<s<q$. Now since $p=n_{1}q-r$, we obtain $m=n_{1}+1$ and $s=r$. This means, after applying slam-dunk twice to the canonical surgery diagram of $L(p,-p-q)$, we get the same surgery coefficient $\frac{q}{r}$ as we did after applying slam-dunk once to the canonical rational surgery diagram of $L(p,-q)$. In this way, we see that deleting the first component of the link in Figure \ref{LD6} and deleting the first two components in Figure \ref{LD7} result in the same framed link. Furthermore, from this computation, we also know, if the first component of the link in Figure \ref{LD6} is $n_{1}$, the second component of the link in Figure \ref{LD7} must be $n_{1}+1$, and the framing of the first component of the link in Figure \ref{LD7} is always $1$. This proves the claim.
  
Secondly, we note Figure \ref{LD7} can be obtained by sliding the first component of Figure \ref{LD6} over a separate unknot with framing $1$ as illustrated below:

\begin{figure}[h]
\caption{}
\label{LD8}
\centering
\includegraphics[width=0.6\textwidth]{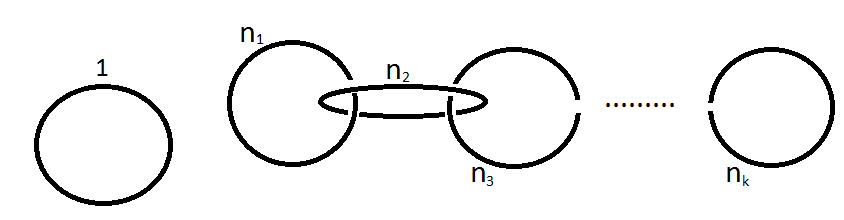}
\end{figure}

\noindent
and the characteristic sublink of Figure \ref{LD8} is the union of the characteristic sublink of $L_{\frac{p}{q}}$ and that separated unknot. Now, if the first component of $L_{\frac{p}{q}}$ is not in the characteristic sublink. Sliding the first component of $L_{\frac{p}{q}}$ over this separated unknot does not change the color of this unknot, so we get the characteristic sublink of $L_{\frac{p}{p+q}}$ contains the first component of $L_{\frac{p}{p+q}}$. Yet it is not possible because $p+q$ is even, and Lemma \ref{LemmaSpinonL2} tells us the first component of $L_{\frac{p}{p+q}}$ must not in the characteristic sublink. We get a contradiction, and therefore the first component of $L_{\frac{p}{q}}$ has to be in the characteristic sublink. The proof is completed. 
\end{proof}

Before returning to the $4$-dimensional cobordism $W$, we recall that, given an integral relative Kirby diagram with a bracketed framed link $L$ and an unbracketed framed knot $K$ in $S^{3}$, the associated $4$-dimensional cobordism is given by $W_{L\cup K}:=M_{L}\times I\cup_{\mathfrak{N}(K)} D^{2}\times D^{2}$, where the attaching diffeomorphism from $S^{1}\times D^{2}$ to $\mathfrak{N}(K)$ is determined by the framing of $K$ (see \cite[p.175-176]{GS} for more on relative Kirby calculus). Now suppose a spin structure, a characteristic sublink $L^{c}$ of $L$, of $M_{L}$ is given. Then this spin structure of $M_{L}$ can be extended over $W_{L\cup K}$ if and only if 
\[\operatorname{lk}(L^{c},K)\equiv \operatorname{lk}(K,K) \hspace*{2em}(\operatorname{mod} 2)\]
(see \cite[p.189-190]{GS}).

Applying slam-dunk to relative Kirby diagram \eqref{LD0}, we obtain the following (integral) relative Kirby diagram:
 
\begin{figure}[h]
\caption{}
\label{LD9}
\centering
\includegraphics[width=0.6\textwidth]{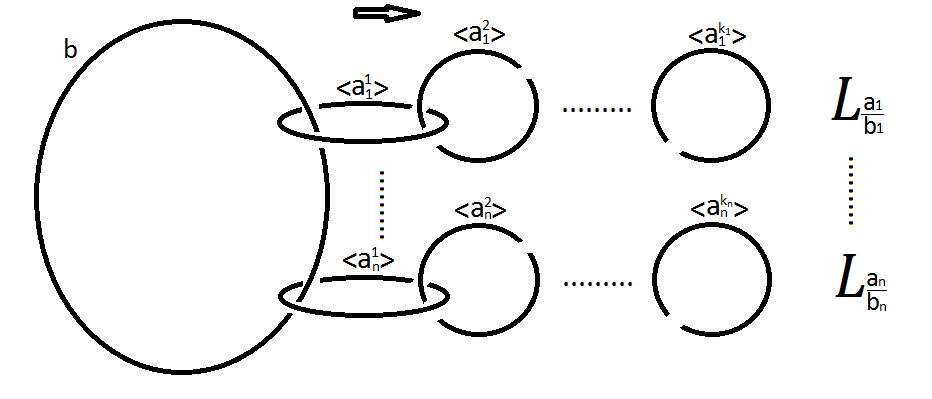}
\end{figure}
\noindent
which also represents the $4$-dimensional cobordism $W$, where $a_{i}^{j}\in\mathbb{Z}$. With this diagram in mind, we want to determine under what conditions $W$ is spin. Firstly, we assume $a_{i}$ is odd for every $i\in\{1,...,n\}$. Hence, without loss of generality, we may assume $b_{i}$ is even for every $i$. In view of Lemma \ref{LemmaSpinonL2}, we know the first component of $L_{\frac{a_{i}}{b_{i}}}$ does not belong to the characteristic sublink, for every $i$. Let $L^{c}$ be the union of the characteristic sublinks of $L_{\frac{a_{i}}{b_{i}}}$, for all $i$, and $K$ be the knot with the framing coefficient $b$ (see Figure \ref{LD1}). Then we get $\operatorname{lk}(L^{c},K)\equiv 0 (\operatorname{mod} 2)$ (see Figure \ref{LD9}). On the other hand, the equality
\[a_{1}...a_{n}(-b+\sum_{i=1}^{n}\frac{b_{i}}{a_{i}})=1\]   
gives $b\equiv 1 (\operatorname{mod} 2)$. Thus, in this case, the spin structure cannot be extended over $W$.

Suppose one of $a_{1},...,a_{n}$ is even, $a_{1}$ say. Then there are two characteristic sublinks for the canonical integral surgery diagram of $\underset{i}{\#}L(a_{i},b_{i})$. We may assume $b_{i}$ is even when $i\neq 1$ and let $L_{1}^{c}$ and $L_{2}^{c}$ be the characteristic sublinks corresponding to the two spin structures on $\underset{i}{\#}L(a_{i},b_{i})$. By Lemma \ref{LemmaSpinonL1}, we may assume $\operatorname{lk}(L_{1}^{c},K)\equiv 1 (\operatorname{mod} 2)$ and $\operatorname{lk}(L_{2}^{c},K)\equiv 0 (\operatorname{mod} 2)$ (see Figure \ref{LD9} and note the only possible component of these characteristic sublinks that has non-trivial linking number with $K$ is the first component of $L_{\frac{a_{1}}{b_{1}}}$). So when $b$ is odd, the spin structure on $\underset{i}{\#}L(a_{i},b_{i})$ corresponding to $L_{1}^{c}$ can be extended over $W$, whereas the spin structure corresponding to $L_{2}^{c}$ is extendable over $W$ when $b$ is even.

We summarize some properties of $W$ in the following proposition:
 
\begin{proposition}\phantomsection\label{SpinstronW}
\begin{enumerate}
\item $\partial W=\ConSum{\mathclap{1\leq i\leq n}}L(a_{i},b_{i})\coprod \Sigma(a_{1},...,a_{n})$ and 
\[\pi_{1}(W)=<x_{1},x_{2},...,x_{n}\mid x_{1}^{a_{1}}=...=x_{n}^{a_{n}}=x_{1}...x_{n}=1>.\] 
\item There is a short exact sequence 
\[0\rightarrow \mathbb{Z}\rightarrow \pi_{1}(\Sigma(a_{1},...,a_{n}))\rightarrow \pi_{1}(W)\rightarrow 0,\]
where the first homomorphism sends the generator of $\mathbb{Z}$ to $h$, and the second homomorphism is induced by the inclusion $\Sigma(a_{1},...,a_{n})\hookrightarrow W$.
\item When the natural numbers $a_{1},...,a_{n}$ are all odd, $W$ is not spinnable.

\item When one of the natural numbers $a_{1},...,a_{n}$ is even, $W$ admits a unique spin structure. 

\item The $\operatorname{spin}^{c}$ structures on $W$ are parameterized by $\mathbb{Z}$. 

\item The $\operatorname{spin}^{c}$ structures on $\underset{i}{\#}L(a_{i},b_{i})$ are parameterized by $2\mathbb{Z}/a_{1}...a_{n}$ when the natural numbers $a_{1},...,a_{n}$ are all odd numbers, whereas when one of the natural numbers $a_{1},...,a_{n}$ is even, they are parameterized by $\mathbb{Z}/2\oplus 2\mathbb{Z}/a_{1}...a_{n}$.

\item When the natural numbers $a_{1},...,a_{n}$ are all odd, there is a $\operatorname{spin}^{c}$ structure on $W$ which restricts to the canonical $\operatorname{spin}^{c}$ structures on $\Sigma(a_{1},...,a_{n})$ and $\underset{i}{\#}L(a_{i},-b_{i})$. By the canonical $\operatorname{spin}^{c}$ structure on a spin manifold, we mean the $\operatorname{spin}^{c}$ structure induced by the spin structure and the trivial complex line bundle.
\end{enumerate}
\end{proposition}
\begin{proof}
Recall first the singular cohomology groups of the lens space $L(p,q)$:
\[H^{\ast}(L(p,q);\mathbb{Z})=\begin{cases}
\mathbb{Z}/p& \ast=2\\
\mathbb{Z}& \ast=0,3  \\
0& \text{otherwise}
\end{cases}.\]
From this, the singular cohomology groups of the connected sum $\underset{i}{\#}L(a_{i},b_{i})$ can be easily deduced---remember the natural numbers $a_{1},...,a_{n}$ are pairwise coprime integers.  
\[H^{\ast}(\underset{i}{\#}L(a_{i},b_{i});\mathbb{Z})=\begin{cases}
\mathbb{Z}/a_{1}...a_{n}& \ast=2\\
\mathbb{Z}& \ast=0,3 \\
0& \text{otherwise}
\end{cases}.\] 
With the universal coefficient theorem, we can further obtain the singular cohomology groups of $\underset{i}{\#}L(a_{i},b_{i})$ with coefficients in $\mathbb{Z}/2$:
\[H^{\ast}(\underset{i}{\#}L(a_{i},b_{i});\mathbb{Z}/2)=\begin{cases}
\mathbb{Z}/2& \ast=0,3  \\
0& \text{otherwise}
\end{cases}\] 
when the product $a_{1}....a_{n}$ is odd;
\[H^{\ast}(\underset{i}{\#}L(a_{i},b_{i});\mathbb{Z}/2)=\begin{cases}
\mathbb{Z}/2& \ast=0,3 \\
\mathbb{Z}/2& \ast=1,2 \\
0 & \text{otherwise}
\end{cases}\]
when the product $a_{1}...a_{n}$ is even.

Using the long exact sequence induced by the pair $(W,\Sigma(a_{1},...,a_{n}))$ or the pair $(W,\underset{i}{\#}L(a_{i},b_{i}))$, we can compute the singular cohomology groups of $W$ with coefficients in both $\mathbb{Z}$ and $\mathbb{Z}/2$:
\[H^{\ast}(W;\mathbb{Z})=\begin{cases}
\mathbb{Z}& \ast=2,3 \\
0& \text{otherwise}
\end{cases};\]
\[H^{\ast}(W;\mathbb{Z}/2)=\begin{cases}
\mathbb{Z}/2& \ast=2,3 \\
0& \text{otherwise}
\end{cases}.\]
In particular, we see $W$ admits a $\operatorname{spin}^{c}$ because the second Stiefel-Whitney class is the mod $2$ reduction of an integral class, in view of the universal coefficient theorem. In fact, every $4$-manifold admits a $\operatorname{spin}^{c}$ structure (see \cite[Remark 5.7.5]{GS}).

Using the computation above we can identify the parameter sets for the sets of $\operatorname{spin}^{c}$ and spin structures on $X=W$ or $\underset{i}{\#}L(a_{i},-b_{i})$. More precisely, the $\operatorname{spin}^{c}$ structures are parameterized by $2H^{2}(X;\mathbb{Z})\oplus H^{1}(X;\mathbb{Z}/2)$ (see \cite[p.392]{LM}), whereas $H^{1}(X;\mathbb{Z}/2)$ parameterizes the set of spin structures when $X$ is spin (see \cite[p.82]{LM}).

To see the last claim, we note that to find a $\operatorname{spin}^{c}$ structure described in Proposition \ref{SpinstronW} on $W$ is equivalent to find a class $x$ in $H^{2}(W;\mathbb{Z})$ such that the following two conditions are fulfilled:\par 
\textbf{i}: Let $H^{2}(W;\mathbb{Z})\rightarrow H^{2}(W;\mathbb{Z}/2)$ be the homomorphism induced by the epimorphism $\mathbb{Z}\rightarrow \mathbb{Z}/2$, then the image of $x$ under this homomorphism in $H^{2}(W;\mathbb{Z}/2)$ is the second Stiefel-Whitney class.\par
\textbf{ii}: Let $H^{2}(W;\mathbb{Z})\rightarrow H^{2}(\underset{i}{\#}L(a_{i},-b_{i}),\mathbb{Z})$ be the homomorphism induced by the inclusion $\underset{i}{\#}L(a_{i},-b_{i})\hookrightarrow W$. Then the image of $x$ under this homomorphism in $H^{2}(\underset{i}{\#}L(a_{i},-b_{i}),\mathbb{Z})$ is trivial.\par
\noindent
It is not difficult to see the first condition says the associated complex line bundle determines a $\operatorname{spin}^{c}$ structure on $W$ and the second condition implies this complex line bundle restricts to the trivial one on the boundary of $W$.

To see how to find such $x$, we first observe the homomorphism
\[ H^{2}(W;\mathbb{Z})\rightarrow  H^{2}(W;\mathbb{Z}/2) \] 
can be identified with the surjective homomorphism $\mathbb{Z}\rightarrow \mathbb{Z}/2$ as $H^{3}(W;\mathbb{Z})=\mathbb{Z}$. Thus, any odd element in $H^{2}(W;\mathbb{Z})$ satisfies the first condition. Secondly, because $H^{3}(W,\underset{i}{\#}L(a_{i},-b_{i});\mathbb{Z})=0$, the homomorphism 
\[ \mathbb{Z}\simeq H^{2}(W;\mathbb{Z})\rightarrow  H^{2}(\underset{i}{\#}{L(a_{i},-b_{i})};\mathbb{Z})\simeq \mathbb{Z}/a_{1}...a_{n}\]
is surjective. Hence, any element that is divided by $a_{1}...a_{n}$ satisfies the second condition. Now the product $a_{1}...a_{n}$ is odd, so we can choose $x$ to be the class corresponding to the element $ka_{1}...a_{n}\in\mathbb{Z}$ with $k$ odd.  

\end{proof}


 
\begin{remark}\phantomsection\label{CobordisminJW} 
Another construction of a $4$-dimensional cobordism whose boundary is homeomorphic to $\ConSum{\mathclap{1\leq i\leq n}}L(a_{i},b_{i})\coprod\Sigma(a_{1},...,a_{n})$ is given in \cite[p.951]{JW}. We denote it by $W_{\operatorname{JW}}$ as it is different from the $4$-dimensional cobordism considered here. It is stated that $W_{\operatorname{JW}}$ is always spin \cite[Lemma $5.1$]{JW}, but the proof given there appears to contain some gaps. In fact, the homomorphism $H^{2}(W_{\operatorname{JW}} ,\mathbb{Z}/2)\rightarrow H^{2}(W_{\operatorname{JW}} \setminus A,\mathbb{Z}/2)$ is not an isomorphism as claimed there \cite[p.951]{JW}. Hence, the fact that $W_{\operatorname{JW}} \setminus A$ is parallelizable does not imply $W_{\operatorname{JW}}$ admits a spin structure. Actually, from the construction given there, the cobordism $W_{\operatorname{JW}}$ should have $H^{2}(W_{\operatorname{JW}},\mathbb{Z})=\mathbb{Z}\oplus \operatorname{Torsions}$, while $H^{2}(W \setminus A,\mathbb{Z})$ is trivial. This can be seen by the Mayer-Vietoris sequence of $(W_{\operatorname{JW}},N(A),W_{\operatorname{JW}}\setminus A)$ and the fact that $N(A)$ the neighborhood of $A$, $N(A)\cap (W_{\operatorname{JW}} \setminus A)$, and $W_{\operatorname{JW}}\setminus A$ are homotopy equivalent to $\coprod\limits_{n} S^{1}$, $T^{2}$, and a $n$-punctured $2$-sphere, respectively. By counting the ranks of their second integral cohomology groups, we see there is at least one free element in $H^{2}(W_{\operatorname{JW}},\mathbb{Z})$. The torsion part of $H^{2}(W_{\operatorname{JW}};\mathbb{Z})$ can be obtained from computing $H^{1}(W_{\operatorname{JW}};\mathbb{Z})$, which is isomorphic to $\mathbb{Z}/a_{1}...a_{n}$. This also shows $W$ and $W_{\operatorname{JW}}$ are not homotopic equivalent. 


\end{remark}

\subsection{$e$-invariants of Seifert homology spheres} 
Following Jones and Westbury's idea, we now compute $e$-invariants of Seifert homology spheres. We restrict our attention to those representations $\rho$ of 
\[\pi_{1}(\Sigma(a_{1},a_{2},...a_{n}))=< h,x_{1},...x_{n}\vert 
 [x_{i},h]=1,x_{1}...x_{n}=h^{-b},x_{i}^{a_{i}}=h^{-b_{i}}, \forall i>\]
in $\operatorname{SL}_{N}(\mathbb{C})$ that have $\rho(h)$ is a scalar multiple of the identity. 

Recall first that Lemma $5.4$ in \cite{JW} says, given a representation $\rho$ of this kind, we have 
\[e_{\ast}[\Sigma,\rho\otimes\bar{\rho}]=2N\mathfrak{Re}(e_{\ast}[\Sigma,\rho]).\]

Observe also, since the tensor product $\rho\otimes \bar{\rho}$ sends $h$ to $id$, $\rho\otimes \bar{\rho}$ can be extended to a representation $\varrho$ of $\pi_{1}(W)$ (see Proposition \ref{SpinstronW}) such that  
\begin{align*}
\varrho\circ \iota_{\ast} &=\rho\otimes\bar{\rho} ;
\\
\varrho\circ \iota^{\prime}_{\ast} &=\bigoplus\limits_{1\leq k,l\leq N}\rho_{s_{k}-s_{l}},
\end{align*}
where $\iota$ and $\iota^{\prime}$ are the inclusions from $\Sigma(a_{1},...,a_{n})$ and $\ConSum{\mathclap{1\leq i\leq n}}L(a_{i},-b_{i})$ into $W$, respectively, and $\rho_{s_{k}-s_{l}}$ is the representation induced by the representations $\rho_{s_{k}(i)-s_{l}(i)}$ on $\pi_{1}(L(a_{i},b_{i}))$, for $1\leq i\leq n$.  Use the fact that the $\tilde{\xi}$-invariant is invariant under spin cobordism (see \cite[Theorem 3.3]{APS2}), we can derive the following formula:


\begin{proposition}\phantomsection\label{Ree} 
Let $\rho$ be a representation of $\pi_{1}(\Sigma(a_{1},...,a_{n}))$ in $\operatorname{SL}_{N}(\mathbb{C})$ of the type $s=(s_{k}(i))_{1\leq i\leq n; 1\leq  k\leq  N}$
with $\rho(h)$ a scalar multiple of the identity matrix. Then 
\[2N\mathfrak{Re}(e_{\ast}[\Sigma(a_{1},a_{2},...a_{n}),\rho])= \sum_{i=1}^{n}\sum_{k=1}^{N}\sum_{l=1}^{N}
\frac{a(s_{k}(i)-s_{l}(i))^{2}}{a_{i}},\]
where $a:=a_{1}...a_{n}$. 
\end{proposition}

\begin{proof}
By Proposition \ref{SpinstronW}, we know the cobordism $W$ constructed there is spin when the product $a=a_{1}...a_{n}$ is even, and when the product $a=a_{1}...a_{n}$ is odd, it admits a $\operatorname{spin}^{c}$ structure that restricts to the canonical $\operatorname{spin}^{c}$ structures on the boundary. Using the fact the $\tilde{\xi}$-invariant is invariant under spin ($\operatorname{spin}^{c}$) cobordism, we obtain the following identities: 
\begin{multline*}
\sum_{1\leq k,l\leq N}\sum_{i=1}^{n}\tilde{\xi}(\rho_{s_{k}(i)-s_{l}(i)},L(a_{i},-b_{i}))=\sum_{1\leq k,l\leq N}\tilde{\xi}(\rho_{s_{k}-s_{l}},\hspace{.5em}\ConSum{\mathclap{1\leq i\leq n}}L(a_{i},-b_{i}))\\
=\tilde{\xi}(\rho\otimes\bar{\rho},\Sigma(a_{1},...,a_{n})).
\end{multline*} 
Plugging in formula \eqref{FL} for $\tilde{\xi}$-invariants of lens spaces, we get 
\begin{multline*}
\tilde{\xi}(\rho\otimes\bar{\rho},\Sigma(a_{1},...,a_{n}))=\sum_{i=1}^{n}\sum_{k=1}^{N}\sum_{l=1}^{N}
(\frac{b^{\prime}_{i}(s_{k}(i)-s_{l}(i))^{2}}{2a_{i}}+\frac{b_{i}^{\prime}(s_{k}(i)-s_{l}(i))}{2})\\
=\sum_{i=1}^{n}\sum_{1\leq k< l \leq N}\frac{b^{\prime}_{i}(s_{k}(i)-s_{l}(i))^{2}}{a_{i}},
\end{multline*}
where $b_{i}^{\prime}$ is chosen in such a way that 
\[b_{i}^{\prime}b_{i}\equiv 1 \hspace{2em}\operatorname{mod} 2a_{j},\] 
and $b_{i}$ is chosen as in Lemma \ref{FL}. Now observe the equality
\[a_{1}\cdots a_{n}(-b+\Sigma(\frac{b_{i}}{a_{i}}))=1\]
implies
\begin{equation}\label{Eq:Formulafore}
a+a_{i}^{2}m=a_{i}b_{i}^{\prime},
\end{equation} 
for some $m\in\mathbb{Z}$. Multiply the denominator and numerator by $a_{i}$:
\[\sum_{i=1}^{n}\sum_{1\leq k< l \leq N}\frac{b^{\prime}_{i}(s_{k}(i)-s_{l}(i))^{2}}{a_{i}}=\sum_{i=1}^{n}\sum_{1\leq k< l \leq N}\frac{a_{i}b^{\prime}_{i}(s_{k}(i)-s_{l}(i))^{2}}{a_{i}^{2}}\]
and plug equation \eqref{Eq:Formulafore} into the summation on the right hand side, we obtain    
\[
\sum_{i=1}^{n}\sum_{1\leq k< l \leq N}\frac{a_{i}b^{\prime}_{i}(s_{k}(i)-s_{l}(i))^{2}}{a_{i}^{2}}
\equiv \sum_{i=1}^{n}\sum_{1\leq k< l \leq N}\frac{a(s_{k}(i)-s_{l}(i))^{2}}{a_{i}^{2}}\hspace*{1em}\operatorname{mod} \mathbb{Z}.
\] 
The assertion then follows from Theorem \ref{Xiandeidentified1} as we can identify the $e$-invariant with the $\tilde{\xi}$-invariant in this case.
\end{proof}

\begin{remark}\label{ThmCDEinJW}
The formula above and the one given in Theorem $C$ in \cite{JW} are basically the same except that, in our formula, there is no minus sign. That means, with $\tilde{\xi}$-invariants of lens spaces and cobordisms given in this paper, we reobtain, up to sign, the formula for the real part of $e$-invariants of $\Sigma(a_{1},..,a_{n})$ in \cite[Theorem $C$]{JW}. Consequently, all the results following Theorem $C$ in \cite{JW}---Theorem $D$ and Theorem $E$ in \cite{JW}---remain intact (except for a sign change).     
\end{remark}
 




\begin{appendices}
\section{The stable homotopy category}\label{SHC}
In this appendix, we shall  review some basic constructions in the category of prespectra $\mathcal{P}$ and explain how the associated homotopy category is isomorphic to the homotopy category of the Adams category $\mathcal{A}$. First we give a digest of the model structure on $\mathcal{P}$ and the homotopy fiber and cofiber constructions in $\mathcal{P}$. Then we review the suspension and loop functors in $\mathcal{P}$ and explain how $\mathcal{P}$ models the stable homotopy category. We shall also see the infinite suspension and infinite loop functors constitute a Quillen adjunction
\[\Sigma^{\infty}:\mathpzc{Top}_{\ast}\leftrightarrows \mathcal{P}:\Omega^{\infty}.\] 
Lastly, we investigate the relation between $\operatorname{Ho}(\mathcal{A})$ and $\operatorname{Ho}(\mathcal{P})$, giving a detailed account of how they are isomorphic. 

The index set used here is $\mathbb{N}\cup\{0\}$ unless otherwise specified.

\subsection{The category of prespectra}\label{CPre}
Most of the results presented in this subsection are taken from \cite{nLPre} and \cite[Part II]{MMSS}. In \cite[Sec.2]{BF}, the simplicial analogue of $\mathcal{P}$ is studied.

\subsubsection*{The model structure}

\begin{definition}[Prespectra]\label{Pres}

\hfill\break
\begin{enumerate}
\item
A prespectrum is a sequence of spaces with structure maps 
\[\{E_{i},\sigma_{i}:S^{1}\wedge E_{i}\rightarrow E_{i+1}\}_{i\in\mathbb{N}\cup\{0\}}.\]
When the adjoint of $\sigma_{i}$
\[\tilde{\sigma}_{i}: E_{i}\rightarrow\Omega E_{i+1}\]
is a weak homotopy equivalence, it is called an $\Omega$-prespectrum. If each component $E_{i}$ is a $CW$-complex and the structure maps are inclusions of subcomplexes, then it is called a $CW$-prespectrum. We usually drop the structure maps from the notation when there is no risk of confusion. 
 
\item 
A map of prespectra $\bm{f}:\mathbf{E}=\{E_{i},\sigma_{i}\} \rightarrow \mathbf{F}=\{F_{i},\sigma^{\prime}_{i}\}$ is a sequence of continuous maps $f_{i}:E_{i}\rightarrow F_{i}$ such that 
\[\sigma^{\prime}_{i}\circ f_{i}=f_{i+1}\circ \sigma_{i},\]
for $i\in\mathbb{N}\cup\{0\}$.  

\item A map of prespectra $\bm{f}:\mathbf{E}\rightarrow\mathbf{F}$ is called a $\pi_{\ast}$-isomorphism if and only if the induced homomorphism 
\[\bm{f}_{\ast}:\pi_{\ast}(\mathbf{E})\rightarrow \pi_{\ast}(\mathbf{F})\]
is an isomorphism, for $\ast\in \mathbb{Z}$, where 
\[\pi_{\ast}(\mathbf{E}):=\operatorname*{colim}_{i;i+\ast\geq 0}\pi_{i+\ast}(E_{i}).\]

\end{enumerate}

We denote the category of prespectra by $\mathcal{P}$.

\end{definition}

\begin{theorem}(\cite{MMSS} and \cite[Theorem 0.70]{nLPre})\label{ModelcategoryofPres}
There is a model structure on the category of prespectra such that 

\begin{enumerate}

\item A map is a weak equivalence if and only if it is a $\pi_{\ast}$-isomorphism.

\item A map $\bm{p}:\mathbf{E}\rightarrow\mathbf{B}$ is a fibration if and only if $p_{i}:E_{i}\rightarrow B_{i}$ is a fibration and the following commutative diagram is a homotopy pullback in $\mathpzc{Top}_{\ast}$:

\begin{center}
\begin{tikzpicture}
\node(Lu) at (0,2) {$E_{i}$};
\node(Ll) at (0,0) {$B_{i}$}; 
\node(Ru) at (2,2) {$\Omega E_{i+1}$};
\node(Rl) at (2,0) {$\Omega B_{i+1}$};

\path[->, font=\scriptsize,>=angle 90]

(Lu) edge  (Ru)  
(Lu) edge node [right]{$p_{i}$}(Ll)
(Ll) edge  (Rl) 
(Ru) edge node [right]{$\Omega p_{i+1}$}(Rl);

\end{tikzpicture}
\end{center}
In particular, a prespectrum is fibrant if and only if it is an $\Omega$-prespectrum.

\item A map of prespectra $\mathbf{A}\rightarrow \mathbf{E}$ is a cofibration if and only if $A_{0}\rightarrow E_{0}$  and 
\[A_{n+1}\cup_{S^{1}\wedge A_{n}}S^{1}\wedge E_{n}\rightarrow E_{n+1}\]
are cofibrations in $\mathpzc{Top}_{\ast}$, for every $n$.
Especially, if a prespectrum $\mathbf{E}=\{E_{i},\sigma_{i}\}$ is cofibrant then $E_{i}$ has the homotopy type of $CW$-complex and $\sigma_{i}$ is a cofibration in $\mathpzc{Top}$. Thus every $CW$-prespectrum is a cofibrant object in this model structure on $\mathcal{P}$. 
  
\end{enumerate} 
\end{theorem}
\begin{proof} 
we only sketch the construction of this model structure and refer the interested readers to \cite[Sec.6-11]{MMSS} and \cite{nLPre} for detailed proofs. \cite{nLPre} contains more details concerning this particular case, while the approach in \cite{MMSS} are more general. 
  Recall the Bousfield-Friedlander theorem \cite[Proposition 0.66]{nLPre} says that, given a right proper model category $\mathcal{C}$ and a Quillen idempotent monad $(Q,\phi)$, which consists of an endofunctor $Q:\mathcal{C}\rightarrow \mathcal{C}$ and a natural transformation $\phi:id\mapsto Q$ (see \cite[Definition 0.62]{nLPre}), one can define a new model structure, denoted by $\mathcal{C}_{Q}$, by demanding a morphism $f$ is a weak equivalence in $\mathcal{C}_{Q}$ if and only if $Q(f)$ is, a cofibration if and only if it is in $\mathcal{C}$, and a fibration if and only if it satisfies the right lifting property against trivial cofibrations in $\mathcal{C}_{Q}$.  

Now observe there is a strict model structure on $\mathcal{P}$ given by level-wise weak homotopy equivalences and level-wise Serre fibrations (see \cite[Definition 0.38]{nLPre}). We also have a Quillen idempotent monad induced by the $\Omega$-prespectrification functor 
\[\bm{Q}:\mathcal{P}\rightarrow \mathcal{P}\] 
and a natural transformation $\bm{\phi}:id\mapsto \bm{Q}$, where we have $\bm{Q}\mathbf{E}$ is an $\Omega$-prespectrum and 
\[\bm{\phi}_{\mathbf{E}}:\mathbf{E}\rightarrow Q\mathbf{E}\] a $\pi_{\ast}$-isomorphism (see \cite[Definition 0.19]{nLPre}). Applying the Bousfield-Friedlander theorem, we obtain the model structure required in the theorem. 

For the description of the cofibrations and fibrations in this model structure, we refer to \cite[Definition 0.38, Definition 0.60]{nLPre} or \cite[Proposition 9.5, Lemma 11.4]{MMSS}. 
\end{proof}
From now on, $\mathcal{P}$ denotes the model category of prespectra described in the theorem above.

\subsubsection*{Homotopy fiber and cofiber}\vspace*{.5em}
\textbf{Homotopy fiber:}
Given a map of prespectra $\bm{f}:\mathbf{E}\rightarrow\mathbf{F}$, the homotopy fiber $\mathbf{Fib}(\bm{f})$ is the prespectrum whose $n$-th component is given by
\[\operatorname{Fib}(f)_{n}:=E_{n}\times_{f_{n}}PF_{n},\] 
where $PF_{n}$ is the mapping space $\mathpzc{Top}_{\ast}(I,F_{n})$ and $0$ is the base point of $I=[0,1]$. This gives us a sequence of prespectra 
\[\mathbf{Fib}(\bm{f})\rightarrow\mathbf{E}\xrightarrow{\bm{f}} \mathbf{F}.\]
Furthermore, if $\mathbf{E}$ and $\mathbf{F}$ are $\Omega$-prespectra, then so is $\mathbf{Fib}(\bm{f})$.

\textbf{Homotopy cofiber:}
Given a map of prespectra $\bm{f}:\mathbf{E}\rightarrow\mathbf{F}$, the homotopy cofiber $\mathbf{Cofib}(f)$ is the prespectra whose $n$-th component is the space $E_{n}\cup_{f_{n}}CF_{n}$ where $CF_{n}=F_{n}\wedge I_{+}$. This also yields a sequence of prespectra 
\[\mathbf{E}\xrightarrow{\bm{f}} \mathbf{F}\rightarrow \mathbf{Cofib}(\bm{f}),\]
and if $\mathbf{E}$ and $\mathbf{F}$ are $CW$-prespectra, then $\mathbf{Cofib}(\bm{f})$ is also a $CW$-prespectrum. These two sequences of prespectra are related by the following lemma (see \cite[p.128-130]{LMS} for the proof):

\begin{lemma}
There is a $\pi_{\ast}$-isomorphism, for every map of prespectra $\bm{f}$,  
\[\mathbf{Fib}(\bm{f})\xrightarrow{\bm{\zeta_{f}}} \Omega^{std}\mathbf{Cofib}(\bm{f})\] 
such that $\bm{\zeta_{f}}$ and its adjoint $\bm{\tilde{\zeta}_{f}}:\Sigma^{std}\mathbf{Fib}\rightarrow \mathbf{Cofib}(\bm{f})$ fit into the following diagram:

\begin{center}
\begin{tikzpicture}
\node (U1) at (3,3){$\Sigma^{std}\Omega^{std} \mathbf{E}$};
\node (U2) at (6,3){$\Sigma^{std}\Omega^{std} \mathbf{F}$};
\node (U3) at (9,3){$\Sigma^{std}\mathbf{Fib}(\bm{f})$};
\node (M0) at (0,1.5){$\mathbf{Fib}(\bm{f})$};
\node (M1) at (3,1.5){$\mathbf{E}$};
\node (M2) at (6,1.5){$\mathbf{F}$};
\node (M3) at (9,1.5){$\mathbf{Cofib}(\bm{f})$};
\node (L0) at (0,0){$\Omega^{std}\mathbf{Cofib}(\bm{f})$};
\node (L1) at (3,0){$\Omega^{std}\Sigma^{std} \mathbf{E}$};
\node (L2) at (6,0){$\Omega^{std}\Sigma^{std} \mathbf{F}$};

\path[->, font=\scriptsize,>=angle 90] 

(U1) edge (U2)
(U2) edge (U3)
(M0) edge (M1)
(M1) edge (M2)
(M2) edge (M3)
(L0) edge (L1)
(L1) edge (L2)

(U1) edge node [right]{$\wr$}(M1)
(U2) edge node [right]{$\wr$}(M2)
(U3) edge node [right]{$\wr$} node [left]{$\bm{\tilde{\zeta}_{f}}$}(M3)
(M0) edge node [right]{$\wr$} node [left]{$\bm{\zeta_{f}}$}(L0)
(M1) edge node [right]{$\wr$}(L1)
(M2) edge node [right]{$\wr$}(L2); 

\end{tikzpicture}
\end{center}
where $\Sigma^{std}$ and $\Omega^{std}$ are the standard suspension and loop functors (see \eqref{Standardsusloop} for the definition).
\end{lemma}

In fact, the map $\mathbf{Fib}(\bm{f})\rightarrow \Omega^{std}\mathbf{Cofib}(\bm{f})$ is natural with respect to $\bm{f}$. 

\begin{lemma}\label{NaturalisoFibandCofib}
Let $\operatorname{Ar}\mathcal{P}$ be the category of maps of prespectra and $\mathbf{Fib}$ and $\mathbf{Cofib}$ the functors from $\operatorname{Ar}\mathcal{P}$ to $\mathcal{P}$ induced by the homotopy fiber and cofiber constructions, respectively. Then the map $\zeta_{f}$ gives a natural transformation 
\[\bm{\zeta}:\mathbf{Fib}\longmapsto \Omega^{std}\mathbf{Cofib}.\]
\end{lemma}
\begin{proof}
Recall the map 
\[\bm{\zeta_{f}}:\mathbf{Fib}(\bm{f})\rightarrow \Omega^{std}(\mathbf{Cofib}(\bm{f}))\] 
is constructed level-wisely. That is given $f:X\rightarrow Y$ a map of spaces, $\zeta_{f}$ is defined as follows
\begin{align*}
\operatorname{Fib}(f)=X\times_{f}(Y,\ast)^{(I,0)}&\rightarrow \Omega\operatorname{Cofib}(f)=\Omega(Y\cup_{f} X\wedge I_{+})\\
(x,\gamma(t))&\mapsto \lambda(t),\\
\text{where   } 
\lambda(t)&=\begin{cases}
\gamma(2t)& t\leq \frac{1}{2},\\
(x,2t-1)& t\geq \frac{1}{2}.\end{cases}
\end{align*} 
With this construction, it is not difficult to see given a commutative square:
\begin{center}
\begin{tikzpicture}
\node(Lu) at (0,2) {$X$};
\node(Ll) at (0,0) {$X^{\prime}$}; 
\node(Ru) at (2,2) {$Y$};
\node(Rl) at (2,0) {$Y^{\prime}$};

\path[->, font=\scriptsize,>=angle 90]

(Lu) edge node [above]{$f$}(Ru)  
(Lu) edge (Ll)
(Ll) edge node [above]{$f^{\prime}$}(Rl) 
(Ru) edge (Rl);
\end{tikzpicture}
\end{center}
The following induce square is also commutative
\begin{center}
\begin{tikzpicture}
\node(Lu) at (0,2) {$\operatorname{Fib}(f)$};
\node(Ll) at (0,0) {$\operatorname{Fib}(f^{\prime})$}; 
\node(Ru) at (3,2) {$\operatorname{Cofib}(f)$};
\node(Rl) at (3,0) {$\operatorname{Cofib}(f^{\prime})$};

\path[->, font=\scriptsize,>=angle 90]

(Lu) edge node [above]{$\zeta_{f}$}(Ru)  
(Lu) edge (Ll)
(Ll) edge node [above]{$\zeta_{f^{\prime}}$}(Rl) 
(Ru) edge (Rl);
\end{tikzpicture}
\end{center}
Hence, $\bm{\zeta}$ is natural.
\end{proof}

The next lemma verify these two constructions indeed give homotopy (co)fiber sequences in the model category $\mathcal{P}$. Let us first recall the definition of homotopy fiber sequences in a model category $\mathcal{M}$.  

\begin{definition}\label{Hocartesian}
Given a commutative diagram in $\mathcal{M}_{f}$, the subcategory of fibrant objects in $\mathcal{M}$:

\begin{center}
\begin{tikzpicture}
\node (Lu) at (0,2){$X$};
\node (Ll) at (0,0){$W$};
\node (Rl) at (2,0){$Z$};
\node (Ru) at (2,2){$Y$};

\path[->, font=\scriptsize,>=angle 90] 

(Lu) edge  (Ru)
(Lu) edge  (Ll)
(Ru) edge node[right] {$f$}(Rl)
(Ll) edge node[above] {$g$}(Rl);
 
\end{tikzpicture}
\end{center}
We say it is a homotopy cartesian square if and only if
there is a factorization $f=p\circ i$ with $p:\tilde{Y}\rightarrow Z$ a fibration and $i:Y\rightarrow \tilde{Y}$ a weak equivalence such that the canonical morphism from $X$ to a pullback of $W\rightarrow Z\leftarrow \tilde{Y}$ (dashed arrow) is a weak equivalence. This can be illustrated as follows:   
\begin{center}
\begin{tikzpicture}
\node (LLuu) at (0,3){$X$};
\node (Lu) at (1,2){$W\times_{Z}\tilde{Y}$};
\node (Ll) at (1,0){$W$};
 
\node (Ruu) at (3,3){$Y$};
\node (Ru) at (3,2){$\tilde{Y}$};

\node (Rl) at (3,0){$Z$}; 

\path[->, font=\scriptsize,>=angle 90] 

(LLuu) edge  (Ruu)
(Lu) edge  (Ll)
 
(Lu) edge (Ru)
(Ruu) edge node[right] {$i$}(Ru)
(Ru) edge node[right] {$p$}(Rl)
(Ll) edge node[above] {$g$}(Rl);

 \draw[dashed,->] (LLuu)--(Lu);
 \draw[->] (LLuu)to [out=-100,in=140] (Ll);
 \draw[->] (Ruu)to[out=-40,in=40] node [right]{$f$}(Rl);
 
\end{tikzpicture}
\end{center}
If $W\rightarrow \ast$ is a weak equivalence, where $\ast$ is the terminal object, then $X\rightarrow Y\rightarrow Z$ is called a homotopy fiber sequence.
The dual notion gives us the definitions of homotopy cocartesian squares and homotopy cofiber sequences. 
\end{definition}

\begin{lemma}\phantomsection\label{HofibercofiberseqinP}
Suppose $\mathbf{E}$ and $\mathbf{F}$ are $\Omega$-prespectra, 
then the sequences of prespectra 
\[\mathbf{Fib}(\bm{f})\rightarrow\mathbf{E}\xrightarrow{\bm{f}} \mathbf{F}\]
is a homotopy fiber sequence in $\mathcal{P}$.

Suppose $\mathbf{E}$ and $\mathbf{F}$ are $CW$-prespectra, 
then the sequence of prespectra
\[\mathbf{E}\xrightarrow{\bm{f}} \mathbf{F}\rightarrow \mathbf{Cofib}(\bm{f})\]
is a homotopy cofiber sequence in $\mathcal{P}$.

\end{lemma}
\begin{proof} 
 

By Theorem \ref{ModelcategoryofPres} and the construction of homotopy fiber $\mathbf{Fib}(\bm{f})$, we know the map 
\[P\mathbf{F}\rightarrow \mathbf{F}\] 
is a fibration in $\mathcal{P}$, and $\mathbf{Fib}(\bm{f})$ is a pullback of 
\[P\mathbf{F}\rightarrow \mathbf{F}\leftarrow\mathbf{E}.\]
In particular, we have the commutative diagram

\begin{center}
\begin{tikzpicture}
\node(Lu) at (0,2) {$\mathbf{Fib}({\bm{f}})$};
\node(Ll) at (0,0) {$P\mathbf{F}$}; 
\node(Ru) at (2,2) {$\mathbf{E}$};
\node(Rl) at (2,0) {$\mathbf{F}$};

\path[->, font=\scriptsize,>=angle 90]

(Lu) edge (Ru)  
(Lu) edge (Ll)
(Ll) edge (Rl) 
(Ru) edge (Rl);
\end{tikzpicture}
\end{center}
is a homotopy cartesian. Since $\ast\rightarrow P\mathbf{F}$ is a $\pi_{\ast}$-isomorphism, the sequence 
\[\mathbf{Fib}(\bm{f})\rightarrow \mathbf{E}\rightarrow \mathbf{F}\]
is a homotopy fiber sequence. 

Similarly, we have 
\[\mathbf{E}\rightarrow \mathbf{E}\wedge I_{+}\]
is a cofibration in $\mathcal{P}$ by Theorem \ref{ModelcategoryofPres}, and $\mathbf{Cofib}(\bm{f})$ is a pushout of the cospan 
\[\mathbf{F}\leftarrow \mathbf{E}\rightarrow \mathbf{E}\wedge I_{+}.\]
Therefore the commutative diagram

\begin{center}
\begin{tikzpicture}
\node(Lu) at (0,2) {$\mathbf{E}$};
\node(Ll) at (0,0) {$\mathbf{F}$}; 
\node(Ru) at (2,2) {$\mathbf{E}\wedge I_{+}$};
\node(Rl) at (2,0) {$\mathbf{Cofib}(\bm{f})$};

\path[->, font=\scriptsize,>=angle 90]

(Lu) edge (Ru)  
(Lu) edge (Ll)
(Ll) edge (Rl) 
(Ru) edge (Rl);
\end{tikzpicture}
\end{center}
is a homotopy cocartesian. Since $\mathbf{E}\wedge I_{+}$ is $\pi_{\ast}$-isomorphic to $\ast$, the sequence 
\[\mathbf{E}\rightarrow \mathbf{F}\rightarrow \mathbf{Cofib}(\bm{f})\]
is a homotopy cofiber sequence

\end{proof}

\begin{remark}
It seems the assumption $\mathbf{E}$ and $\mathbf{F}$ are $\Omega$-prespectra in the first assertion (resp. $\mathbf{E}$ and $\mathbf{F}$ are $CW$-prespectra in the second), is redundant as $\mathcal{P}$ might actually be a proper and simplicial model category as its simplicial analogue is (see \cite[Theorem $2.3$]{BF}). 
\end{remark}

\subsubsection*{Stability of $\mathcal{P}$}\label{StabilityofP}
Observe the infinite loop functor 
\[\Omega^{\infty}:\mathcal{P}\rightarrow  \mathpzc{Top}_{\ast}\]
which associates a prespectra $\mathbf{E}$ with its zero component $E_{0}\in \mathpzc{Top}_{\ast}$ and the infinite suspension functor   
\[\Sigma^{\infty}: \mathpzc{Top}_{\ast}\rightarrow \mathcal{P}.\]
constitute a Quillen adjunction
\[\Sigma^{\infty}: \mathpzc{Top}_{\ast}\leftrightarrows \mathcal{P}:\Omega^{\infty}\]
(see \cite{nLPre} for more details).
Recall also that, given a pointed model category $\mathcal{M}$, we have the (canonically) induced suspension functor 
\[\Sigma: \operatorname{Ho}(\mathcal{M})\rightarrow \operatorname{Ho}(\mathcal{M}),\]
and the (canonically) induced loop functor
\[\Omega: \operatorname{Ho}(\mathcal{M})\rightarrow \operatorname{Ho}(\mathcal{M}),\]
which are a homotopy pushout of $\ast\leftarrow C \rightarrow \ast$, and a homotopy pullback of $\ast\rightarrow C\leftarrow \ast$, respectively. $\Omega$ and $\Sigma$ are uniquely determined up to isomorphisms in $\operatorname{Ho}(\mathcal{M})$.

A model category $\mathcal{M}$ is a stable model category if and only if the induced suspension functor $\Sigma$ and loop functor $\Omega$ are inverse equivalences on $\operatorname{Ho}(\mathcal{M})$.

\textbf{Standard suspension and loop functors:}
Now in the model category of prespectra, the induced suspension functor and loop functor can be realized by the standard suspension $\Sigma^{std}$ and loop functor $\Omega^{std}$ constructed as below. Let $\mathbf{E}:=\{E_{k},\sigma_{k}\}$ and $\tilde{\sigma}_{k}: E_{k}\rightarrow \Omega E_{k+1}$ be the adjoint of $\sigma_{k}$. Then we define
\begin{align} 
\Sigma^{std}\mathbf{E}&:=\{E_{k}\wedge S^{1}, \sigma_{k}\wedge S^{1}\}\nonumber \\
\Omega^{std}\mathbf{E}&:=\{\mathpzc{Top}_{\ast}(S^{1},E_{k}),\sigma^{\prime}_{k}\}, \label{Standardsusloop}  
\end{align}
where  
\begin{multline*}
\sigma^{\prime}_{k}:= S^{1}\wedge \mathpzc{Top}_{\ast}(S^{1},E_{k})\xrightarrow{(c,id)} \mathpzc{Top}_{\ast}(S^{1},S^{1}\wedge E_{k})\\
 \xrightarrow{\mathpzc{Top}_{\ast}(S^{1},\sigma_{k})} \mathpzc{Top}_{\ast}(S^{1},E_{k+1}),
\end{multline*}
and $(c,id)(t,\phi):=(c_{t},\phi)$ with $(c_{t},\phi)(s):=(t,\phi(s))$ for $s\in S^{1}$.
These constitute a Quillen adjunction
\[\Sigma^{std}:\mathcal{P}\leftrightarrows \mathcal{P}: \Omega^{std}.\]
One can prove this by a similar argument used in \cite[Lemma 0.72]{nLPre}---although the lemma there is for the ``alternative" suspension which is different from the standard suspension, the approach can be applied to this case.

On the other hand, we have the $k$-fold shifting functor which associates a prespectrum $\mathbf{E}$ with another prespectrum $\mathbf{E}[k]$ whose $n$-th component is $E_{n+k}$ when $n+k\geq 0$, and $\ast$ otherwise. These readily give us the following Quillen equivalence
\[(-)[-1]:\mathcal{P}\leftrightarrows \mathcal{P}:(-)[1]\]

\textbf{Alternative suspension and loop functors:} 
Even though theoretically the construction of the standard suspension functor $\Sigma^{std}$ is the natural one, it suffers some disadvantages. For example, it is hard to compare $\Sigma^{std}$ with the shifting functor $(-)[1]$. For this reason, the alternative suspension and loop functors, denoted by $\Sigma^{alt}$ and $\Omega^{alt}$ are introduced. They are given as follows: 
\begin{align*}
\Sigma^{alt}\mathbf{E}&:=\{S^{1}\wedge E_{k},S^{1}\wedge \sigma_{k}\}\\
\Omega^{alt}\mathbf{E}&:=\{\mathpzc{Top}_{\ast}(S^{1},E_{k}),\mathpzc{Top}_{\ast}(S^{1},\tilde{\sigma}_{k})\} 
\end{align*}
where $\mathbf{E}:=\{E_{k},\sigma_{k}:S^{1}\wedge E_{k}\rightarrow E_{k+1}\}$ is a prespectrum and $\tilde{\sigma}_{k}$ is the adjoint of $\sigma_{k}$. One can check they constitute a Quillen adjunction
\[\Sigma^{alt}:\mathcal{P}\leftrightarrows \mathcal{P}:\Omega^{alt},\]
and there are natural transformations
\begin{align*}
\Sigma^{alt}(-) &\mapsto (-)[1]\\
(-)[-1]         &\mapsto \Omega^{alt}(-).
\end{align*}  
such that $\Sigma^{alt}\mathbf{E}\rightarrow \mathbf{E}[1]$ and 
$\mathbf{E}[-1]\rightarrow \Omega^{alt}\mathbf{E}$ are $\pi_{\ast}$-isomorphisms, for every $\mathbf{E}\in\mathcal{P}$. Therefore, $\Sigma^{alt}$ and $\Omega^{alt}$ actually constitute a Quillen equivalence
\[\Sigma^{alt}:\mathcal{P}\leftrightarrows \mathcal{P}:\Omega^{alt}.\]

\textbf{The relation between $\Sigma^{std}$ and $\Sigma^{alt}$:}
Though there seems no direct way to compare $\Sigma^{std}$ and $\Sigma^{alt}$---the most natural map $E_{k}\wedge S^{1}\rightarrow S^{1}\wedge E_{k}$ does not give us a map in $\mathcal{P}$ as it is not compatible with the structure maps of the prespectra (see the explanation in \cite[Remark 0.34]{nLPre}), it is possible to compare them in $\operatorname{Ho}(\mathcal{P})$. In fact, it is proved that there is a natural transformation 
from $\Sigma^{std}$ to $\Sigma^{alt}$ in $\operatorname{Ho}(\mathcal{P})$ \cite[Lemma 0.81]{nLPre}. Notice both of them, as Quillen left adjoints, descend to functors of the homotopy category $\operatorname{Ho}(\mathcal{P})$. In this way, we see $\mathcal{P}$ is a stable model category. We summarize the discussion above is the following proposition: 
 
\begin{proposition}\phantomsection\label{QadjunctionPresandspace}
\begin{enumerate}
\item

There is a diagram of Quillen's adjunctions (below) with functors at the bottom are Quillen equivalences.

\begin{center}
\begin{tikzpicture}[scale=1.2]
\node(Lu) at (0,2) {$\mathpzc{Top}_{\ast}$};
\node(Ll) at (0,0) {$\mathcal{P}$}; 
\node(Ru) at (4,2) {$\mathpzc{Top}_{\ast}$};
\node(Rl) at (4,0) {$\mathcal{P}$};
\node(l) at (0,1){\tiny $\dashv$};
\node(r) at (4,1){\tiny $\dashv$};
\node(u) at (2,2){\tiny $\perp$}; 
\node(u) at (2,0){\tiny $\perp$};

\draw[transform canvas={yshift=1ex},->] (Lu) --(Ru) node[above,midway] {$\Sigma$}  ;
\draw[transform canvas={yshift=-1ex},->](Ru) --(Lu) node[below,midway] {$\Omega$};

\draw[transform canvas={yshift=1ex},->] (Ll) --(Rl) node[above,midway] {$\Sigma^{alt/std},\sim$};
\draw[transform canvas={yshift=-1ex},->](Rl) --(Ll) node[below,midway] {$\Omega^{alt/std},\sim$};

\draw[transform canvas={xshift=-1ex},->] (Ll) --(Lu) node[left,midway] {$\Sigma^{\infty}$};
\draw[transform canvas={xshift=1ex},->] (Lu) --(Ll) node[right,midway] {$\Omega^{\infty}$};

\draw[transform canvas={xshift=-1ex},->] (Rl) --(Ru) node[left,midway] {$\Sigma^{\infty}$} ;
\draw[transform canvas={xshift=1ex},->] (Ru) --(Rl) node[right,midway] {$\Omega^{\infty}$};
 
\end{tikzpicture}
\end{center}

\item
The model category $\mathcal{P}$ is a stable model category as we have
\[\Sigma=\Sigma^{std}:\operatorname{Ho}(\mathcal{P})\leftrightarrows \operatorname{Ho}(\mathcal{P}):\Omega^{std}=\Omega\] are inverse equivalences. Especially, every homotopy fiber sequence is a homotopy cofiber sequence and vice verse. 
The distinguished triangles in $\operatorname{Ho}(\mathcal{P})$ are the closure of homotopy cofiber sequences under isomorphisms, and there is the long exact sequence 
\begin{multline*}
...\rightarrow [\mathbf{X},\Omega \mathbf{B}]_{\operatorname{Ho}(\mathcal{P})}\rightarrow [\mathbf{X},\mathbf{F}]_{\operatorname{Ho}(\mathcal{P})}\rightarrow [\mathbf{X},\mathbf{E}]_{\operatorname{Ho}(\mathcal{P})}\\
\rightarrow [\mathbf{X},\mathbf{B}]_{\operatorname{Ho}(\mathcal{P})} \rightarrow [\mathbf{X},\Sigma \mathbf{F}]_{\operatorname{Ho}(\mathcal{P})}\rightarrow...,
\end{multline*}
for every homotopy (co)fiber sequence 
\[\mathbf{F}\rightarrow \mathbf{E}\rightarrow \mathbf{B}\] and any $\mathbf{X}\in\operatorname{Ho}(\mathcal{P})$.

\end{enumerate}
\end{proposition}

\begin{proof}
The discussion preceding the theorem is an outline. A detailed proof for the first assertion can be found in \cite[Theorem 0.84]{nLPre}, while for the second assertion, we refer to \cite[Chap.6-7]{Ho} and \cite[Prop.0.110,0.116]{nLPre}. 
\end{proof} 


\subsection{Relation with the Adams category}
\label{AandPcomparison} 
Here we describe how the Adams category $\mathcal{A}$ is related to the category of prespectra $\mathcal{P}$. To distinguish, we reserve the term ``$CW$-spectra" for objects in $\mathcal{A}$ while using the term ``$CW$-prespectra" for the $CW$-objects in $\mathcal{P}$. We also use $[-,-]_{\operatorname{Ho}(\mathcal{A})}$ and $[-,-]_{\operatorname{Ho}(\mathcal{P})}$ to distinguish the homotopy classes of maps in $\mathcal{A}$ and $\mathcal{P}$. Recall maps of $CW$-spectra are eventually-defined maps and a map is a $\pi_{\ast}$-isomorphism if and only if it is a homotopy equivalence (see \cite[Chap.8]{Sw} for the definitions). Note also, without loss of generality, we can assume $CW$-spectra are indexed over $\mathbb{N}\cup\{0\}$.

\begin{lemma}
For every $CW$-spectrum $\mathbf{E}$ in $\mathcal{A}$, one can construct an $\Omega$-$CW$-spectrum $\bar{\mathbf{E}}$ such that there is a homotopy equivalence of $CW$-spectra 
\[\bm{\phi}_{\mathbf{E}}:\mathbf{E}\rightarrow \bar{\mathbf{E}}.\]
\end{lemma}
\begin{proof}
When $\mathbf{E}$ is already an $\Omega$-$CW$-spectrum, we let $\bar{\mathbf{E}}$ be $\mathbf{E}$ and $\bm{\phi}$ identity. Otherwise, we use the construction described in \cite[Chap II, Prop.1.21]{Ru}. Alternatively, one can first use the functor $\bm{Q}$ in \cite[Definition 0.19]{nLPre}, and then apply the $CW$-approximation as in \cite[II, Prop. 1.21]{Ru}.
\end{proof}

With this lemma, we obtain
\begin{proposition}
The inclusion 
\[\operatorname{Ho}(\Omega\mathcal{A})\hookrightarrow \operatorname{Ho}(\mathcal{A})\]
is an equivalence of the categories, where $\Omega\mathcal{A}$
is the full subcategory of $\Omega$-$CW$-spectra.
\end{proposition}
\begin{proof}
The functor in the other direction can be constructed by the lemma above. In more details, it sends a $CW$-spectrum $\mathbf{E}$ to $\bar{\mathbf{E}}$ and a map from $\mathbf{E}$ to $\mathbf{F}$ to a map from $\bar{\mathbf{E}}$ to $\bar{\mathbf{F}}$ via the isomorphism
\[[\mathbf{E},\mathbf{F}]_{\operatorname{Ho}(\mathcal{A})}\xrightarrow{\sim} [\bar{\mathbf{E}},\bar{\mathbf{F}}]_{\operatorname{Ho}(\mathcal{A})},\]
which is induced by the following cospan of isomorphisms
\[[\bar{\mathbf{E}},\bar{\mathbf{F}}]_{\operatorname{Ho}(\mathcal{A})}\xrightarrow{\sim} [\mathbf{E},\bar{\mathbf{F}}]_{\operatorname{Ho}(\mathcal{A})}\xleftarrow{\sim} [\mathbf{E},\mathbf{F}]_{\operatorname{Ho}(\mathcal{A})}.\]
In this way, we see the homotopy equivalence 
\[\bm{\phi}_{\mathbf{E}}:\mathbf{E}\rightarrow \bar{\mathbf{E}}\]
in the preceding lemma is natural and induces the natural isomorphism required.
\end{proof}

On the other hand, via the model structure of $\mathcal{P}$, we have the following observation
\begin{proposition}
The inclusion 
\[\operatorname{Ho}(\mathcal{P}_{cf})\hookrightarrow \operatorname{Ho}(\mathcal{P})\]
induces an equivalence of categories, where $\mathcal{P}_{cf}$ is the full subcategory of cofibrant-fibrant objects in $\mathcal{P}$.
\end{proposition}
\begin{proof}
This follows from the definition of the homotopy category of a model category.
\end{proof}

The next lemma about $CW$-approximation is proved in \cite{nLPre}.
\begin{lemma}
For any prespectrum $\mathbf{E}$, one can construct a $CW$-prespectrum $\hat{\mathbf{E}}$ and a map of prespectra 
\[\bm{\nu}_{\mathbf{E}}:\hat{\mathbf{E}}\rightarrow \mathbf{E},\]
such that it is a level-wise weak homotopy equivalence. Moreover, if $\mathbf{E}$ is already a $CW$-prespectrum, one can choose $\hat{\mathbf{E}}=\mathbf{E}$ and $\bm{\nu}_{\mathbf{E}}=id_{\mathbf{E}}$.
\end{lemma}
\begin{proof}
By induction and the $CW$-approximation theorem (see \cite[Prop.0.95]{nLPre}).
\end{proof}

\begin{theorem}\phantomsection\label{RelationbetweenAandP}
There are inverse equivalences of categories 
\[\mathcal{F}:\operatorname{Ho}(\Omega\mathcal{A})\leftrightarrows \operatorname{Ho}(\mathcal{P}_{cf}):\mathcal{G}.\]
\end{theorem}

\begin{proof}
\textbf{Step 1} (Construction of the functor $\mathcal{F}$): Recall a map of $CW$-spectra $\mathbf{E}\rightarrow \mathbf{F}$ is a pair $(\mathbf{E}^{\prime},\bm{f}^{\prime})$ where $\mathbf{E}^{\prime}$ is a cofinal subspectrum and $\bm{f}^{\prime}:\mathbf{E}^{\prime}\rightarrow \mathbf{F}$ a map of prespectra. Because $\mathbf{E}^{\prime}$ is cofinal subspectrum of $\mathbf{E}$ and both of them are $CW$-prespectra, the following induced homomorphism
\[[\mathbf{E},\mathbf{F}]_{\operatorname{Ho}(\mathcal{P})}\xrightarrow{\theta} [\mathbf{E}^{\prime},\mathbf{F}]_{\operatorname{Ho}(\mathcal{P})}\]
is an isomorphism (see Remark \ref{Morenotionsonhomotopy} or \cite[Theorem 14.4.8]{MP}). Thus we can define $\mathcal{F}((\mathbf{E}^{\prime},\bm{f}^{\prime}))$ to be $\theta^{-1}(\mathbf{E}^{\prime},\bm{f}^{\prime})$. To see it is well-defined, we recall two maps of $CW$-spectra $(\mathbf{E}^{\prime},
 \bm{f}^{\prime});(\mathbf{E}^{\prime\prime},\bm{f}^{\prime\prime}):\mathbf{E}\rightarrow\mathbf{F}$ are homotopic if and only if 
there exists a pair $(\mathbf{E}^{\prime\prime\prime}\wedge I_{+},H)$ such that $\mathbf{E}^{\prime\prime\prime}$ is a cofinal subspectrum of $\mathbf{E}^{\prime}\cap \mathbf{E}^{\prime\prime}$ and a map of prespectra
\[\bm{H}:\mathbf{E}^{\prime\prime\prime}\wedge I_{+}\rightarrow \mathbf{F}\]
that restricts to $\bm{f}^{\prime}\vert_{\mathbf{E}^{\prime\prime\prime}\times\{0\}}$ and $\bm{f}^{\prime\prime}\vert_{\mathbf{E}^{\prime\prime\prime}\times\{1\}}$.
Then the following commutative diagram of isomorphisms 
\begin{center} 
\begin{tikzpicture}
\node (L)  at (0,0){$[\mathbf{E},\mathbf{F}]_{\operatorname{Ho}(\mathcal{P})}$};
\node (Mr) at (2,1){$[\mathbf{E}^{\prime},\mathbf{F}]_{\operatorname{Ho}(\mathcal{P})}$};
\node (Ml) at (2,-1){$[\mathbf{E}^{\prime\prime},\mathbf{F}]_{\operatorname{Ho}(\mathcal{P})}$};
\node (R)  at (4,0){$[\mathbf{E}^{\prime\prime\prime},\mathbf{F}]_{\operatorname{Ho}(\mathcal{P})}$};

\path[->, font=\scriptsize,>=angle 90] 

(L)  edge (Mr)
(L)  edge (Ml)
(Mr) edge (R)
(Ml) edge (R);
\end{tikzpicture}
\end{center}
implies 
\[\mathcal{F}((\mathbf{E}^{\prime},\bm{f}^{\prime}))=\mathcal{F}((\mathbf{E}^{\prime\prime},\bm{f}^{\prime\prime})).\]

Suppose we are given $(\mathbf{E}^{\prime},\bm{f}^{\prime}):\mathbf{E}\rightarrow \mathbf{F}$ and $(\mathbf{F}^{\prime},\bm{g}^{\prime}):\mathbf{F}\rightarrow \mathbf{G}$. We want to show the following equality
\[\mathcal{F}((\mathbf{F}^{\prime},\bm{g}^{\prime})\circ(\mathbf{E}^{\prime},\bm{f}^{\prime}))= \mathcal{F}((\mathbf{F}^{\prime},\bm{g}^{\prime}))\circ \mathcal{F}((\mathbf{E}^{\prime},\bm{f}^{\prime})).\]
To see this, we first note, without loss of generality, one can assume $\bm{f}^{\prime}$ factors through $\mathbf{F}^{\prime}$. Then the identity above follows quickly from the commutative diagram below:
\begin{center}
\begin{tikzpicture}
\node(Lu) at (0,2) {$\mathbf{E}^{\prime}$};
\node(Ll) at (0,0) {$\mathbf{E}$};
\node(Mu) at (2,2) {$\mathbf{F}^{\prime}$};
\node(Ml) at (2,0) {$\mathbf{F}$};
\node(Ru) at (4,2) {$\mathbf{G}$};
\node(Rl) at (4,0) {$\mathbf{G}$};

\path[->, font=\scriptsize,>=angle 90] 

(Ll) edge node [above]{$\mathcal{F}(\mathbf{E}^{\prime},\bm{f}^{\prime})$}(Ml)
(Ml) edge node [above]{$\mathcal{F}(\mathbf{F}^{\prime},\bm{g}^{\prime})$}(Rl)
(Lu) edge node [above]{$\bm{f}^{\prime}$}(Mu)
(Mu) edge node [above]{$\bm{g}^{\prime}$}(Ru);

\draw [double equal sign distance] (Ru) to (Rl);  
\draw [right hook->] (Lu) to (Ll);
\draw [right hook->] (Mu) to (Ml); 
\end{tikzpicture}
\end{center}

\textbf{Step 2} (Construction of $\mathcal{G}$): The functor $\mathcal{G}$ sends a cofibrant $\Omega$-prespectrum $\mathbf{E}$ to $\hat{\mathbf{E}}$, the $\Omega$-$CW$-prespectrum given in the lemma preceding the theorem. Then observe following cospan of isomorphisms 
\[[\mathbf{E},\mathbf{F}]_{\operatorname{Ho}(\mathcal{P})}\xrightarrow{\bm{\nu}_{\mathbf{E}}^{\ast}}[\hat{\mathbf{E}},\mathbf{F}]_{\operatorname{Ho}(\mathcal{P})}\xleftarrow{\bm{\nu}_{\mathbf{F},\ast}}[\hat{\mathbf{E}},\hat{\mathbf{F}}]_{\operatorname{Ho}(\mathcal{P})}\]
and define $\mathcal{G}(\bm{f})$ to be $\bm{\nu}_{\mathbf{F},\ast}^{-1}\circ \bm{\nu}_{\mathbf{E}}^{\ast}(\bm{f})$ for any map of prespectra $\bm{f}:\mathbf{E}\rightarrow \mathbf{F}$. Note it is clear $\mathcal{G}(\bm{f})$ depends only on the homotopy class of $\bm{f}$. The following identity is also not hard to verify  
\[\mathcal{G}(\bm{g}\circ \bm{f})=\mathcal{G}(\bm{g})\circ\mathcal{G}(\bm{f})\]
as it ensues from the commutative diagram
\begin{center}
\begin{tikzpicture}
\node(Lu) at (0,2) {$\mathbf{E}$};
\node(Ll) at (0,0) {$\hat{\mathbf{E}}$};
\node(Mu) at (2,2) {$\mathbf{F}$};
\node(Ml) at (2,0) {$\hat{\mathbf{F}}$};
\node(Ru) at (4,2) {$\mathbf{G}$};
\node(Rl) at (4,0) {$\hat{\mathbf{G}}$};

\path[->, font=\scriptsize,>=angle 90] 

(Ll) edge node [above]{$\mathcal{G}(\bm{f})$}(Ml)
(Ml) edge node [above]{$\mathcal{G}(\bm{g})$}(Rl)
(Lu) edge node [above]{$\bm{f}$}(Mu)
(Mu) edge node [above]{$\bm{g}$}(Ru)
(Lu) edge node [right]{$\bm{\nu}_{\mathbf{E}}$} (Ll)
(Mu) edge node [right]{$\bm{\nu}_{\mathbf{F}}$} (Ml)
(Ru) edge node [right]{$\bm{\nu}_{\mathbf{G}}$} (Rl);  
\end{tikzpicture}
\end{center}

\textbf{Step 3:} The natural isomorphism from $\mathcal{F}\circ \mathcal{G}$ to $id$ is given by the level-wise homotopy equivalence 
\[\bm{\nu}_{\mathbf{E}}:\hat{\mathbf{E}}\rightarrow \mathbf{E},\] 
as the cospan of isomorphisms:
\[[\mathbf{E},\mathbf{F}]_{\operatorname{Ho}(\mathcal{P})}\xrightarrow{\bm{\nu}_{\mathbf{E}}^{\ast}}[\hat{\mathbf{E}},\mathbf{F}]_{\operatorname{Ho}(\mathcal{P})}\xleftarrow{\bm{\nu}_{\mathbf{F},\ast}}[\hat{\mathbf{E}},\hat{\mathbf{F}}]_{\operatorname{Ho}(\mathcal{P})},\]
implies that, given a map $\bm{f}:\mathbf{E}\rightarrow\mathbf{F}$, the following diagram commutes
\begin{center}
\begin{tikzpicture}
\node(Lu) at (0,2) {$\mathbf{E}$};
\node(Ll) at (0,0) {$\hat{\mathbf{E}}$}; 
\node(Ru) at (2,2) {$\mathbf{F}$};
\node(Rl) at (2,0) {$\hat{\mathbf{F}}$};

\path[->, font=\scriptsize,>=angle 90]

(Lu) edge node [above]{$\bm{f}$}(Ru)  
(Lu) edge node [right]{$\bm{\nu}_{\mathbf{E}}$}(Ll)
(Ll) edge node [above]{$\mathcal{F}\circ\mathcal{G}(\bm{f})$}(Rl) 
(Ru) edge node [right]{$\bm{\nu}_{\mathbf{F}}$}(Rl);

\end{tikzpicture}
\end{center}
Lastly, the identity assigning map provides the natural isomorphism between $\mathcal{G}\circ\mathcal{F}$ and $id$. More precisely, one assigns a $CW$-spectrum $\mathbf{E}$ to the identity map $id_{\mathbf{E}}\in [\mathbf{E},\mathbf{E}]_{\operatorname{Ho}(\mathcal{A})}$. 
To see why it works, we recall the definition of $\mathcal{F}((\mathbf{E}^{\prime},\bm{f}^{\prime}))$ which is the map of prespectra from $\mathbf{E}$ to $\mathbf{F}$ such that its image under the map
\[[\mathbf{E},\mathbf{F}]_{\operatorname{Ho}(\mathcal{P})}\rightarrow [\mathbf{E}^{\prime},\mathbf{F}]_{\operatorname{Ho}(\mathcal{P})}\]
is $\bm{f}^{\prime}$, and it is equivalent to say the following diagram commutes
\begin{center}
\begin{tikzpicture}
\node(Lu) at (0,2) {$\mathbf{E}^{\prime}$};
\node(Ll) at (0,0) {$\mathbf{E}$}; 
\node(Ru) at (2,2) {$\mathbf{F}$};
\node(Rl) at (2,0) {$\mathbf{F}$};

\path[->, font=\scriptsize,>=angle 90]

(Lu) edge node [above]{$\bm{f}^{\prime}$}(Ru)  
 
(Ll) edge node [above]{$\mathcal{G}\circ\mathcal{F}(\bm{f})$}(Rl) 
(Ru) edge node [right]{$id_{\mathbf{F}}$}(Rl);

\draw[right hook->] (Lu) to (Ll);
\end{tikzpicture}
\end{center}
Since the inclusion of a cofinal subspectrum represents the identity map in $\operatorname{Ho}(\mathcal{A})$, we see the identity-assigning map yields the desired natural isomorphism.

\end{proof}

\begin{remark}\phantomsection\label{Morenotionsonhomotopy}
Here we recall some notions about homotopy in a model category $\mathcal{M}$: Let $\operatorname{Hom}_{\mathcal{M}}(X,Y)$ be the set of morphisms from $X$ to $Y$ in $\mathcal{M}$, then one can use cylinder objects of $X$ to define left homotopy and use path objects of $Y$ to define right homotopy. They are, however, not equivalence relations in general. Nevertheless, if $X$ is cofibrant, left homotopy is an equivalence relation on $\operatorname{Hom}_{\mathcal{M}}(X,Y)$; and, if $Y$ is fibrant, right homotopy is an equivalence relation on $\operatorname{Hom}_{\mathcal{M}}(X,Y)$ (see \cite[Proposition 14.3.9]{MP}). We denote the left and right homotopy classes of morphisms from $X$ to $Y$ by $[X,Y]_{l}$ and $[X,Y]_{r}$, respectively. Then we have the following variant of Whitehead theorem (see \cite[Proposition 14.3.14]{MP}):
\begin{enumerate}
\item Given $p:Y\rightarrow Z$ an acyclic fibration and $X$ a cofibrant object, then $p_{\ast}:[X,Y]_{l}\rightarrow [X,Z]_{l}$ is a bijection.
\item Given $i:X\rightarrow Y$ an acyclic cofibration and $Z$ a fibrant object, then $i^{\ast}:[Y,Z]_{r}\rightarrow [X,Z]_{r}$ is a bijection.
\end{enumerate}

The notions of left and right homotopy are not the same in general, but when $X$ is cofibrant, we have $f,g:X\rightarrow Y$ are left homotopic implies $f,g$ are right homotopic; and when $Y$ is fibrant, $f,g$ are right homotopic implies $f,g$ are left homotopic (see \cite[Proposition 14.3.11]{MP}). Especially, if $X$ is cofibrant and $Y$ is fibrant, then left homotopy and right homotopy induce the same equivalence relation on $\operatorname{Hom}_{\mathcal{M}}(X,Y)$, and, furthermore, combing with the variant of Whitehead theorem above, one obtains the following bijections 
\begin{equation}\label{CFobjectandhomotopy}
[R(X),Q(Y)]_{r}=[R(X),Q(Y)]_{l}=:[X,Y]_{\operatorname{Ho}(\mathcal{M})}\xrightarrow{\sim} [X,Y]_{l}=[X,Y]_{r},
\end{equation}   
where $R(X)$ and $Q(Y)$ are the fibrant and cofibrant replacements of $X$ and $Y$, respectively. Apply this to the model category $\mathcal{P}$ and assume $\mathbf{E}$ is a $CW$-prespectrum and $\mathbf{F}$ an $\Omega$-$CW$-prespectrum. Then Theorem \ref{RelationbetweenAandP} implies the canonical homomorphism
\[[\mathbf{E},\mathbf{F}]_{l}\rightarrow [\mathbf{E},\mathbf{F}]_{\operatorname{Ho}(\mathcal{A})}\]
is an isomorphism. 
\end{remark}

 
\end{appendices}

{
\center \textbf{Acknowledgment}

This paper contains and improves some results of the author's thesis. The author would like to express his gratitude to his supervisor Sebastian Goette for his advice, guidance and patience. Thanks also to Ulrich Bunke, Wolfgang Steimle, and Matthias Wendt for their insightful comments and suggestions. The author thanks Jørgen Olsen Lye for reading part of the early draft carefully and the English corrections. The project is financially supported by the DFG Graduiertenkolleg 1821 "Cohomological Methods in Geometry".

}
\phantomsection 
\addcontentsline{toc}{section}{Bibliography}
    
\bibliographystyle{alpha}
 
\bibliography{Reference}

\begin{thebibliography}{CCMT84}

\bibitem[AAS14]{ASS}
P.~Antonini, S.~Azzali, and G.~Skandalis.
\newblock Flat bundles, von {N}eumann algebras and {K}-theory with
  $\mathbb{R/Z}$-coefficients.
\newblock {\em J. K-theory}, 13(2):275--303, 2014.

\bibitem[Ada71]{Ad}
J.~F. Adams.
\newblock A variant of {E}. {H}. {B}rown's representability theorem.
\newblock {\em Topology}, 1971.

\bibitem[Ada74]{Ad1}
J.~F. Adams.
\newblock {\em Stable Homotopy and Generalized Homology}.
\newblock Chicago Lectures in Mathematics. University of Chicago Press, 1974.

\bibitem[APS75a]{APS1}
M.~F. Atiyah, V.~K. Patodi, and I.~M. Singer.
\newblock Spectral asymmetry and riemannian geometry i.
\newblock {\em Math. Proc. Camb. Phil. Soc.}, 1975.

\bibitem[APS75b]{APS2}
M.~F. Atiyah, V.~K. Patodi, and I.~M. Singer.
\newblock Spectral asymmetry and riemannian geometry ii.
\newblock {\em Math. Proc. Camb. Phil. Soc.}, 1975.

\bibitem[APS76]{APS3}
M.~F. Atiyah, V.~K. Patodi, and I.~M. Singer.
\newblock Spectral asymmetry and riemannian geometry iii.
\newblock {\em Math. Proc. Camb. Phil. Soc.}, 1976.

\bibitem[Ber82]{Be2}
A.~J. Berrick.
\newblock The plus construction and fibration.
\newblock {\em Q.J. Math.}, 33:149--157, 1982.

\bibitem[BF78]{BF}
A.~K. Bousfield and E.~M. Friedlander.
\newblock Homotopy theory of ${\Gamma}$-spaces, spectra, and bisimplicial sets.
\newblock {\em Lecture Notes in Math.}, 1978.

\bibitem[BNV16]{BNV}
U.~Bunke, T.~Nikolaus, and M.~Völkl.
\newblock Differential cohomology theories as sheaves of spectra.
\newblock {\em J. Homotopy Relat. Struct.}, 11(1):1--66, 2016.

\bibitem[Bun14]{Bu1}
Ulrich Bunke.
\newblock A regulator for smooth manifolds and an index theorem.
\newblock \href{arXiv:1407.1379}{arXiv:1407.1379}, 2014.

\bibitem[Bun17]{Bu2}
Ulrich Bunke.
\newblock Foliated manifolds, algebraic {K}-theory, and a secondary invariant.
\newblock \href{arXiv:1507.06404}{arXiv:1507.06404}, 2017.

\bibitem[CCMT84]{CCMT}
J.~Caruso, F.R. Cohen, J.P. May, and L.R. Taylor.
\newblock {J}ames maps, {S}egal maps, and the {K}ahn-{P}riddy theorem.
\newblock {\em Trans. Amer. Math. Soc.}, 1984.

\bibitem[Don78]{Don}
Harold Donnelly.
\newblock Eta invariants for {G}-spaces.
\newblock {\em Indiana Univ. Math. J.}, 27(6):889--918, 1978.

\bibitem[Gil84a]{Gi}
P.~B. Gilkey.
\newblock The eta invariant and the {K}-theory of odd dimensional spherical
  space forms.
\newblock {\em Invent. Math.}, 1984.

\bibitem[Gil84b]{Gi2}
P.~B. Gilkey.
\newblock {\em Invariance theory, The Heat Equation, and The Atiyah-Singer
  index theorem}, volume~11 of {\em Mathematics Lecture Series}.
\newblock Publish or Perish, Inc., 1984.

\bibitem[GS99]{GS}
R.~E. Gompf and A.~I. Stipsicz.
\newblock {\em 4-manifolds and Kirby Calculus}, volume~20 of {\em Graduate
  Studies in Mathematics}.
\newblock American Mathematical Society, Providence, RI, 1999.

\bibitem[Hat]{Hat4}
Allen Hatcher.
\newblock Notes on basic 3-manifold topology.

\bibitem[Hov99]{Ho}
Mark Hovey.
\newblock {\em Model Categories}, volume~63 of {\em Mathematical Surveys and
  Monographs}.
\newblock American Mathematical Society, 1999.

\bibitem[HV78]{HV}
J.-C. Hausmann and P.~Vogel.
\newblock The plus construction and lifting maps from manifolds.
\newblock {\em Proc. Symp. Pure. Math.}, 32:67–76, 1978.

\bibitem[JW95]{JW}
J.~D.~S. Jones and B.~W. Westbury.
\newblock Algebraic {K}-theory, homology spheres and the $\eta$-invarinat.
\newblock {\em Topology}, 1995.

\bibitem[Kar87]{Ka1}
Max Karoubi.
\newblock {\em Homologie Cyclique et K-théorie}, volume 149.
\newblock Astérisque, 1987.

\bibitem[Kar90]{Ka2}
Max Karoubi.
\newblock Th\'eorie g\'en\'erale des classes caract\'eristiques secondaires.
\newblock {\em K-Theory}, 4:55--87, 1990.

\bibitem[Kar94]{Ka3}
Max Karoubi.
\newblock Classes caract\'eristiques de fibr\'es feuillet\'es, holomorphes ou
  alg\'ebriques.
\newblock {\em K-Theory}, 8:153--211, 1994.

\bibitem[Kir78]{Ki1}
Robion Kirby.
\newblock A calculus for framed links in ${S^3}$.
\newblock {\em Invent. Math.}, 45(1):35--36, 1978.

\bibitem[Kos07]{Ko}
A.~A. Kosinski.
\newblock {\em Differential Manifolds}.
\newblock Dover publications, Inc., 2007.
\newblock Original edition published by Academic press in 1993.

\bibitem[Kri11]{Kr}
Tobias Krieg.
\newblock Kamber-tondeur formen f\"ur lokal nilpotente zusammenh\"ange.
\newblock Diplomarbeit, Albert-Ludwigs-Universität Freiburg, 2011.

\bibitem[LM89]{LM}
H.~B. Lawson and M.-L. Michelsohn.
\newblock {\em Spin Geometry}.
\newblock Princeton Mathematical Series. Princeton University Press, Princeton,
  NJ, 1989.

\bibitem[LMS86]{LMS}
Jr. L.G.Lewis, J.P. May, and M.~Steinberger.
\newblock {\em Equivariant Stable Homotopy Theory}, volume 1213 of {\em Lecture
  Notes in Mathematics}.
\newblock Springer-Verlag, Berlin, 1986.
\newblock With contributions by J. E. McClure.

\bibitem[MMSS01]{MMSS}
M.~A. Mandell, J.~P. May, S.~Schwede, and B.~Shipley.
\newblock Model categories of diagram spectra.
\newblock {\em Proc. London Math. Soc.}, 82:441--512, 2001.

\bibitem[MP12]{MP}
J.~P. May and K.~Ponto.
\newblock {\em More Concise Algebraic Topology. Localization, Completion, and
  Model categories.}
\newblock Chicago Lectures in Mathematics. University of Chicago Press,
  Chicago, IL, 2012.

\bibitem[Rol90]{Rol}
Dale Rolfsen.
\newblock {\em Knots and Links}, volume~7 of {\em Mathematics Lecture Series}.
\newblock Publish or Perish, Inc., 1990.
\newblock Corrected reprint of the 1976 original.

\bibitem[Rud08]{Ru}
Y.~B. Rudyak.
\newblock {\em On Thom Spectra, Orientability, and Cobordism}.
\newblock Springer Mongraphs in Math, 2008.

\bibitem[Sav02]{Sa}
Nikolai Saveliev.
\newblock {\em Invariants for Homology 3-spheres}, volume 140 of {\em
  Encyclopaedia of Mathematical Sciences}.
\newblock Springer-Verlag, Berlin, 2002.

\bibitem[Sch]{nLPre}
Urs Schreiber.
\newblock
  \href{https://ncatlab.org/nlab/show/model+structure+on+topological+sequential+spectra}{The
  model structure on topological sequencial spectra}.

\bibitem[Sch11]{Sch}
Jan Schl\"uter.
\newblock {\em Lokal filtrierte {V}ektorb\"undel als geometrische
  {I}nterpretation der algebraischen {K}-{T}heorie}.
\newblock PhD thesis, Albert-Ludwigs-Universit\"at Freiburg, 2011.
\newblock
  \href{https://freidok.uni-freiburg.de/data/8162}{https://freidok.uni-freiburg.de/data/8162}.

\bibitem[Sus84]{Su}
A.~A. Suslin.
\newblock On the {K}-theory of local field.
\newblock {\em J. Pure Appl. Algebra}, 34:301--318, 1984.

\bibitem[Swi02]{Sw}
Robert Switer.
\newblock {\em Algebraic Topology-Homotopy and Homology}.
\newblock Classics in Mathematics. Springer-Verlag, Berlin, 2002.

\bibitem[Wan17]{Wang2}
Y.-S. Wang.
\newblock An approximation of the $e$-invariant in the stable homotopy
  category.
\newblock \href{arXiv:1707.03453 [math.KT]}{arXiv:1707.03453 [math.KT]}, 2017.

\bibitem[Wei13]{We1}
C.~A. Weibel.
\newblock {\em The K-book}.
\newblock Graduate Studies in Mathematics, 145. American Mathematical Society,
  Providence, RI, 2013.

\end{thebibliography}

\end{document}